\documentclass[12pt,a4paper]{article}

\usepackage{amsmath}
\usepackage{amssymb}
\usepackage{latexsym}
\usepackage{srcltx}
\usepackage{graphicx}
\usepackage{epsfig}
\usepackage[ruled,vlined]{algorithm2e}
\usepackage{tikz}
\usepackage{bm}
\usetikzlibrary{arrows}
\usepackage{hyperref}

\usepackage{rotating}
\usepackage{caption}
\usepackage{subcaption}
\usepackage{color}
\newcommand{\co}{{\mathbb C}}

\newcommand{\re}{{\mathbb R}}
\newcommand{\n}{{\mathbb N}}

\newcommand{\cA}{{\cal{A}}}
\newcommand{\cD}{{\cal{D}}}
\newcommand{\cO}{{\cal{O}}}
\newcommand{\cV}{{\cal{V}}}
\newcommand{\cL}{{\cal{L}}}

\newcommand{\cR}{{\cal{R}}}

\newcommand{\cC}{{\cal{C}}}
\newcommand{\cP}{{\cal {P}}}
\newcommand{\cM}{{\cal{M}}}
\newcommand{\cF}{{\cal{F}}}
\newcommand{\cK}{{\cal{K}}}

\newcommand{\hr}{{ \rho}}

\newtheorem{theorem}{Theorem}
\newtheorem{prop}{Proposition}

\newtheorem{cor}{Corollary}
\newtheorem{remark}{Remark}
\newtheorem{ex}{Example}
\newtheorem{defi}{Definition}

\date{}

\author{Antonio Cicone
\thanks{Dipartimento di Ingegneria Scienze Informatiche e
Matematica, INdAM and DEWS, L'Aquila, Italy
{e-mail: \tt\small antonio.cicone@univaq.it}}, 
Nicola Guglielmi  
\thanks{Dipartimento di Ingegneria Scienze Informatiche e
Matematica and Gran Sasso Science Institute, L'Aquila, Italy
{e-mail: \tt\small nicola.guglielmi@univaq.it}}
and Vladimir Yu.~Protasov
\thanks{Department of Mechanics and Mathematics of Moscow State University,
and Faculty of Computer Science of  National Research University Higher School of Economics,
Moscow, Russia,  {e-mail: \tt\small
v-protassov@yandex.ru}}}

\title{Linear dynamical systems on graphs
\thanks{The first author is a Marie Curie fellow of the INdAM and he acknowledges support by National Group for Scientific Computation (GNCS -- INdAM) ``Progetto giovani ricercatori 2014'', and by Istituto Nazionale di Alta Matematica (INdAM) ``INdAM Fellowships in Mathematics and/or Applications cofunded by Marie Curie Actions''. The second author acknowledges support by National Group for Scientific Computation (GNCS -- INdAM) ``Progetti di ricerca 2015''. The third author is supported by RFBR grants nos
14-01-00332 and 16-04-00832, and by the grant of Dynasty foundation}}

\begin{document}
\maketitle

\begin{abstract}
We consider linear dynamical systems with a structure of a multigraph. The vertices are associated to
linear spaces and the edges correspond to linear maps between those spaces.
We analyse the asymptotic growth of trajectories (associated to paths along the multigraph),  the
stability and the stabilizability problems. This generalizes the classical linear switching systems
and their recent extensions to Markovian systems, to systems generated by regular languages, etc. We show that an arbitrary system can be  factorized 
into several irreducible systems on strongly connected multigraphs.
For the latter systems, we prove the existence of invariant (Barabanov) multinorm and
derive a method of its construction. The method works for a vast majority of systems and finds the joint spectral radius (Lyapunov exponent). Numerical examples are presented and applications to the study of fractals,  attractors, and multistep methods for ODEs are discussed.

\smallskip

\smallskip

\begin{flushright}
\noindent  \textbf{AMS 2010} {\em subject classification 15A60, 37B25, 15-04}
\end{flushright}

\end{abstract}

\section{Introduction}\label{s-intro}

Linear switching systems draw much attention in the literature due to their applications in
the electronic engineering, dynamical systems, control theory, etc. A discrete
linear switching system (LSS) has the form
\begin{equation}\label{discr}
x(k) \ = \ A(k)\, x(k-1)\, , \quad k \in \n\, ,
\end{equation}
where $\{x(i)\}_{i =0}^{\infty}$ is a sequence of points ({\em trajectory}) in the Euclidean space $\re^d$,
$\{A(k)\}_{k\in \n}$ is a sequence of matrices ({\em switching law}) chosen independently from a given compact set
$\cA$.  The system is {\em stable} if $x(k) \to 0$ as $k \to \infty$ for every initial point $x(0) \in \re^d$
and for every switching law. It is well known that the stability is equivalent to the inequality $\rho(\cA) < 1$, where
\begin{equation}\label{jsr0}
\rho(\cA) \ =\ \lim_{k\to \infty} \, \max_{A(i) \in \cA, \, i = 1, \ldots , k}\, \bigl\|\, A(k)\cdots A(1)\, \bigr\|^{1/k}\,   \end{equation}
is the {\em joint spectral radius} (JSR) of the matrix family~$\cA$. This limit always exists and does not depend on the
matrix norm. The joint spectral radius is equal to the Lyapunov exponent of the system, which is the maximal
exponent of the asymptotic growth of trajectories: $\sup\, \frac{\log \|x_k\|}{\log k}\, = \, \rho$, where the supremum is
 over all trajectories with $\|x_0\| = 1$. If $\cA$ consists of one matrix, then JSR becomes its usual spectral radius (the largest modulus of eigenvalues of a matrix).

The joint spectral radius originated in 1960 with Rota and Strang~\cite{RS60} and
has found countless applications in various fields, from functional analysis to discrete mathematics and theory of formal languages
(see~\cite{GP13, J09} for properties and applications). The computation of JSR is an algorithmically hard problem even for finite families~$\cA$.
Nevertheless, there are several practically efficient methods for estimating the JSR~\cite{DL92, G, PJ, PJB} or even
for its precise computation~\cite{GP13, MR14} for wide classes of matrices. The theoretical base of many methods is the
so-called {\em invariant norm} called also the {\em Barabanov norm}. This is a norm in~$\re^d$ such that
$\|x\| \, = \, \rho(\cA)\, \max\limits_{A \in \cA}\|Ax\|$ for all $x \in \re^d$. Every irreducible family
of matrices (i.e., the matrices from~$\cA$ do not share a common nontrivial  invariant subspace) possesses
an invariant norm~\cite{B1}. If $\cA$ is reducible, then it can always be factorized, which makes the problem
of JSR computation to be equivalent to several analogous  problems in smaller dimensions~\cite{BW, P96}.

Recently many authors introduced and analysed {\em constrained switching systems}, where not all switching laws are possible but only those satisfying certain  stationary constraints~\cite{D, K, OPJ, PhJ1, PhJ2, SFS, WRDV}. The concept slightly varies in different papers, but in general can be described as follows: given a directed graph $G$ with edges labelled by matrices from a family~$\cA$
(one matrix may correspond to several edges). We call a trajectory  {\em admissible} if it is realized as a path along the graph. This leads to the concept of stability of constrained system
(all admissible trajectories converge to zero) and the corresponding (constrained) joint spectral radius, where
the maximum in~(\ref{jsr0}) is computed only over admissible sequences. Methods for evaluating the lower bounds and upper bounds for the constrained JSR were presented in~\cite{K, PhJ1}. We shell discuss them in more detail in Section~\ref{s-appl}.
\smallskip

The aim of this paper is to develop a method of precise computation of JSR for constrained systems by evaluating a
piecewise linear convex invariant Lyapunov function. That function is actually  a norm defined by a special convex polytope.  We first prove several theoretical results (Sections~\ref{s-constr} -- \ref{s-norms}) based on which we derive Algorithm~1 for JSR computation and for constructing invariant polytopes (Section~\ref{s-algor}).
As we shell see in examples (Section~\ref{sec:IllEx})  and in statistics of numerical experiments (Section~~\ref{sec:IllEx}) the method
 is indeed able to find  precisely the joint spectral radius for a vast majority of constrained systems. For general systems, it works efficiently in dimensions at least up to~$20$, for positive systems, it works much faster and is applicable in higher dimensions (at least up to~$100$).  For the classical (unconstrained) systems, the Invariant polytope algorithm was introduced and analyzed in~\cite{GP13}. That algorithms proved its efficiency for many problems, both numerical and theoretical. In particular, it helped solving several long-standing open problems in combinatorics, theory of formal languages, theory of wavelets and subdivisions~\cite{GP13, GP16, PG}. However,
 a naive attempt to extend  the Invariant polytope algorithm directly to constrained systems offers a strong resistance.
This extension needs  a well-developed theoretical base, which is presented in the following Sections~\ref{s-constr} -- \ref{s-norms}.

First of all, we slightly generalize the concept of constrained systems. We consider a directed multigraph $G$ with
linear spaces~$\{L_i\}_{i=1}^n$ (maybe, of different dimensions) associated to its vertices.
For each pair of vertices $g_i$ and $g_j$, there is a finite (maybe empty) collection of edges leading from
$g_i$ to $g_j$ identified with linear operators acting from $L_i$ to $L_j$. Thus we have a family of spaces
and linear operators mapping them to each other according to the multigraph~$G$. To every
path starting at some vertex $g_s$ and going along successive edges of $G$, and to every point $x_0 \in g_s$
we naturally associate a trajectory of the point~$x_0$. Thus, we obtain a linear dynamical system on the multigraph~$G$.
We use the short terminology {\em system on a graph},  although it is actually defined on a multigraph. All definitions and
properties of this construction are given in the next section. To study the stability, growth of trajectories, and JSR of such systems we realize the following plan:
\smallskip

1. First, we  prove a factorization theorem that reduces a system on an arbitrary multigraph to several smaller systems
on strongly connected multigraphs. This allows us to restrict the whole theory to the case of strongly connected
multigraph~(Section~\ref{s-connect}).
\smallskip

2. Then we introduce the concept of reducible and irreducible systems. We show that every reducible system
can be factorized to several irreducible ones of smaller dimensions~(Section~\ref{s-factor}).  Thus, we have the same situation as for the classical (unconstrained) systems. However, the notion of reducibility cannot be extended directly
to constrained systems (see an example in~\cite{PhJ2}) and requires a significant modification.
Several properties of irreducible systems are established. In particular, we show that every such a system is non-defective,
i.e., its trajectories grow not faster than $\|x_k\| \le C \rho^k\, , \, k \in \n$, where $\rho$ is the joint spectral radius and $C > 0$ is a constant (Section~\ref{s-irred}). Working with practical examples and 
with randomly generated matrices in Section~\ref{sec:IllEx} we observe a surprising phenomenon. 
Reducible systems, which are very exceptional in the usual (unconstrained) case, becomes 
usual for the constrained systems. According to our statistics given in Section~\ref{sec:IllEx}
in dimension~$d=3$ about $20 \%$ of randomly generated systems are reducible, while for 
$d = 20$, this ratio growth to $60 \%$. Thus, reducible systems  become dominant in high dimensions. 
 This makes our technique for factorizing reducible systems to be important in most practical cases.   
\smallskip

3. We restrict the theory to irreducible systems on strongly connected multigraphs.
The next step is to introduce the concept of extremal and invariant multinorms. An invariant (Barabanov) multinorm
is a collection of norms $\|\cdot \|_i$ in the spaces~$L_i, \, i=1, \ldots , n$, respectively such that
$\, \max\limits_{A_{ji}} \|A_{ji}x\|_j\, = \, \rho\, \|x\|_i$ for every $i$ and $x \in L_i$, where the maximum is over all operators associated to all outgoing  edges from $g_i$. For the extremal multinorm, the definition is the same,
but with inequality $\le$ instead of equality~(Section~\ref{s-norms}). We prove that an every irreducible system
possesses an invariant multinorm~(Theorem~\ref{th20}).
\smallskip

4. Based on these theoretical results we elaborate an algorithm that computes the joint spectral radius
and constructs a piecewise-linear extremal mutlinorm~(Section~\ref{s-algor}). A criterion of its convergence
is provided by Theorem~\ref{th60}. The efficiency, even in relative high dimensions, is demonstrated
in examples and in the statistics of numerical experiments~(Section~\ref{sec:IllEx}). Moreover, in Section~\ref{s-bar} we show how to construct Barabanov piecewise-linear norm by that algorithm.
\smallskip

This is a summary of the main results. In addition, we consider several special cases of the general construction:
Markovian systems, systems defined by regular languages, etc.~(Section~\ref{s-cases}). We derive  some corollaries, such as an improved Berger-Wang formula which sharpens the classical formula even in the usual (unconstrained) case; we estimate the rate of growth
of defective ({\em marginally unstable}) systems (Section~\ref{s-corol}). In Section~\ref{s-special} we discuss some special cases, as the one
of positive systems.
Finally, in Section~\ref{s-appl} we consider possible applications of our results to linear switching  systems, automata, fractals, attractors of
hyperbolic dynamical systems, consensus problems, and stability of multistep methods in ODEs.
\smallskip

We use the following notation. For two points $x, y \in \re^d$, we write
$x \ge y$ ($x > y$) if the vector $x-y$ is nonnegative (respectively, strictly positive).
As usual, the positive orthant $\re^d_+$ is the set of nonegative vectors.
For a given set $M \subset \re^d$ we denote by ${\rm co}(M)$ its convex hull and by
${\rm absco}(M) = {\rm co}\{M, -M\}$ the symmetrized convex hull. For $M \subset \re^d_+$ we
denote
\begin{equation}\label{pmco}
{\rm co}_{-}(M)\ = \ \Bigl\{\, x - y \ \Bigl|\  x \in {\rm co}\, (M) , \ y \ge 0\Bigr\}\ ; \quad
{\rm co}_{+}(M)\ = \ \Bigl\{\, x + y \ \Bigl|\  x \in {\rm co}\, (M) , \ y \ge 0\Bigr\}\, .
\end{equation}
Note that ${\rm co}_{+}(M)$ is always unbounded, whenever $M$ is nonempty.
If $M$ is finite then ${\rm co}_{-}(M)$ is called {\em infinite polytope}.
The set ${\rm co}_{+}(M)$ in this case is a polytope (in a usual sense) containing the
polytope ${\rm co} (M)$.
The sign $\asymp$ denotes as usual the asymptopic equivalence of two values (i.e., equivalence up to multiplication by a constant).

\medskip

\section{The general construction}\label{s-constr}

 We have a directed multigraph $G$
with $n$ vertices $g_1, \ldots , g_n$. Sometimes, the vertices will be denoted by their numbers.
To each  vertex~$i$ we associate a linear space $L_i$ of dimension $d_i < \infty$.
If the converse is not stated, we assume $d_i \ge 1$.  The set of spaces
$L_1, \ldots , L_n$ is denoted by $\cL$. For each vertices $i, j \in G$
(possibly coinciding),
there is a set $\ell_{ji}$ of edges
from $i$ to $j$. Each edge from $\ell_{ji}$ is identified with a linear operator $A_{ji}: \, L_i \to L_j$.
The family of those operators (or edges) is denoted by $\cA_{ji}$. If  $\ell_{ji} = \emptyset$, i.e., there are
no edges from $i$ to $j$, then $\cA_{ji} = \emptyset$. Thus, we have a family of spaces~$\cL$
and a family of operators-edges  $\cA = \cup_{i, j}\cA_{ji}$ that act between these spaces according to the
multigraph $G$. This {\em triplet $\, \xi = (G, \cL, \cA)$}  of the multigraph, spaces, and operators  will be called {\em system}.
A path $\alpha$ on the multigraph $G$ is a sequence of connected subsequent edges, its length (number of edges)
is denoted by $|\alpha|$. The length of the empty path is zero. To every path $\alpha$ along
vertices $i_1\rightarrow i_2 \rightarrow \cdots \rightarrow i_{k+1}$ that consists of edges (operators)
$A_{i_{s+1}i_s} \in \cA_{i_{s+1}i_s}, \, s = 1, \ldots , k$, we associate the corresponding product (composition) of operators $\Pi_{\alpha} = A_{i_{k+1}i_k}\cdots A_{i_2i_1}$.
Note that $|\alpha| = k$. Let us emphasize that a path is not a sequence of vertices but edges. If $G$ is a graph, then any path is uniquely defined by the sequence of its vertices, if
$G$ is a multigraph, then there may be many paths corresponding to the same sequence of vertices.
If the path is closed ($i_1 = i_{k+1}$), then $\Pi_{\alpha}$ maps the space
$L_{i_1}$ to itself. In this case $\Pi_{\alpha}$ is given by a square matrix, and possess eigenvalues, eigenvectors and the spectral radius $\rho(\Pi_{\alpha})$, which is the maximal modulus of its eigenvalues. The set of all
closed paths will be denoted by $\cC (G)$. For an arbitrary $\alpha \in \cC(G)$ we denote by
$\alpha^k = \alpha \ldots \alpha$ the $k$th power of $\alpha$. A closed path is called {\em simple} if it is not a power
of a shorter path.

In what follows we assume all the sets $\ell_{ji}$ and the corresponding sets of operators~$\cA_{ji}$ are finite.
This assumption is for the sake of simplicity; all our results are easily extended to the case of  arbitrary compact
sets~$\cA_{ji}$.

\begin{defi}\label{d3}
If every  space $L_i$ on the multigraph $G$ is equipped with a norm $\|\cdot \|_i$, then
the collection of norms $\|\cdot \|_i, \,
i = 1, \ldots , n$, is called  a multinorm.
The norm of an operator $A_{ji} \in \cA_{ji}$ is defined as
$\, \|A_{ji}\| \, = \sup\limits_{x \in L_i, \|x\|_i = 1}\|A_{ji}x\|_j$.
\end{defi}
Note that the notation $\|x\|_i$ assumes that $x \in L_i$.
In the sequel we suppose that our multigraph $G$ is equipped with some multinorm
$\{ \|\cdot\|_i\}_{i=1}^n$.
We denote that multinorm by $\|\cdot\|$ and sometimes use the
short notation $\|x\| = \|x\|_i$ for $x \in L_i$. Thus, we drop the index of the norm if it is clear
to which space $L_i$ the point $x$ belongs to.

For a given $x_0 \in L_i$ and for an infinite path $\alpha$ starting at the vertex $i$, we
 consider the {\em trajectory} $\{x_k\}_{k \ge 0}$ of the system along this path. Here
 $x_{k} = \Pi_{\, \alpha_{k}}\, x_0$, where $\alpha_{k}$ is a prefix of~$\alpha$ of length~$k$.

 \begin{defi}\label{d10}
The system $\xi$ is called stable if every its trajectory tends to zero as $\, k \to \infty$.
 \end{defi}
 As in the classical case of usual unconstrained discrete systems, the stability is decided in terms of the joint spectral radius, which in
 this case is modified as follows:
 \begin{defi}\label{d20}
 The joint spectral radius (JSR) of a triplet $\xi = (G, \cL, \cA)$ is
 \begin{equation}\label{jsr}
 \hr(\xi) \ = \ \lim_{k \to \infty}\, \max_{|\alpha| = k} \|\Pi_{\alpha}\|^{\, 1/k}\, .
 \end{equation}
 \end{defi}
 Note that the function $\varphi(k) = \max_{|\alpha| = k} \|\Pi_{\alpha}\|$
 possesses the property
 $$
 \varphi(k+l) \ \le \ \varphi(k)\, \varphi(l), \quad k, l \in \n\, ,
 $$
  hence by the well-known Fekete lemma~\cite{Fek} ,
 the limit $\lim_{k \to \infty}\varphi(k)^{1/k}$ exists and coincides with $\inf_{k \in \n}\varphi(k)^{1/k}$.
 This ensures that the joint spectral radius is well defined. Moreover, for every $r, k \in \n$,  we have the double inequality
 \begin{equation}\label{ineq}
\max_{\alpha \in \cC(G), \, |\alpha| = r} \, \rho\bigl(\Pi_{\alpha}\bigr)^{\, 1/r}  \ \le \ \hr(\xi) \ \le \
 \max_{|\alpha| = k} \, \bigl\|\Pi_{\alpha}\|^{\, 1/k}\, .
 \end{equation}
The right hand side of this inequality tends to $\hr(\xi)$ as $k \to \infty$, this follows from the
definition of the joint spectral radius. The upper limit of the left hand side as $k \to \infty$ is also
 equal to~$\hr(\xi)$. In the classical case ($m$ operators in one space) this fact is known as
 the {\em Berger-Wang formula}~\cite{BW}. Recently it was generalized by Dai~\cite{D} and Kozyakin~\cite{K}
 to the Markovian systems (see subsection~3.2). We extend it to general systems $\xi = (G, \cL, \cA)$
 and establish an improved version of this formula in Theorem~\ref{th40}. The right hand side of
 inequality~(\ref{ineq})  for $k=1$ implies that the joint spectral radius never exceeds the
 maximal norm of operators from~$\cA$, and this does not depend on the multinorm introduced
 for the system~$\xi$.   This leads to the alternative definition of JSR. In the classical
  case, this result is well-known and originated with Rota and Strang~\cite{RS60} and Elsner~\cite{E95}.
 It is extended to systems on graphs in a straightforward manner, we give its proof for convenience of the reader.
  \begin{prop}\label{p5}
 The joint spectral radius is the greatest lower bound of numbers $\lambda \ge 0$
 for which there exists a  multinorm $\|\cdot \| = \{\|\cdot \|_i\}_{i=1}^n$
 on~$G$ such that   $\|A_{ji}\| \le \lambda\, , \ A_{ji} \in \cA$.
 \end{prop}
 {\tt Proof. } If $\|A_{ji}\| \le  \lambda\, , \ A_{ji} \in \cA$, then from inequality~(\ref{ineq})
 for $k=1$, it follows that $\hr(\xi) \le  \lambda$. Conversely, assume $\hr(\xi) \le  \lambda$. We need to show that
for every $\mu > \lambda $, there is a multinorm such that  $\|A_{ji}\| \le  \mu\, , \ A_{ji} \in \cA$. Let $\tilde \cA = \mu^{-1}\cA$.
Clearly, all trajectories of the system $\tilde \xi = (G, \cL, \tilde \cA)$ tend to zero, because $\rho(\tilde \xi) < 1$. Hence, those trajectories  are uniformly bounded. Consequently, the function
$f(x) = \sup_{|\alpha| \ge 0}\, \|\tilde \Pi_{\, \alpha}\, x\|$ (the supremum is taken over all paths~$\alpha$ along~$G$ starting at~$x$)
is bounded for all $x \in L_i, i = 1, \ldots , n$.
 This function is positive, symmetric, and positively homogeneous. It is convex as a supremum of convex functions. Hence, $f$ is a norm.
 For each operator $\tilde A_{ji}$, every path starting at the point $A_{ji}x \in L_j$  is a part of the corresponding path
 starting at~$x$, hence, $f(x) \ge f(\tilde A_{ji}x) =  \mu^{-1} f(A_{ji}x)$. Thus,
 $f(A_{ji}x) \le \mu f(x)$ for all $x \in L_i$, hence the operator norm of $A_{ji}$ does not exceed~$\mu$.

   {\hfill $\Box$}
\medskip

The infimum in Proposition~\ref{p5} is not necessarily attained. If it is, then the corresponding norm is called {\em extremal} (Definition~\ref{d60}).
In Section~7 we are going to see that, similarly  to the classical case (unconstrained systems),
an extremal norm exists at least for irreducible systems.

An immediate consequence of Proposition~\ref{p5} is that the joint spectral radius
is responsible for the simultaneous contractibility of all the operators $A_{ij}$.
\begin{prop}\label{p6}
The following properties of a system are equivalent:

1) $\, \rho(\xi) < 1$;

2) There exists a multinorm $\|\cdot \|$ and a number $q < 1$ such that
$\|A_{ji}\| < q$ for all $A_{ji} \in \cA$.
 \end{prop}
 Thus, $\rho(\xi) < 1$ precisely when each space $L_i$ can be equipped with a norm
 $\|\cdot\|$ so that all operators $A_{ji}$ are contractions. We use this property
 in Section~\ref{s-appl} for applications to fractals and to dynamical systems.

 \noindent {\tt Proof. } If $2)$ holds, then obviously $\rho (\xi) \le q < 1$. If 1) holds, then
 taking arbitrary $q \in (\rho, 1)$ and applying Proposition~\ref{p5} for $\lambda = q$ we obtain a multinorm
 such that $\|A_{ji}\| < q$ for all $A_{ji} \in \cA$.

   {\hfill $\Box$}
\medskip

Similarly to the classical case, the joint spectral radius measures the stability of the system.
  \begin{prop}\label{p10}
 A system is stable if and only if $\hr(\xi) < 1$.
 \end{prop}
 So, for systems on graphs we have the same situation as in the classical case.
 The proof is also similar, we give it in Section~7.

\section{Special cases}\label{s-cases}

Before we establish the main properties of systems on multigraphs, we spot several important special cases.

\subsection{The classical (unconstrained) case: $m$ operators in one space.}\label{ss-classic}

If $G$ has one vertex and $m$ edges (loops) connecting that vertex with itself, we
obtain the classical case: $m$ operators $A_1, \ldots , A_m$ act in one space $\re^d$.
In this case the notions of trajectories, the joint spectral radius, invariant norms, etc., are  the same
as those elaborated in the extensive literature on asymptotics of matrix products~(see~\cite{J09, GP13} for
reviews).

\subsection{Markovian systems}\label{ss-markov}

We are given a family $\cA = \{A_1, \ldots, A_m\}$ of operators acting in the space~$\re^d$.
Let $\cD$ be a subset of the set of $m^2$ ordered pairs $\{(j, i) \ | \ j, i = 1, \ldots , m\}$.
The {\em Markovian system} consists of all {\em admissible} products of operators from $\cA$, i.e.,
 products that avoid subproducts
$A_jA_i\, , \, (j,i) \in \cD$. In other words, the operator $A_j$ cannot follow the operator $A_i$ in any product, whenever
$(j,i) \in \cD$. The {\em Markovial joint spectral radius} originated in~\cite{D}, it is defined as the usual joint spectral radius, but over a set of admissible matrix products.

The Markovian systems can be put in our framework as follows. We consider the graph $G$ with $m$
vertices, each vertex $i$ is associated to the space~$L_i = \re^d$; the edge $i\to j$ exists if and only if $(j, i) \notin \cD$, this edge corresponds to the operator $A_j$. Thus, all edges from $\cD$ are excluded from the graph, the remaining edges are arranged as follows: all incoming edges of the vertex $j$ correspond to the same operator $A_j\, , j = 1, \ldots , m$.

In Theorem~\ref{th40} we slightly improve the main results of~Dai~\cite{D} and Kozyakin~\cite{K} and extend them from the Markovian systems to
general systems.

\subsection{Identifying several vertices}\label{ss-ident}

Consider the trivial Markovian system when the set of prohibited links $\cD$ is empty.
In this case the Marovian joint spectral radius coincides with the usual joint spectral radius.
A question arises whether it is possible to treat this case without considering $m$ copies of the space
$\re^d$ as vertices of the graph $G$ (which is a clique in this case)  and to manage  with one space~$\re^d$
as in the classical case.   The answer is affirmative. This can be done by the procedure of {\em identifying
vertices that have  the same sets of outgoing edges}. In many practical cases this significantly simplifies
the analysis of general systems $\xi = (G, \cL, \cA)$.

We consider a general system $\xi$ with a multigraph $G$. If two its vertices $i_1$ and $i_2$ satisfy the following three conditions:

1. The associated spaces $L_{i_1}$ and $L_{i_2}$ have the same
    dimension;

2. They have the same set of outgoing edges, i.e., for every $j = 1, \ldots , n$, we have~${\ell_{ji_1} = \ell_{ji_2}}$;

 3. $\cA_{ji_1} = \cA_{ji_2}$, i.e., there is a basis  in the space $L_{i_1}$ and a basis in the space $L_{i_2}$
 such that

 \ \ the operators  from $\cA_{ji_1}$ are written  by the same matrices as the corresponding operators

 \ from $\cA_{ji_2}$.
\smallskip

\noindent Then the vertices $i_1$ and $i_2$ can be identified. They are replaced by one vertex $i$.
The sets~$\cA_{ji}$ of its outgoing vertices  is the same as those of vertices $i_1$ and $i_2$.
The set of incoming vertices is the union of those of $i_1$ and $i_2$.
The norm in $L_i$ can be chosen arbitrarily. For example, the pointwise maximum of norms in $L_{i_1}$ and $L_{i_2}$.
We obtain the multigraph $G'$ with $n-1$ vertices and the corresponding system~$\xi'$.

Every path $\alpha$ on $G$ is naturally identified with a path $\alpha'$ on $G'$
by replacing both $i_1$ and $i_2$ by $i$. This establishes the correspondence between trajectories of
$\xi$ and $\xi'$. In particular, those trajectories have the same asymptotics as $k \to \infty$.
Therefore, $\hr(\xi') = \hr(\xi)$.

For example, if the Markovian system has no prohibited links, then all $m$ ist vertices can be identified, and we obtain the classical system, with one space and $m$ operators.

\subsection{Maximal growth of  trajectories avoiding prohibited words}\label{ss-words}

We have a family $\cA = \{A_1, \ldots, A_m\}$ of operators acting in~$\re^d$ and a finite set
$\cD$ of words of the $m$-ary alphabet $\{1, \ldots , m\}$. This is the dictionary of prohibited words.
We consider a discrete system with operators from~$\cA$ and with trajectories avoiding those prohibited words.
In particular, we are interested in the  exponent of the maximal asymptotic growth of those trajectories.
This exponent is the joint spectral radius of the family~$\cA$ along the products
that avoid words from~$\cD$. In particular, the system is stable
if and only is this value is smaller than one.

These systems can be put in our framework as follows. Let $l\ge 2$ be the maximal length of words from~$\cD$.
The vertices of the graph $G$ are all $m$-ary words of length~$l-1$ that avoid subwords from~$\cD$
(of lengths smaller than~$l$, if they exist).
There is an edge from the  word (vertex) $\beta$ to $\gamma$ if and only if the following two conditions are satisfied:
\smallskip

1) the prefix of length $(l-2)$ of $\beta$ is the suffix of $\gamma$;
\smallskip

2) the word $\gamma_1 \, \beta = \gamma \beta_{l-1}$ of length~$l$ is not from $\cD$
($\gamma_1$ and $\beta_{l-1}$ are the first and the last letters of $\gamma$ and $\beta$ respectively).
\smallskip

If these conditions are fulfilled, then there is a unique edge $\beta\to \gamma$, it is associated to the operator $A_{\gamma_1}$. Note that if $l=2$, then condition $1)$ is always fulfilled. If $\cD = \emptyset$, i.e., there are no prohibited words,
 then the graph $G$ has $m^{l-1}$ vertices. Each vertex $\beta$ has exactly $m$ outgoing edges, to the
 vertices $1\beta, \ldots , m\beta$, and $m$ incoming edges, from the vertices $\beta1, \ldots , \beta_m$.
 Hence, for any set of prohibited words, the graph $G$ has at most $m^{l-1}$ vertices and at most $m$ incoming and $m$ outgoing edges for each vertex. So, the graph $G$ has at most $2^{l}$ edges.
 All the spaces $L_i$ at the vertices of $G$ are copies of~$\re^d$.

 To every word $c = c_1c_2\ldots c_N$ of the $m$-ary alphabet avoiding subwords from $\cD$, the corresponding path
 $\beta_1\to \beta_2\to \ldots \beta_{N}$ along~$G$ is naturally associated as follows:
 $\beta_k = c_k\ldots c_{k+l-2}\, , \, k = 1, \ldots , N+2-l$.
 This path corresponds to the product $A_{c_{N+2-l}}\cdots A_1$.

 Thus, we have the triplet~$\xi$.
The joint spectral radius of the family~$\cA$ along the products
avoiding words from $\cD$ is equal to $\hr(\xi)$.

The Markovian systems is a special case of this construction when $l=2$.

\begin{ex}\label{ex10}
{\em Let us have two operators $A_1$ and $A_2$ and one prohibited word $A_1A_2A_1$. So, we are interested in the
maximal asymptotic growth of products $A_1^{k_s}A_2^{r_s} \cdots A_1^{k_1}A_2^{r_1}$ with $r_i \ge 2$
for all $i$. In this case, $G^{(1)}$ has four vertices $A_1A_1, A_1A_2, A_2A_1, A_2A_2$ and eight edges.
The edge $A_2A_1 \to A_1A_2$ associated to the operator~$A_1$ is omitted because of the prohibited word $A_1A_2A_1$, all other seven edges are kept  (see fig \ref{fig:GraphG1}).
The JSR along infinite paths of this graph equals to the JSR along those products.

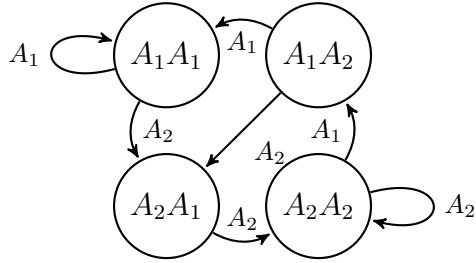
\begin{figure}[ht]
\begin{center}
\begin{tikzpicture}[->,>=stealth',shorten >=1pt,auto,node distance=2cm,
                    thick,main node/.style={circle,draw,font=\sffamily\bfseries}]

  \node[main node] (1) {$A_1 A_1$};
  \node[main node] (2) [right of=1] {$A_1A_2$};
  \node[main node] (3) [below of=1] {$A_2A_1$};
  \node[main node] (4) [below of=2] {$A_2 A_2$};

  \path[every node/.style={font=\sffamily\footnotesize}]
    (1) edge [bend right] node {$A_2$} (3)
        edge [loop left] node [left] {$A_1$} (1)
    (2) edge [bend right] node {$A_1$} (1)
        edge node {$A_2$} (3)
    (3) edge [bend right] node {$A_2$} (4)
    (4) edge [loop right] node {$A_2$} (4)
        edge [bend right] node {$A_1$} (2);

\end{tikzpicture}
\caption{ $G^{(1)}$ Graph}\label{fig:GraphG1}
\end{center}
\end{figure}

If $\cD = \{A_1A_2A_1, A_1^2\}$, then we have the set of products $A_1A_2^{r_s} \cdots A_1A_2^{r_1}$ with $r_i \ge 2$
for all~$i$. The graph $G^{(2)}$ has three vertices $A_1A_2, A_2A_1, A_2A_2$ (the vertex $A_1A_1$ has been omitted)
and four edges (one of the five edges, $A_2A_1 \to A_1A_2$, is omitted), see fig \ref{fig:GraphG2}.

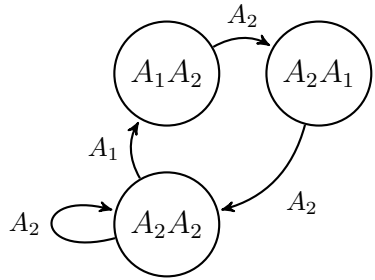
\begin{figure}[ht]
\begin{center}
\begin{tikzpicture}[->,>=stealth',shorten >=1pt,auto,node distance=2cm,
                    thick,main node/.style={circle,draw,font=\sffamily\bfseries}]

  \node[main node] (1) {$A_1 A_2$};
  \node[main node] (2) [right of=1] {$A_2 A_1$};
  \node[main node] (3) [below of=1] {$A_2 A_2$};

  \path[every node/.style={font=\sffamily\footnotesize}]
    (1) edge [bend left] node {$A_2$} (2)
    (2) edge [bend left] node {$A_2$} (3)
    (3) edge [loop left] node {$A_2$} (3)
        edge [bend left] node {$A_1$} (1);

\end{tikzpicture}
\caption{ $G^{(2)}$ Graph}\label{fig:GraphG2}
\end{center}
\end{figure}
}

 \end{ex}

\section{Strong connectivity}\label{s-connect}

First of all, let us  show that the analysis of  general
systems can be reduced to systems with strongly connected multigraphs~$G$. Recall that a multigraph is strongly connected if
for an arbitrary pair of vertices, there is a path from one to the other.

If $G$ is not strongly connected, then there is a closed submultigraph $G_1 \subset G$, for which
 all vertices reachable from $G_1$ belong to $G_1$. Let $G_2 = G\setminus G_1$
 be a complement of $G_1$ in~$G$. Denote by $\xi_i$ the restriction
 of the triplet $\xi$ to the submultigraph $G_i, i = 1,2$.
 The triplet $\xi_1$ contains only vertices from $G_1$ with the corresponding family of linear spaces~$\cL_1  \subset \cL$,  the edges of the multigraph $G$ connecting vertices from $G_1$ and the family of operators~$\cA_1 \subset \cA$ corresponding to those edges. The same with $\xi_2$.
 \begin{prop}\label{p40}
 If $G$ is  not strongly connected, then for every its closed submultigraph~$G_1\subset G$,
 we have $\hr(\xi) = \max\, \{\hr(\xi_1), \hr(\xi_2)\}$.
 \end{prop}
{\tt Proof.} If a path starts at a vertex of $G_1$, then it is contained
in $G_1$. Therefore,  any path~$\alpha$ on the multigraph $G$ is a concatenation $\alpha_1\alpha_2$, where
$\alpha_k$ is on  $G_k, \, k = 1,2$. Hence,  $\Pi_{\, \alpha} =
\Pi_{\, \alpha_1} A_{ij}\Pi_{\, \alpha_2}$, where the edge the operator $A_{ji} \in \cA_{ji}$
corresponds to an edge $l_{ji} \in \ell_{ji}$ connecting the two paths.
Let $\hr_k = \hr(\xi_k), k=1,2$, and $\bar \hr = \max\{\hr_1, \hr_2\}$.
For any $\varepsilon > 0$ we have
$\|\Pi_{\, \alpha_k}\| \le C (\hr_k + \varepsilon)^{\, |\alpha_k|}, \, k = 1,2$.
Since $|\alpha| = |\alpha_1| + |\alpha_2|-1$ and $\|A_{ji}\| \le C$, we have
 $\|\Pi_{\, \alpha}\| \le C^3 (\bar \hr + \varepsilon)^{\, |\alpha|-1}$. Hence $\hr \le \bar \hr + \varepsilon$
 for every $\varepsilon > 0$, and so $\hr \le \bar \hr$. On the other hand, obviously $\hr \ge \hr_k$, for each
 $k = 1,2$,  and therefore $\hr \ge \bar \hr$.

   {\hfill $\Box$}
\medskip

An elementary induction shows that for a not strongly connected multigraph,
there is a disjoint  partition of its vertices $G = \bigsqcup_{\, i=1}^{\, r} G_i$ such that
each submultigraph $G_i$ is strongly connected and $G_i$ is not reachable from $G_j$, whenever $i > j$.
Proposition~\ref{p40} yields
\begin{cor}\label{c5}
If $G$ is not strongly connected then $\hr(\xi) = \max \, \{\hr(\xi_1), \ldots , \hr(\xi_r)\}$.
\end{cor}
Thus, the analysis of the asymptotic properties of an arbitrary  system $\xi = (G, \cL, \cA)$ is reduced
to the same problem for the
systems $\xi_i$ that are characterized by strongly connected multigraphs
Our next assumption
concerns irreducibility of the system.

\section{Irreducibility}\label{s-irred}

 Consider an arbitrary triplet $\xi = (G, \cL, \cA)$. Here the case of trivial spaces $L_i = \{0\}$ is allowed for some (but not all) $i$.
 A triplet
$\xi' = (G, \cL', \cA')$ is {\em embedded} in $\xi$, if $L_i' \subset L_i$
for each $i$ and every operator $A_{ji}' = A_{ji}|_{L_i'}$ maps $L_i'$ to $L_j'$, whenever
$l_{ji} \in G$.
The embedding is strict if $L_i'$ is a proper subspace of $L_i$ at least for one~$i$. Thus,
an embedded triplet has the same multigraph and smaller spaces at the vertices.
\begin{remark}\label{r10}
{\em Actually, one could extend this definition allowing the embedded triplet $\xi'$ to omit some edges from~$G$,
i.e., to have a submultigraph $G' \subset G$ instead of the same multigraph~$G$. However, instead of eliminating an edge $l_{ij} \in \ell_{ji}$ we can set the corresponding $A_{ji}'$
to be the null operator. In this case, all trajectories passing through the eliminated edges vanish and can therefore be ignored. So, for the sake of simplicity, we always assume that the embedded triplet $\xi'$ has the same
multigraph $G$.}
\end{remark}

\begin{defi}\label{d40}
A triplet $\xi = (G, \cL, \cA)$ is reducible if it has a strictly embedded triplet. Otherwise, it is called irreducible.
\end{defi}
There are three remarkable properties of irreducible triples:
\smallskip

1) for an arbitrary initial vector $x \in L_j$,
its images $\Pi_{\alpha}\, x$ span all the spaces $L_1, \ldots , L_n$;
\smallskip

2) $\hr (\xi) > 0$, provided all the spaces $L_i$ are nontrivial;
\smallskip

3) if $\hr(\xi) = 1$, then all trajectories are uniformly bounded;
\smallskip

4) there is an invariant (Barabanov) multinorm.
\smallskip

\noindent We prove these properties and then, in Theorem~\ref{th30}, show that a general triplet
can be factorized to several irreducible ones of smaller dimensions. This will imply that the whole analysis can be focused on irreducible triples. Properties 1), 2), and  3) are established in this section, the proof of 4) is postponed to Section~\ref{s-norms}.

We begin with Property 1), which is characteristic of an irreducible family: the orbits of any nonzero element
$x \in L_i$ span all the spaces $L_1, \ldots , L_n$. Consider the set  of all paths $\alpha$ on $G$
from a vertex $i$ to $j$. For arbitrary $x \in L_i$ we denote by $\cO_j(x) = \{\Pi_{\alpha} x\ | \
\alpha\, : \,  i \to j\}$  the set of points from $L_j$ to which all the corresponding products $\Pi_{\alpha}$
map the point $x$. Thus, $\cO_j(x)$ is the complete {\em orbit} of the point $x \in L_i$ in the space $L_j$.
\begin{prop}\label{p50}
 A triplet $\xi = (G, \cL, \cA)$ is irreducible if and only if for every $i$ and for an arbitrary point $x \in L_i, x\ne 0$,
its orbits $\cO_j(x), \, j = 1, \ldots , n$,  are all full-dimensional, i.e., each orbit $\cO_j(x)$ spans the
corresponding space $L_j$.
\end{prop}
{\tt Proof.} Let $L_j' = {\rm span}\, \cO_j(x)$. We have $A_{kj}L_j' \subset L_k'$ for all pairs $j, k \in \{1, \ldots , n\}$
and for all $A_{kj} \in \cA_{kj}$.
Besides, the space $L_i'$ is nontrivial, since it contains $x$. Hence, $\xi' = (G, \cL', \cA|_{\cL'})$ is an embedded triplet.
By the irreducibility, it cannot be strictly embedded, consequently $\, L_j' = L_j$ for all $j$.

   {\hfill $\Box$}
\medskip

Now we are ready to establish 2).

\begin{prop}\label{p55}
 If a triplet $\xi = (G, \cL, \cA)$ is irreducible and all the spaces $L_i$ are nontrivial,
 then $\hr(\xi) > 0$.
\end{prop}
{\tt Proof.} By Proposition~\ref{p50}, for an arbitrary nonzero element $x \in L_1$,
there are products $\Pi_{i}: L_1 \to L_1$ such that the vectors $\{\Pi_ix\}_{i=1}^{{\rm dim}\, L_1}$ span
$x$, i.e, $x = \sum_i \alpha_i\Pi_ix$. On the other hand, Proposition~\ref{p5} yields that for any $\varepsilon > 0$,
there is a multinorm $\|\cdot \|$ such that $\|\Pi_i\| < \varepsilon$ for all $i$, and hence
$\|x\| < \varepsilon \sum_i |\alpha_i|$. Taking $\varepsilon$ small enough, we come to the contradiction.

   {\hfill $\Box$}
\medskip

\begin{defi}\label{d50}
A triplet $\xi$ is {\em non-defective} if there is a constant $C$ such that
$\|\Pi_{\, \alpha}\| \, \le \,  C\, \rho^{\, |\alpha|} $ for all paths $\alpha$ along $G$.
\end{defi}
If all the spaces $L_i$ are nontrivial, then $\hr  > 0$ (Proposition~\ref{p55}). After normalization, when $\hr = 1$, the non-defectivity means that all admissible products of operators from $\cA$
are uniformly bounded. The crucial fact about the non-defectivity of irreducible systems is well-known for the classical case (unconstrained systems). We are now extending it to arbitrary triples.
\begin{theorem}\label{th10}
An irreducible triplet with nontrivial spaces is non-defective.
\end{theorem}
{\tt Proof.} After normalization it can be assumed that $\hr = 1$.
For  arbitrary $i$, we consider the unit sphere
$S_i = \{x \in L_i \ | \ \|x\| = 1\}$ in the space~$L_i$. For each $k \in \n$, we
denote by $U_{i,k}$ the set of points $x \in L_i$ such that there exists
a path $\alpha$ of length $k$ starting at the vertex $i$ such that
$\|\Pi_{\, \alpha} x\| > 2$. Note that all those sets are open in $S_i$. Consider two possible cases.

If $\cup_{k \in \n} U_{i,k} \, = \, S_i$ for each $i = 1, \ldots , n$, then by the
compactness of the spheres~$S_i$, the open cover $\cup_{k \in \n} U_{i,k}$
admits a finite subcover $\cup_{k \le N_i} U_{i,k} = S_i$.
This implies that for every $i$ and for every $x \in S_i$,
there is a path $\alpha = \alpha (x)$ of length at most
$N_i$ such that $\|\Pi_{\, \alpha} x\| > 2$.
Let $N = \max\limits_{i=1, \ldots , n}N_i$.
Then, for every
$x \in L_i, x\ne 0$,  there is a path $\alpha = \alpha (x), \, |\alpha| \le N$, such that $\|\Pi_{\, \alpha} x\| >  2\, \|x\|$. Starting with arbitrary
$x_1 \in S_1$ we successively build a sequence $\{x_k\}_{k \in \n}$
such that for every $j \in \n$,  we have $\|x_{j+1}\| > 2\|x_j\|$ and there is a path of length at most
$N$ from $x_j$ to $x_{j+1}$. Therefore, $\|x_{k+1}\| > 2^{\, k}$ and $x_{k+1}$ is obtained from $x_1$
by multiplying with a product $\Pi_k$ of length at most $Nk$. Hence, $\hr \ge \lim\limits_{k \to \infty}\|\Pi_k\|^{1/kN} > 2^{\, 1/N}$, which contradicts to the assumption $\hr = 1$.

Otherwise, if for some $q \in \{1, \ldots, n\}$, the union $\cup_{k \in \n} U_{q,k}$ does not cover the sphere
$S_q$, then there exists $z \in S_q$ which does not belong to any of the sets $U_{q, k}$.
This means that for every path $\alpha$ staring at the vertex $q$ we have  $\|\Pi_{\, \alpha} \, x\| \le 2$.
Hence, the point $z$ has a bounded orbit. For every $i$, denote by $M_i$ the set of points from
$L_i$ that have bounded orbits. Observe several properties of the sets $\cM = \{M_i\}_{i = 1, \ldots , n}$.
\smallskip

1. Every $M_i$ is a linear subspace of $L_i$. Indeed, if $x, y \in M_i$,
then for any linear combination $ax + by, \, a, b \in \re$, and for every product $\Pi$, we have
$\|\Pi (ax + by)\| \le |a| \|\Pi x\| + |b| \|\Pi y\|$. Hence, if $x$ and $y$ have bounded orbits, then so
does $ax + by$.
\smallskip

2. If the set of edges $\ell_{ji}$ is nonempty, then $A_{ji}M_i \subset M_j$ for every $A_{ji} \subset \cA_{ji}$.
Indeed, if $x\in M_i$, then norms of elements of all trajectories starting at $x$ are uniformly bounded by some constant.
Hence, all trajectories starting at the point $A_{ji}x$ are also bounded by the same constant. Consequently,
$A_{ji}x \in M_j$.
\smallskip

3. The subspace $M_q$ is nontrivial, since it contains~$z$.

\noindent

Thus, we have a nontrivial triplet $(G, \cM, \cA|_{\cM})$ embedded into $(G, \cL, \cA)$.
If $M_j \ne  L_j$ at least for one $j$, then the triplet $(G, \cL, \cA)$ is reducible, which is impossible.
Otherwise, if $M_i = L_i$ for all $i$, then all points from these spaces have bounded orbits.
Take some $i$ and consider an orthonormal basis $e_1, \ldots , e_{d_i}$ of the space $L_i$.
If the orbit of each $e_s$ is bounded by a constant $C$, then for for every
$x = \sum_s x^{s}e_s \in S_i$ we have $\sum_s |x^{s}| \le \sqrt{d_i}$, and hence the orbit of
$x$ is bounded by the constant $C \sqrt{d_i}$. Hence, for all paths $\alpha$ starting at the vertex $i$,
we have $\|\Pi_{\, \alpha}\| \le C \sqrt{d_i}$. For $d = \max_i d_i$, we obtain that
all admissible products of operators from $\cA$ are bounded by norm by $C \sqrt{d}$.

   {\hfill $\Box$}

\section{Factorization of reducible systems}\label{s-factor}

In the classical case (unconstrained systems), the main advantage of using irreducible systems is that the general case can be solved by several irreducible ones of smaller dimensions.
This is done by a simultaneous factorization of all matrices of the family to an upper block-diagonal form. That is why the irreducible case can be considered as the basic one.
What situation do we have for the general triplet~$\, \xi = (G, \cL, \cA)$ ? We are going to see that
a reducible triplet can always be factorized, which splits the problem into several irreducible problems of smaller dimensions.

Let a triplet $(G, \cL, \cA)$ be reducible and have an embedded triplet
$(G, \cL^{(1)}, \cA^{(1)})$.  For every~$i$ we take an arbitrary subspace $L^{(2)}_i$ that
complements  $L^{(1)}_i$ to $L_i$.
 Denote ${\rm dim}\, L^{(s)}_i = d^{(s)}_i, \, s = 1,2$.
Thus, $L_i = L^{(1)}_i \oplus L^{(2)}_i\, $ and $\, d^{(1)}_i + d^{(2)}_i = d_i$, $\, i = 1, \ldots , n$.
Consider arbitrary vertices $i, j$ and an operator $A_{ji}$ from the family $\cA_{ji}$.
Let $A_{ji}^{(1)}= A_{ji}|_{L^{(1)}_i}$ be the operator from $L^{(1)}_i$ to $L^{(1)}_j$ and let
$A_{ji}^{(2)} = P_j\circ A_{ji}|_{L^{(2)}_i}$ be the
operator  from $L^{(2)}_i$ to $L^{(2)}_j$, which is the composition of $A_{ji}|_{L^{(2)}_i}$
and of the operator $P_j$ of projection of the space $L_j$ to its subspace $L^{(2)}_j$ parallel to
the subspace $L^{(1)}_j$. Each vector $x \in L_j$ has a unique representation
$x = x^{(1)} + x^{(2)}$ with $x^{(s)} \in L_j^{(s)}, \, s = 0,1$. The projection $P_j$ is defined as
$P_jx = x^{(2)}$.

 For each $j$, we take arbitrary bases of subspaces
$L^{(1)}_j$ and $L^{(2)}_j$, their union is a basis for~$L_j$. In this basis,
the projection $P_j$ is written by $d_j\times d_j$-matrix
\begin{equation}\label{Pj}
P_{j} \  = \
\left(
\begin{array}{cc}
0 & 0 \\
0 & I_{d_j^{(2)}}\,
\end{array}
\right)
\end{equation}
with two square diagonal blocks: zero matrix of size $d_j^{1}$ and the identity matrix~$I_{d_j^{(2)}}$ of size~$d_j^{(2)}$. In the same basis,
each matrix $A_{ji} \in \cA_{ji}$ has the following block upper-triangular form:
\begin{equation}\label{upper2}
A_{ji} \  = \
\left(
\begin{array}{cc}
A_{ji}^{(1)} & D_{ji} \\
0 & A_{ji}^{(2)}\, .
\end{array}
\right)
\end{equation}
Here $A_{ji}^{(1)}$ is a $d_j^{(1)}\times d_i^{(1)}$-matrix, $A_{ji}^{(2)}$ is a $d_j^{(2)}\times d_i^{(2)}$-matrix.
It is easy to see that for any path $i_1\to \ldots \to i_{k}$ along the multigraph $G$,
any product $A_{i_{k}i_{k-1}}\cdots A_{i_{2}i_{1}}$,
for arbitrary choice $A_{i_{s+1}i_s} \in \cA_{i_{s+1}i_s}$ for each $s=1, \ldots , k-1$,   also has block upper-triangular form~(\ref{upper2})
with the diagonal blocks of sizes $d_{i_{k}}^{(1)}\times d_{i_1}^{(1)}$ and
$d_{i_k}^{(2)}\times d_{i_1}^{(2)}$ respectively.
We denote $\xi^{(1)} = (G, \cL^{(1)}, \cA^{(1)}), \, \xi^{(2)} = (G, \cL^{(2)}, \cA^{(2)})$
\begin{theorem}\label{th30}
Every reducible triplet~$\xi$ can be factorized into two triples
$\xi^{(1)}$ and $\xi^{(2)}$ of smaller total dimensions. All matrices $A_{ji}$
of the family $\cA$ are factorized in the form~(\ref{upper2}). For the joint spectral radii, we have
\begin{equation}\label{reduc2}
\hr(\xi) \ =\ \max \, \bigl\{ \, \hr(\xi^{(1)})\, , \, \hr(\xi^{(2)})\,  \bigr\}\, .
\end{equation}
\end{theorem}
{\tt Proof.} We have proved all claims of the theorem except for the equality~(\ref{reduc2}).
Denote $\hr_i = \hr(\xi^{(i)}), i=1,2$, and $\bar \hr = \max \{\hr_1, \hr_2\}$.
Observe that in the $L_1$-norm, all products of matrices from $\cA$ are bigger than
the corresponding products to their submatrices from $\cA^{(1)}$. Consequently, $\hr\ge \hr_1$.
Similarly, $\hr \ge \hr_2$, and hence $\hr \ge \max \{\hr_1, \hr_2\} = \bar \hr$.
To establish the inverse inequality we take
an arbitrary path $\alpha = i_1\to \ldots \to i_{k+1}$. For each $s = 1, \ldots , k+1$
 we denote $\alpha_s^{-} = i_1 \to \cdots \to i_s$ (if $s=1$, the path is empty) and
 $\alpha_s^{+} = i_{s+1} \to \cdots \to i_{k+1}$ (if $s=k$, the path is empty).
 Take arbitrary $\varepsilon > 0$. The product
 $\Pi_{\, \alpha}$ has the same upper triangular block form~(\ref{upper2}). In the upper
 diagonal block it has the product $\Pi_{\, \alpha}^{(1)}$ whose norm does not exceed $C (\rho_1+\varepsilon)^{k}$.
 In the lower  diagonal block it has the product $\Pi_{\, \alpha}^{(2)}$ whose norm does not exceed $C (\rho_2+\varepsilon)^k$. Both these values do not exceed $C (\bar \rho +\varepsilon)^k$. Finally,
 the off-diagonal block is equal to
 \begin{equation}\label{marg-sum}
 \sum_{s=1}^{k}\ \Pi_{\, \alpha^{+}_{s}}^{(1)}\ D_{i_{s+1}i_{s}}\ \Pi_{\, \alpha^{-}_{s}}^{(2)}\, .
 \end{equation}
 the norm of the $s$th  term is bounded below by $C(\rho_1+ \varepsilon)^{s-1} \|D_{i_{s+1}i_{s}}\|
 \, C(\rho_2+ \varepsilon)^{k-s}$. Estimating both $\rho_1$ and $\rho_2$ from above by $\bar \rho$ and
 all $\|D_{i_{s+1}i_{s}}\|$ by $C$, we obtain the upper bound $C^3 (\bar \rho+ \varepsilon)^{k-1}$.
Hence, the norm of the the off-diagonal block in the product $\Pi_{\, \alpha}$
does not exceed $C^3 k\, (\bar \rho+ \varepsilon)^{k-1}$. Thus,
for every path $\alpha$ of length $k$, we have $\|\Pi_{\, \alpha}\| \le C_0 k (\bar \rho+ \varepsilon)^{k}$,
where $C_0$ does not depend on $\alpha$. Taking the power $1/k$ and a limit as $k \to \infty$, we see that
$\rho \le \bar \rho + \varepsilon$. Since this holds for every $\varepsilon$, we have $\rho \le \bar \rho$.

  {\hfill $\Box$}
\medskip

If we have a reducible triplet $\xi$, then applying Theorem~\ref{th30} several times, we obtain
\begin{cor}\label{c10}
Every reducible triplet~$\xi$ can be  factorized as a sum of $r\ge 2$ irreducible triples
$\xi^{(1)}, \ldots , \xi^{(r)}$ of smaller total dimensions. All matrices $A_{ji}$
of the family $\cA$ are factorized in the form
 \begin{equation}\label{blocks}
A_{ji}  \quad = \quad \left(
\begin{array}{cccccc}
A^{(1)}_{ji} & * &  \ldots &  * \\
0 & A^{(2)}_{ji}&  * & \vdots \\
\vdots & {} &  \ddots  & * \\
0 &  \ldots &  0 & A^{(r)}_{ji}
\end{array}
\right)\ ,
\end{equation}
where the matrix $A^{(s)}_{ji}$ in the $s$th diagonal block represents the family
$\cA^{(s)}$ of the irreducible triplet $\xi^{(s)} = (G, \cL^{(s)}, \cA^{(s)})$.
 For the joint spectral radii, we have
\begin{equation}\label{reducr}
\hr(\xi) \ =\ \max \, \bigl\{ \, \hr(\xi^{(1)}), \, \ldots ,  \, \hr(\xi^{(r)})\,  \bigr\}\, .
\end{equation}
\end{cor}

\begin{remark}\label{r12}
{\em Another concept of irreducibility of triplets was suggested in~\cite{PhJ2}, where it was shown that it also sufficient for non-defectivity. Definition from~\cite{PhJ2} involves the set of all cycles of~$G$, which is finite, but may be very large. This made it possible to prove theoretical decidability of irreducible systems, although its practical use is difficult for some graphs. Our concept has an advantage that it allows us to factorize an arbitrary system to several irreducible ones, exactly as in the classical case of usual (unconstrained) systems.
This extends most of  methods from irreducible systems to all systems. In particular, the problem of computing or estimating the joint spectral radius is completely reduced to that case by means of formula~(\ref{reducr})}.
\end{remark}

\begin{remark}\label{r15}
{\em Note that the problem of deciding irreducibility is algorithmically hard even for usual (unconstrained) systems of two matrices (see~\cite{AP} and references therein). Nevertheless, for our method of JSR computation (Section~\ref{s-algor}) this problem usually does not offer any resistance. If the system is defective (Definition~\ref{d50}) and the Invariant polytope algorithm does not converge, then we can make a step-by-step construction of an embedded system, which reduces the problem
to two similar problem of smaller total dimensions (Remark~\ref{r100}).}
\end{remark}

\section{Extremal multinorms and invariant multinorms}\label{s-norms}

The next crucial property of irreducible triplets is the existence of extremal and invariant multinorms.
Again, for the classical case (with one space and an invariant norm instead of multinorm)
this fact is well-known, it originated with Barabanov in~\cite{B1} a dual fact was independently proved in~\cite{P96}.
We are going to extend Barabanov's theorem for
all triplets $\xi = (G, \cL, \cA)$.

\begin{defi}\label{d60}
A multinorm $\|\cdot \| = \{\|\cdot \|_i\}_{i=1}^n$ is extremal if for every $i$ and $x \in L_i$, we have
\begin{equation}\label{extr}
\max_{A_{ji} \in \cA_{ji}, \, j=1, \ldots , n} \, \|A_{ji}x\|_j \ \le \ \hr\, \|x\|_i\,  .
\end{equation}
A multinorm is called invariant, or Barabanov, if   for every $i=1, \ldots , n$ and $x \in L_i$, we have
\begin{equation}\label{invar}
\max_{A_{ji} \in \cA_{ji}, \, j=1, \ldots , n} \, \|A_{ji}x\|_j  \ = \ \hr\, \|x\|_i\, .
\end{equation}
\end{defi}
An invariant multinorm is also extremal. So, it will suffice to prove the existence results
for invariant multinorms. On the other hand, extremal multinorms are sufficient to compute the JSR.
The class of extremal multinorms is much wider and they are easier to find or to estimate in practice.

The invariance property of a Barabanov multinorm remains valid after multiplication of all operators $A_{ji} \in \cA$
by the same constant.
Hence, it can always be assumed than  our system is normalized
so that $\hr=1$.
The multinorm
$\|\cdot \|$ is invariant if
for every point $x \in L_i$, the maximal norm of its images $\|A_{ji}x\|_j$ over all edges
going from the vertex $i$ is equal to
$\|x\|_i$. Let $B_i$ and $S_i$ be the unit ball and the unit sphere of the invariant norm in $L_i$.
Then for every $i$ and $x \in S_i$,  all images $A_{ji}x$ lie inside the corresponding
balls $B_j$ and at least one of them lies  on the sphere $S_j$.
In what follows we work with a multinorm $\|\cdot\|$ and drop the index $i$ of each concrete
norm $\|\cdot\|_i$ of the space $L_i$
(see the remark after Definition~\ref{d3}).

\begin{theorem}\label{th20}
An irreducible triplet possesses an invariant multinorm.
\end{theorem}
{\tt Proof.} First, we omit all vertices with zero-dimensional spaces $L_i$, along with
all their incoming and outgoing vertices. This does neither change irreducibility nor the JSR.
 Thus, we assume $d_i \ge 1$ for all $i$. By Proposition~\ref{p55}, $\hr(\xi) > 0$, hence,
after normalization it can be assumed that $\rho = 1$.
For every $i$ and $x \in L_i$, we denote
$f(x) = \limsup\limits_{|\alpha| \to \infty} \|\Pi_{\, \alpha}\, x\|$. By Theorem~\ref{th10}, the function
$f(x)$ is bounded. It is convex being a pointwise  upper limit of convex functions. Obviously, $f$ is symmetric and positively homogeneous. Furthermore, it possesses the invariance property:
 $\max\limits_{l_{ji} \in \ell_{ji}, \, j=1, \ldots , n} f(A_{ji}x)  \ = \ f(x), \,
x \in L_i, \, i=1, \ldots , n$.
It remains to show that $f$ is a norm, i.e., that $f(x) > 0$ for all $x\ne 0$.
Note that if an upper limit of a nonnegative sequence is zero, then that sequence tends to zero.
Thus, $f(x) = 0$ implies $\lim\limits_{|\alpha| \to \infty} \|\Pi_{\, \alpha}\, x\| = 0$.
Let $M_i$ be the set of points $x \in L_i$ satisfying this equality.
It is shown easily that $M_i$ is a linear subspace of $L_i$ and moreover, $A_{ji}M_i \subset M_j$, for all $i, j$,
and all $A_{ji} \in \cA_{ji}$. Hence the irreducibility yields that either $M_i = L_i$ for all $i$ or
$M_i = \{0\}$ for all $i$. In the former case, for all $x$, we have $\Pi_{\, \alpha}\, x \to  0$
as $|\alpha| \to \infty$.
Take some~$i$ and apply this assertion to all elements of an orthonormal basis $e_1, \ldots , e_{d_i}$ of the space $L_i$. If $\|\Pi_{\, \alpha}\, e_s\| < \varepsilon $ for all $s = 1, \ldots , d_i$, whenever $|\alpha| > N$,
then $\|\Pi_{\, \alpha}\| < d_i\varepsilon$, whenever $|\alpha| > N$.
Thus, $\|\Pi_{\, \alpha}\| \to  0$
as $|\alpha| \to \infty$, which contradicts to the assumption $\rho = 1$. Finally, in the latter case, when
$M_i = \{0\}$ for all $i$, we have $f(x) > 0$ for all $x\ne 0$, i.e., $f$ is a norm. Due to its invariance property, it is a desired invariant norm.

  {\hfill $\Box$}
\medskip

\begin{remark}\label{r50}
{\em If a system $\xi = (G, \cL, \cA)$ has an invariant multinorm, then
after possibly identifying two vertices $i_1, i_2$ of the multigraph $G$ (see subsection 3.3)
we set the norm in the new space $L_i$ to be the maximum of ``old''
norms in $L_{i_1}$ and $L_{i_2}$:  $\, \|x\|_{i} = \max \bigl\{\|x\|_{i_1}\, , \, \|x\|_{i_2} \bigr\}$. It is checked easily that the
new multinorm is also invariant.
}
\end{remark}
We are now able to prove Proposition~\ref{p10} on the stability of  linear systems.
\smallskip

{\tt Proof of Proposition~\ref{p10}.} The sufficiency is easy. If $\rho(\xi) < 1$, then for an arbitrary $q > \rho (\xi)$ and for an arbitrary multinorm, there exists a constant $C$ such that $\|P_{\alpha}\| \le C \, q^{|\alpha|}$ for all paths along~$G$. Fix some $q < 1$.
For every trajectory $\{x_k\}_{k \ge 0}$,  we have $\|x_k\| \le C q^|\alpha|\, \|x_0\| \to 0$ as $k \to \infty$, hence the system is stable.

To establish the necessity it suffices to show that if $\rho(\xi) = 1$, then there is at least one
trajectory that does not converge to zero. One can show even more: there is a trajectory such that $\|x_k\| \ge C, \, k \in \n$, where $C > 0$ is some constant. In view of Corollary~1, the multigraph~$G$ has a strongly connected sub-multigraph~$G'$
with $\rho(\xi') = 1$. So, we can consider only trajectories along~$G'$. Hence, without loss of generality we
assume that~$G$ is strongly connected. Furthermore, if $\xi$ is reducible, then it can be factorized to the form~(\ref{blocks})
with irreducible blocks~$\cA_{ji}^{(s)}, \, s = 1, \ldots , r$, and at least one of these blocks
has joint spectral radius one (Corollary~\ref{c10}. If this block has a trajectory bounded away from zero, then the corresponding
trajectory of the original system possesses the same property. Hence, it suffices to realize the proof for that irreducible block. This allows us to assume without loss of generality that $\xi$ is irreducible, in which case   Theorem~\ref{th20}
provides an invariant multinorm $\|\cdot\|$. Taking arbitrary $x_0 \in L_{i_0}, \, \|x_0\| = 1$, and applying recursively the invariance property we obtain a trajectory $\{x_k\}_{k \ge 0}$ such that $x_k \in L_{i_k}$ and $\|x_{k+1}\| =
\max\limits_{A_{i_{k+1}i_{k}} \in A_{i_{k+1}i_{k}}} \|A_{i_{k+1}i_{k}}x_{k}\| = \|x_k\| = 1$
for all~$k$ (by the strong connectivity, the maximum is taken over a nonempty set). Thus, $\|x_k\| = 1$
for all $k$ which completes the proof.

  {\hfill $\Box$}

\section{Two corollaries}\label{s-corol}

Before we turn into the algorithmic part we observe two corollaries of Theorem~\ref{th20}.
The first one is Theorem~\ref{th40} below that slightly improves the Berger-Wang formula
that expresses the joint spectral radius by the spectral radii of products.
To the best of our knowledge, this improved version is new even for the classical (unconstrained) case.
The second corollary concerns the issue of marginal instability of a system, i.e., possible
growth of trajectories in case $\rho(\xi) = 1$. Here we generalize some results that are known in
the classical case.

\subsection{An improved Berger-Wang formula}\label{ss-bw}

Applying Theorems~\ref{th20} and~\ref{th30} we can establish the following improved
Berger-Wang formula that sharpens the results of Dai~\cite{D} and Kozyakin~\cite{K}.
\begin{theorem}\label{th40}
For any triplet we have
\begin{equation}\label{BerW}
\limsup\limits_{k \to \infty} \quad (\rho(\xi))^{-k}\, \max\limits_{\alpha \in \cC(G), |\alpha| = k}\, \rho(\Pi_{\, \alpha})\, = \, 1\, .
\end{equation}
\end{theorem}
Clearly, this assertion is stronger than the Berger-Wang formula. If, for instance,
$\rho(\xi) = 1$, then the upper limit of
$\max\limits_{\alpha \in \cC(G), |\alpha| = k}\rho (\Pi_{\, \alpha})$ as $\, |\alpha| \to \infty\, $ is one. This, of course, implies
that the upper limit of $\max\limits_{\alpha \in \cC(G), |\alpha| = k}\bigl( \rho (\Pi_{\, \alpha})\bigr)^{1/k}$
is one (the Berger-Wang formula), but not vice versa.

Let us emphasize that assertion~(\ref{BerW}) holds for all triplets, including reducible ones and those not
strongly connected. An analogous statement for the norm $\|\Pi_{\, \alpha}\|$ instead
of the spectral radius $\rho(\Pi_{\, \alpha})$ holds only for special norms (i.e., for extremal norms). For
reducible triplets the upper limit of $(\rho(\xi))^{-k}\|\Pi_{\, \alpha}\|$ may be infinite, in which case
there is no norm possessing property~(\ref{BerW}).
\smallskip

{\tt Proof.} After normalization it can  be assumed that  $\rho(\xi) = 1$. Consider first the case when
$\xi$ is irreducible. By Theorem~\ref{th20} there exists an invariant norm. Hence, there
are infinite trajectories $\{x_k\}_{k \in \n}$ such that $\|x_k\| = 1$ for all $k$.
Infinitely many  points $x_k$ belong to one space $L_i$ and, due to compactness of the unit sphere,
there is a convergent  subsequence $x_{k_s}$ as $s \to \infty$. Fix some $\varepsilon > 0$.
For arbitrary $\delta > 0$,  there is $N = N(\delta)$ such that
$\|x_{k_{s}} - x_{k_{s+1}}\| < \delta$, whenever $s > N$.
On the other hand, since the points $x_{k_{s}}$ and  $x_{k_{s+1}}$
belong to one trajectory, it follows that there is a product $\Pi$ such that
$x_{k_{s+1}} = \Pi x_{k_s}$ (see~\cite[Lemma 2]{E95}). Thus, $\|(\Pi - I)x_{k_s}\| < \delta\, , \, \| x_{k_s}\| = 1$ and
$\|\Pi\| \le 1$. This implies that $\rho(\Pi) > 1-\varepsilon$, whenever $\delta$ is small enough.
 Thus, there are closed paths of the multigraph $G$ such that the spectral radii of the corresponding products
 are arbitrarily close to one. This proves~(\ref{BerW}) for irreducible triples.

If $\xi$ is reducible, then by Corollary~\ref{c10},
$\xi$ is a sum of irreducible triplets $\xi^{(1)}, \ldots , \xi^{(r)}$, and $\xi^{(s)} = 1$ for some of them.
As shown above, assertion~(\ref{BerW}) holds for the triplet $\xi^{(s)}$.
This means that for any $\varepsilon > 0$,
 there are arbitrarily long closed paths for which $\rho(\Pi^{(s)}_{\, \alpha}) > 1-\varepsilon$.
 On the other hand, in factorization~(\ref{blocks}) we have $\rho(\Pi_{\, \alpha}) =
\max\limits_{j=1, \ldots , r} \rho(\Pi^{(j)}_{\, \alpha}) \ge \rho(\Pi^{(s)}_{\, \alpha})$.
Thus, $\rho(\Pi_{\, \alpha}) > 1-\varepsilon$, which concludes the proof.

  {\hfill $\Box$}

\subsection{Marginal instability}\label{ss-marg}

A system is called {\em marginally stable} if the corresponding normalized system,
for which $\rho = 1$, has uniformly bounded trajectories. According to Theorem~\ref{th10},
an irreducible system is marginally stable. For reducible systems, the phenomenon of marginal instability may appear, even for the case of one matrix (when it has Jordan blocks corresponding to the largest by modulo eigenvalues).
For general reducible systems, the growth of trajectories is at most  polynomial and the power does not exceed the number of blocks
in factorization~(\ref{blocks}) with the maximal joint spectral radius. The following theorem extends the
results in~\cite{P06, CMS} from the classical case to arbitrary systems.
\begin{theorem}\label{th50}
For any triplet $\xi = (G, \cL, \cA)$, there is a constant $C_1 > 0$ such that
for every $x_1 \in L_i, \,i = 1, \ldots, n,$ there exists a trajectory $\{x_k\}_{k \in \n}$ such that
\begin{equation}\label{marg1}
\|x_k\| \ \ge \ C_1 \, \rho^k\, , \qquad k \in \n\, ,
\end{equation}
and there is a constant $C_2 > 0$ such that
\begin{equation}\label{marg2}
\max_{|\alpha| = k} \|P_{\alpha}\|\ \le \  C_2\, k^{\, r_1-1}\, \rho^k\, , \qquad k \in \n \, ,
\end{equation}
where $r_1$ is the total number of diagonal  blocks $\cA^{(1)}, \ldots , \cA^{(r)}$
in factorization~(\ref{blocks}) with $\rho(\xi^{(i)}) = \rho(\xi)$.
\end{theorem}
{\tt Proof.} Let $\rho(\xi) = 1$. Assume $\xi$ is irreducible. Then it possesses an invariant norm (Theorem~\ref{th20}), for which is suffices to prove~(\ref{marg1}), since all norms in a finite-dimensional space are equivalent.
From the definition of invariant norm it follows easily that there exists an infinite trajectory
$\{x_k\}_{k \in \n}$ such that $\|x_k\| = 1$ for all~$k \in \n$, which completes the proof for an irreducible
triplet. If $\xi$ is reducible, then consider its factorization~(\ref{blocks}) with irreducible systems
$\xi^{(i)}$ corresponding to the diagonal blocks.  By Corollary~\ref{c10}, at least one of them,
say, $\xi^{(s)}$ has the joint spectral radius one. Hence, there is a trajectory
$\{x_k^{(s)}\}_{k \in \n}$ of this system such that $\|x_k^{(s)}\| \ge  C_1$ for all~$k \in \n$.
On the other hand, for the corresponding trajectory $\{x_k\}_{k \in \n}$ of the full system~$\xi$,
each element $x_k^{(s)}$ is a projection of $x_k$ onto a subspace $L_i^{(s)}$
parallel to the other subspaces $L_i^{(t)}, \, t \ne s$. Hence $\|x_k\| \ge C_0\, \|x_k^{(s)}\| $,
where the constant $C_0$ does not depend on $x_k$. This proves~(\ref{marg1}).

We establish~(\ref{marg2}) for $r=2$ blocks in factorization~(\ref{blocks}), the case of general $r$
then follow by induction. We assume that the multinorms in both $\cL^{(1)}$ and $\cL^{(2)}$ are invariant.
If $r_1 = 2$, i.e., $\rho_1 = \rho_2 = 1$, the norm of the sum~(\ref{marg-sum}) is bounded above by
$\sum_{s=1}^k \|D_{i_{s+1}i_{s}}\| \le k\, C$, and we arrive at~(\ref{marg2}) with $r_1 = 2$.
If $\rho_1 = 1, \rho_2  = q < 1$, then this norm is bounded by $\sum_{s=1}^k q^{k-s}\, \|D_{i_{s+1}i_{s}}\| \, \le
\, \frac{1}{1-q} \, C$, and we arrive at~(\ref{marg2}) with $r_1 = 1$.

  {\hfill $\Box$}

\section{The Invariant polytope algorithm}\label{s-algor}

We give a short description of the Invariant polytope algorithm for exact computation
of the constrained JSR and for constructing an extremal polytopic multinorm. The main approach is very similar to the classical single space
case $(n=1)$ elaborated in detail in~\cite{GP13} (the basic idea traces back to the papers~\cite{P96} and~\cite{GWZ05}.)
Then we make the formal description and provide a criterion for its
convergence within finite time (Theorem~\ref{th60}).

A product of matrices $\Pi_{\alpha}$ corresponding to a path $\alpha$ is called {\em spectrum maximizing product} (in short, s.m.p.) if $[\rho(\Pi_{\alpha})]^{1/|\alpha|} = \rho(\xi)$. Inequality~(\ref{ineq}) shows that
we always have $[\rho(\Pi_{\alpha})]^{1/|\alpha|} \le \rho(\xi)$. So, an s.m.p. is a product for which
this inequality becomes equality. Even in the classical (single-space) case an s.m.p. may not exist~\cite{BTV}.

The idea of the algorithm is to select a {\em canditate s.m.p.} $\Pi_{\alpha}$ and prove that it is actually a real s.m.p.
by constructing an extremal polytope multinorm for $\xi$.
\smallskip

Due to the numerical computation we make use of a tolerance {\tt tol} in the computation which establishes whether
a vector is internal or external to a polytope.
\bigskip

Given the triplet $\xi=(G, \cL, \cA)$, where $G$ is a graph with $n$ nodes, $\cA = \{ A_{ji}^s \}$ where $i,\ j=1,\ \ldots, n$, $s =1,\ \ldots, N_{ji}$ and $N_{ji}$ is the number of edges connecting a node $i$ to the node $j$. We assume that
\begin{enumerate}
  \item $G$ is strongly connected
  \item $\Pi_\alpha$ is a candidate s.m.p. for $\xi$ of length $|\alpha|=N_\alpha$.
  \item the largest by modulo eigenvalue of $\Pi_\alpha$ is real.
  \item $\xi$ is irreducible (this assumption is for the sake of simplicity and can be omitted, see Remark~\ref{r100}).
\end{enumerate}

\begin{algorithm}[H] \label{algo}
\DontPrintSemicolon
\KwData{triplet $\xi=(G, \cL, \cA)$, $\Pi_\alpha$, tol}
\KwResult{$\hr(\xi)$, $\cV$}
\nl Scale the set $\cA$ and get ${\widetilde{\cA}} = \{ \widetilde{A}_{ji}^s=\rho(\Pi_\alpha)^{-1/|\alpha|}\,A_{ji}^s\}_{i,j,s}$ so that $\hr(\widetilde \xi)\geq 1$, with $\widetilde \xi=(G,\ \cL,\ \widetilde{\cA},\ )$, and $\rho(\widetilde\Pi_\alpha)=1$ \;
\nl Compute the leading eigenvectors $ {v}$ of $\widetilde\Pi_\alpha$ normalized with $\| {v}\|_2=1$ and of its $N_\alpha-1$ cyclic permutations. \;
\nl Set $k=0$ \;
\nl Define $\cV^{(0)} = \bigl\{V_i^{(0)}\bigr\}_{i=1}^n$
and $\cR^{(0)} = \bigl\{R_i^{(0)}\bigr\}_{i=1}^n$
where, for all $i$, $R_i^{(0)} = V_i^{(0)} \subset L_i$ contains the eigenvectors $\{ {v}_j\}$ which belong to the space~$L_i$ based on the path $\alpha$\footnote{See the illustrative examples in Section \ref{sec:IllEx}}.\;
\While{$R_i^{(k)} \neq \emptyset$ for at least one $i=1,\ \ldots,\ n$}{
\nl Set $k=k+1$\;
\For{$i = 1, \ldots , n$}
{Set $V_i^{(k)} = V_i^{(k-1)},  \, R_i^{(k)} = \emptyset$ \;
\For{\textbf{\textrm{all}} $v \in R_i^{(k-1)}$}{
\For{\textbf{\textrm{all}} edges from $i$ to $j$}{
\For{$s = 1, \ldots , N_{ji}$}{
Set $P_{j}^{(k)}={\rm absco} \, \left\{V_{j}^{(k-1)}\right\}$\;
\If{$\|\widetilde A_{ji}^s v\|_{P_{j}^{(k)}}\geq 1-\rm{tol}$}
{add $\widetilde A_{ji}^s v$ to the sets $V_j^{(k)}$ and $R_j^{(k)}$ \; }
}}}}}
$\Pi_\alpha$ is an s.m.p. for $\xi$, and $\hr(\xi)=\rho(\Pi_\alpha)^{1/{|\alpha|}}$ \;
$\cV = \cV^{(k)} = \left\{V_i^{(k)}\right\}_{i=1}^n$ is the set of vertices of the polytope extremal multinorm\;
\caption{The algorithm for computing the constrained JSR for a triplet $\xi=(G, \cL, \cA)$}
\end{algorithm}

If Algorithm \ref{algo} terminates after $N$th iteration, then we have the family of invariant polytopes~$\{P_j^{(N)}\}_{j=1}^n$
such that $P_j^{(N)} = {\rm co} \{\widetilde A_{ji} P_j^{(N)} \ | \ A_{ji} \in \cA_{ji}, \ i = 1, \ldots , n\}$.
The Minkowski norm defined by those polytopes is extremal, and hence $\rho(\widetilde \xi ) = 1$, which proves that
$\rho(\xi) = [\rho(\Pi_{\alpha})]^{1/|\alpha|}$.
\smallskip

Regarding the assumptions, first of all we observe that if the graph is not strongly connected we can always find a disjoint partition of its vertices so that $G = \bigsqcup_{\, i=1}^{\, r} G_i$ where each submultigraph $G_i$ is strongly connected and $\hr(\xi) = \max \, \{\hr(\xi_1), \ldots , \hr(\xi_r)\}$, ref. Corollary \ref{c5}.

To identify a candidate  s.m.p. $\Pi_\alpha$ for the triplet $\xi = (G, \cL, \cA)$, we fix some number $l_0$ and look among all the simple closed path $\alpha$ on $G$, with $|\alpha| = l \le l_0$, for the maximal value $\rho_{\alpha} = [\rho(\Pi_{\alpha})]^{1/|\alpha|}$. As mentioned in the assumptions, in this work we assume the leading eigenvalue $\lambda$ of~$\Pi_\alpha$, which is the largest by modulo eigenvalue, to be real. We observe that the ideas and the algorithm proposed in this paper extend as they are to the complex case.

Furthermore we assume that there exists a unique path $\alpha$ in $G$ such that $\rho_{\alpha}$ is maximal. It is easy to construct graphs where such path is not unique by allowing the same sequence of matrices in different paths of the graph. However in order to handle such cases a proper balancing of the vectors involved is required. We plan to study such kind of problems in the future.

If we consider the scaled set ${\widetilde{\cA}}$, as described in step 1 of Algorithm~1, the candidate s.m.p. becomes $\widetilde \Pi_\alpha$ with leading eigenvalue $\lambda=1$. The eigenvector $ {v}_1$ corresponding to $\lambda$ is called {\em leading eigenvector}. Let $ {v}_j =\widetilde \Pi_{\alpha_{j-1}} {v}_1, \, j= 2, \ldots , l$, where $\alpha_s$ is the prefix of the path $\alpha$ of length~$s$ (the first $s$ edges of $\alpha$). Thus, $\{ {v}_j\}_{j=1}^{N_\alpha}$ are the leading eigenvectors of all cyclic permutations of $\widetilde \Pi_\alpha$.

Following the path $\alpha$ we can assign each eigenvector $ {v}_j$ to the corresponding space $L_i$. Examples of this procedure are given in Section \ref{sec:IllEx}.

\begin{remark}\label{r100}
{\em Regarding the irreducibility, we observe that actually Algorithm~1 does not use it. If it produces
full-dimensional invariant polytopes, then we are done. The only trouble may occur in the case
when some polytopes~$P_i^{(N)}$ are not full dimensional. In this case the triple~$\xi$ is reducible, and we can proceed as follows.

In $k$th iteration we compute the number $d^{(k)}$ which is the sum of dimensions of the linear spans $L^{(k)}_i$ of the sets $V_i^{(k)}$ over $i=1, \ldots , n$. If $d^{(k)} = d^{(k-1)} \, < \, \sum_{j=1}^n {\rm dim}\, L_j$, then $L^{(k)}_i = L^{(k-1)}_i$ for all $i$. In this case, the triplet $\xi' = (G, \cL^{(k)}, \cA|_{\cL^{(k)}})$ is strictly embedded into $\xi$, and hence $\xi$ is reducible. Using Theorem~\ref{th30} we make the reduction to two triples $\xi'$ and $\xi''$ of smaller dimensions. We stop the algorithm and apply it to $\xi'$ and $\xi''$ separately. If $d^{(k)} >  d^{(k-1)}$ or if $d^{(k-1)} = \sum_{j=1}^n {\rm dim}\, L_j$, then we simply continue the iterations. A complete analysis of such procedure is out of the scope of this work.}
\end{remark}

If the algorithm terminates within finite time, then it proves that the
chosen candidate is indeed an s.m.p. and gives the corresponding polytope extremal norm. Although  there are simple examples,
when Algorithm~1 does not terminate within finite time, this phenomenon is believed to be rare in practice.
In the single-space case, all numerical experiments made with randomly generated matrices and with matrices from applications, Algorithm~1 did terminate in finite time providing an invariant polytope (see~\cite{GP13} for examples and statistics).

The following theoretical criterion ensures the convergence of Algorithm~1. It generalizes Theorem~4 from~\cite{GP13}
proved for the single-space case and uses the notion of dominant product.
\begin{defi}\label{d70}
Let $\xi = (G, \cL,  \cA)$ be an arbitrary triplet. A closed simple path $\alpha$ along $G$ and the
corresponding products $\Pi_{\alpha}$ are called dominant,  if there is a constant
$q < 1$ such that the spectral radii of products of the normalized family $\widetilde \cA = [\rho(\Pi_{\alpha})]^{-1/|\alpha|}\cA$ corresponding to all other simple paths which are not cyclic permutations of~$\alpha$, are smaller than~$q$.
\end{defi}

In view of Theorem~\ref{th40}, if the path $\alpha$ is dominant, then $\rho(\xi) = [\Pi_{\alpha}]^{1/|\alpha|}$.
Thus, a dominant product is always an  s.m.p., but, in general, not vice versa.

\smallskip

\begin{theorem}\label{th60}  For a given triplet $\xi = (G, \cL, \cA)$ and for a given
initial path $\alpha$ (candidate s.m.p.), Algorithm~1 terminates within finite time if and only if~$\alpha$ is dominant and
the leading eigenvalue of $\Pi_{\alpha}$ is unique and simple.
\end{theorem}
The proof is actually the same as the proof of the single-space case in~\cite[Theorem~4]{GP13}, and we omit it.

\subsection{Illustrative examples}\label{sec:IllEx}

We start with two examples to show in details how the method works, whereas in Example \ref{ex:three} we run some statistics to show the performance of the algorithm.

\begin{ex}\label{ex:one}

Given the triplet $\xi=(G, \cL, \cA)$ with $\cL=\left\{L_i\right\}_{i\in\{1,\ \ldots,\ 4\}}$, $L_i=\re^2, i=1,\ \ldots,\ 4$,
$$\cA=\left\{A_i\right\}_{i\in\{1,\ \ldots,\ 4\}}
            =\left\{-I_2,\
            \left[
              \begin{array}{cc}
                0 & 1 \\
                -1 & -1 \\
              \end{array}
            \right],\
            \left[
              \begin{array}{cc}
                -1 & 1 \\
                -1 & 0 \\
              \end{array}
            \right],\
            \left[
              \begin{array}{cc}
                1 & 2 \\
                0 & 1 \\
              \end{array}
            \right]
\right\}$$

and $G$ given in Figure \ref{fig:GraphG}

\begin{figure}[ht]
\begin{center}
\begin{tikzpicture}[->,>=stealth',shorten >=1pt,auto,node distance=2cm,
                    thick,main node/.style={circle,draw,font=\sffamily\bfseries}]

  \node[main node] (1) {$L_{1}$};
  \node[main node] (2) [below of=1] {$L_2$};
  \node[main node] (3) [right of=2] {$L_3$};

  \path[every node/.style={font=\sffamily\footnotesize}]
    (1) edge [bend right] node [below left] {$A_1,A_3$} (2)
        edge  node  [right]  {$A_1$} (3)
        edge [loop left] node [left] {$A_1,A_3$} (1)
    (2) 
        edge node[above right] {$A_2$} (1)
        edge [loop left] node {$A_2$} (2)
        edge [bend right] node {$A_2$} (3)
    (3) edge [bend right] node[above] {$A_4$} (1);
\end{tikzpicture}
\caption{ Graph $G$ }\label{fig:GraphG}
\end{center}
\end{figure}

By an exhaustive search among all closed path of length $l_0\leq 10$ we identify the candidate s.m.p. $ \Pi = A_3     A_2     A_3     A_4     A_1     A_4  A_2$ corresponding to the closed path $\alpha$ shown in Figure \ref{fig:pathAlpha}, $|\alpha|=7$.

We scale the set of matrices to get ${\widetilde{\cA}} = \{ \widetilde{A}_i={A_i}\rho(\Pi)^{-\frac{1}{7}}\}_{i=1}^{4}$, so that $\hr(\widetilde\xi)\geq 1$, with $\widetilde \xi=({\widetilde{\cA}},\ G,\ \cL)$, and $\rho(\widetilde {\Pi})= 1$.

We denote the leading eigenvectors of the candidate s.m.p. $\widetilde \Pi$ and its cyclic permutations by $\left\{ {v}_i\right\}_{i\in\{1, \ldots 7\}}$. Assuming that $ {v}_1 =\widetilde \Pi  {v}_1$, then $ {v}_2=\widetilde{A}_2   {v}_1,\  {v}_3=\widetilde{A}_4   {v}_2,\ {v}_4=\widetilde{A}_1  {v}_3,\  {v}_5=\widetilde{A}_4   {v}_4,\  {v}_6=\widetilde{A}_3   {v}_5,\  {v}_7=\widetilde{A}_2   {v}_6,\  {v}_1=\widetilde{A}_3   {v}_7$.

We have ${v}_3,\ {v}_5, \ {v}_7\in L_1$, ${v}_1,\ {v}_6\in L_2$, and ${v}_2,\ {v}_4\in L_3$, see fig.  \ref{fig:Ex1_Step0}.

\tikzset{%
   peer/.style={draw,circle,black,bottom color=white, top color= white, text=black, minimum width=25pt},
   superpeer/.style={draw, circle,  left color=white, text=black, minimum width=25pt},
   point/.style = {fill=black,inner sep=1pt, circle, minimum width=5pt,align=right,rotate=60},
   }

\def \n {7}

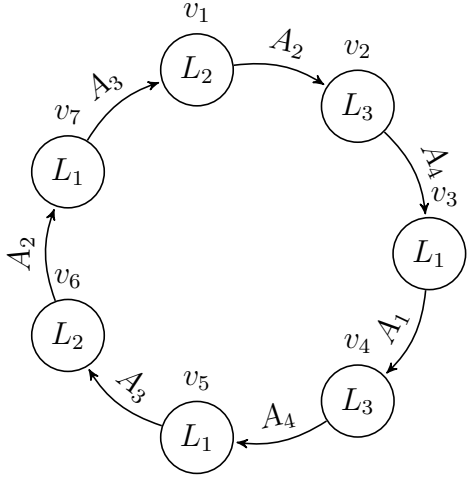
\begin{figure}[ht]
\begin{center}
\begin{tikzpicture}[->,>=stealth',shorten >=1pt,auto,node distance=2.8cm, semithick]
\node[superpeer,label={$\quad {v}_3$}]    (A)   at ({360/\n * (1 - 1)}:2.5){$L_1$};
\node[peer,label={$ {v}_2$}]         (B)   at ({360/\n * (2 - 1)}:2.5){$L_3$};
\node[peer,label={$ {v}_1$}]         (C)   at ({360/\n * (3 - 1)}:2.5){$L_{2}$};
\node[peer,label={$ {v}_7$}]         (E)   at ({360/\n * (4 - 1)}:2.5){$L_1$};
\node[peer,label={$ {v}_6$}]         (A1)  at ({360/\n * (5 - 1)}:2.5){$L_2$};
\node[peer,label={$ {v}_5$}]         (B1)  at ({360/\n * (6 - 1)}:2.5){$L_1$};
\node[peer,label={$ {v}_4$}]         (C1)  at ({360/\n * (7 - 1)}:2.5){$L_3$};

\path (B) edge [color=black, bend left=20, above, sloped] node[] {$ A_4$}(A);

\path (C) edge [color=black, bend left=20, above, sloped] node[] {$ A_2$}(B);

\path (E) edge [color=black, bend left=20, above, sloped] node[] {$ A_3$}(C);

\path (A1) edge [color=black, bend left=20, above, sloped] node[]{$ A_2$} (E);

\path (A) edge [color=black, bend left=20, above, sloped] node[] {$ A_1$} (C1);

\path (B1) edge [color=black, bend left=20, above, sloped] node[] {$ A_3$}(A1);

\path (C1) edge [color=black, bend left=20, above, sloped] node[] {$ A_4$}(B1);

\end{tikzpicture}
\caption{Closed path $\alpha$ of the candidate s.m.p. $\Pi$ for the triplet $ \xi$}\label{fig:pathAlpha}
\end{center}
\end{figure}

\begin{figure}
\centering
\begin{minipage}[b]{.3\linewidth}
  \centering
  \centerline{\includegraphics[width=\linewidth]{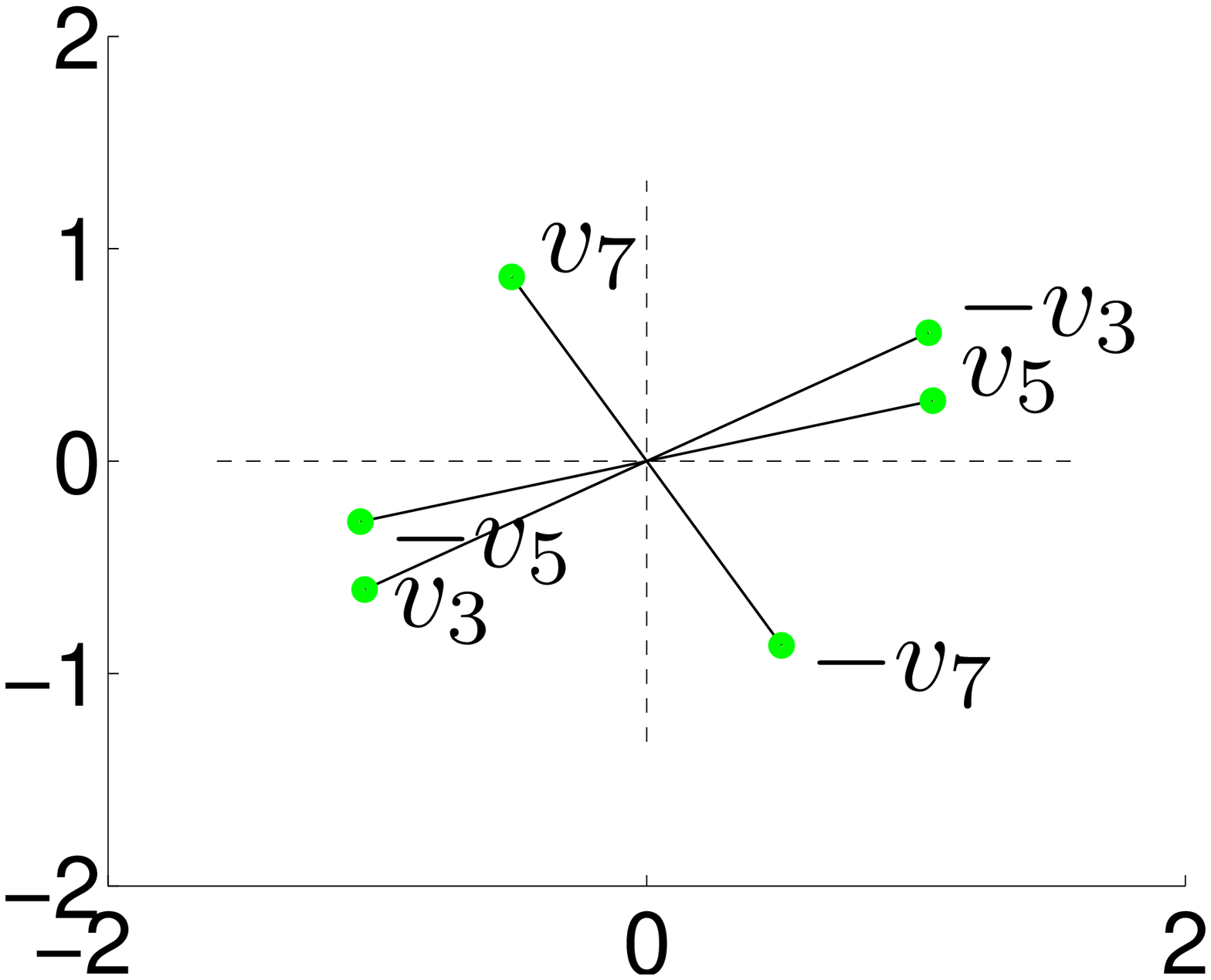}}
\small\centerline{$L_{1}$}\medskip
\end{minipage}
\hfill
\begin{minipage}[b]{0.3\linewidth}
  \centering
  \centerline{\includegraphics[width=\linewidth]{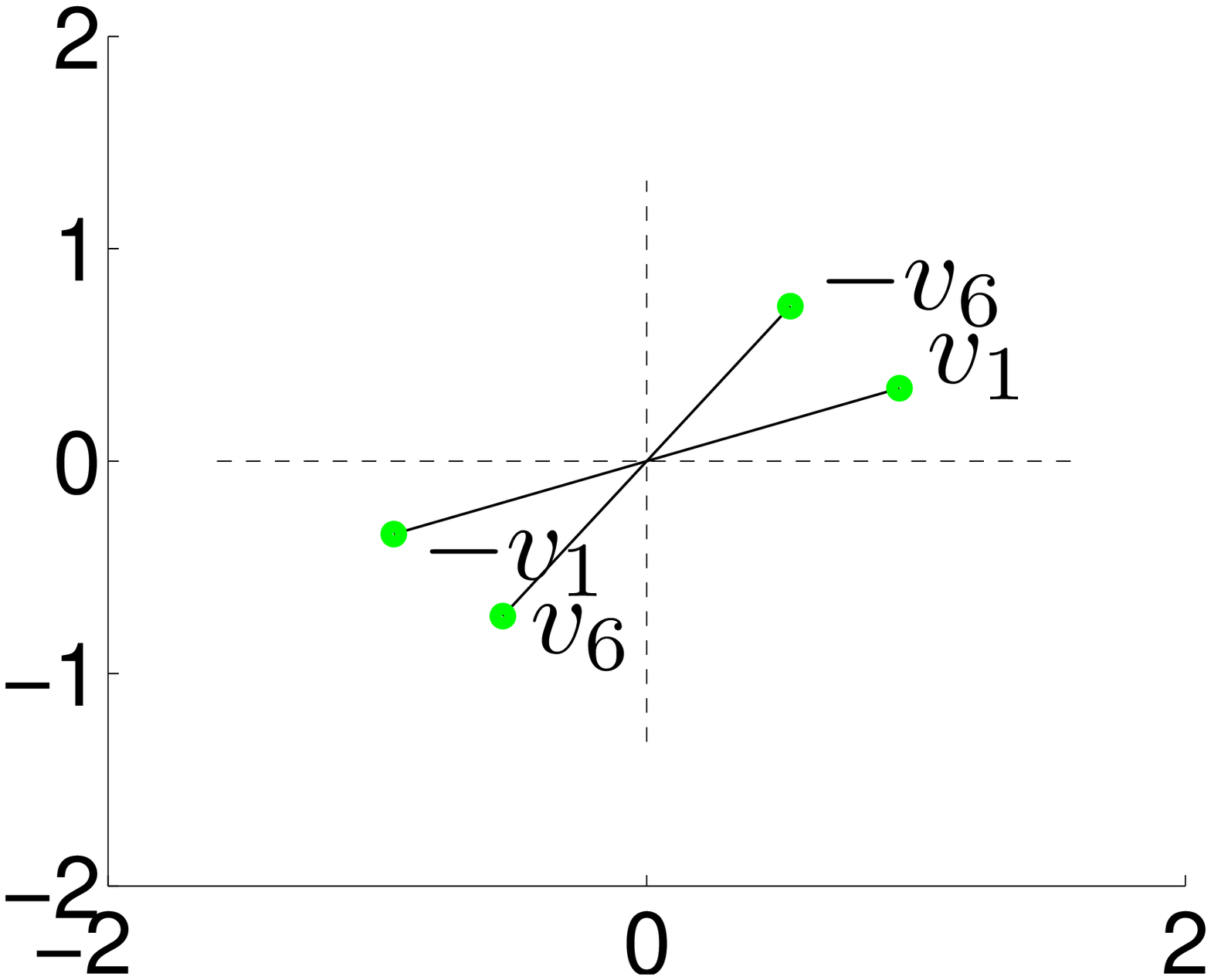}}
\small\centerline{$L_{2}$}\medskip
\end{minipage}
\hfill
\begin{minipage}[b]{0.3\linewidth}
  \centering
  \centerline{\includegraphics[width=\linewidth]{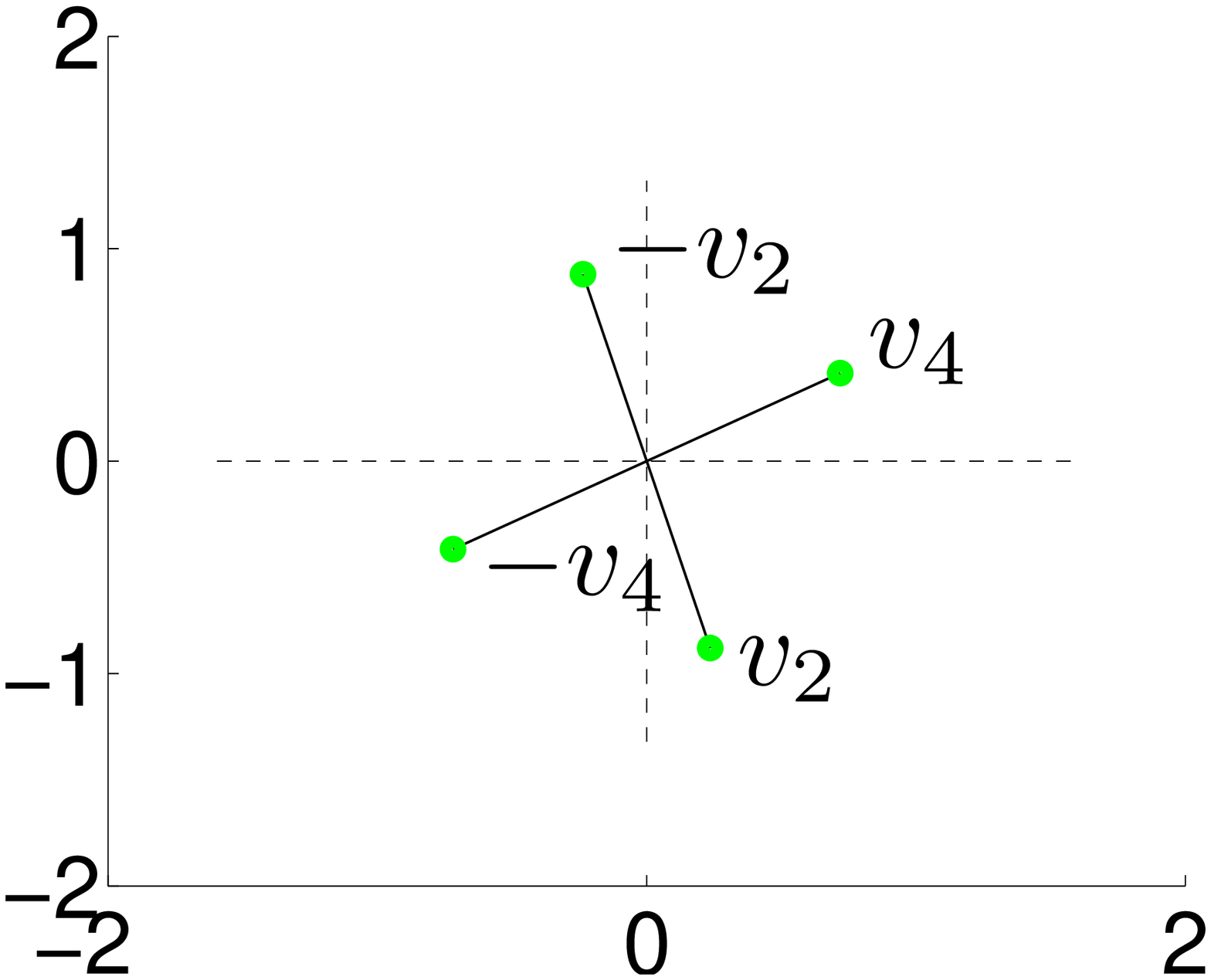}}
\small\centerline{$L_{3}$}\medskip
\end{minipage}
%
\caption{Leading eigenvectors of $\widetilde\Pi$ and its cyclic permutations each of them assigned to the corresponding space.}\label{fig:Ex1_Step0}
\end{figure}

\begin{figure}
\centering
\begin{minipage}[b]{.3\linewidth}
  \centering
  \centerline{\includegraphics[width=\linewidth]{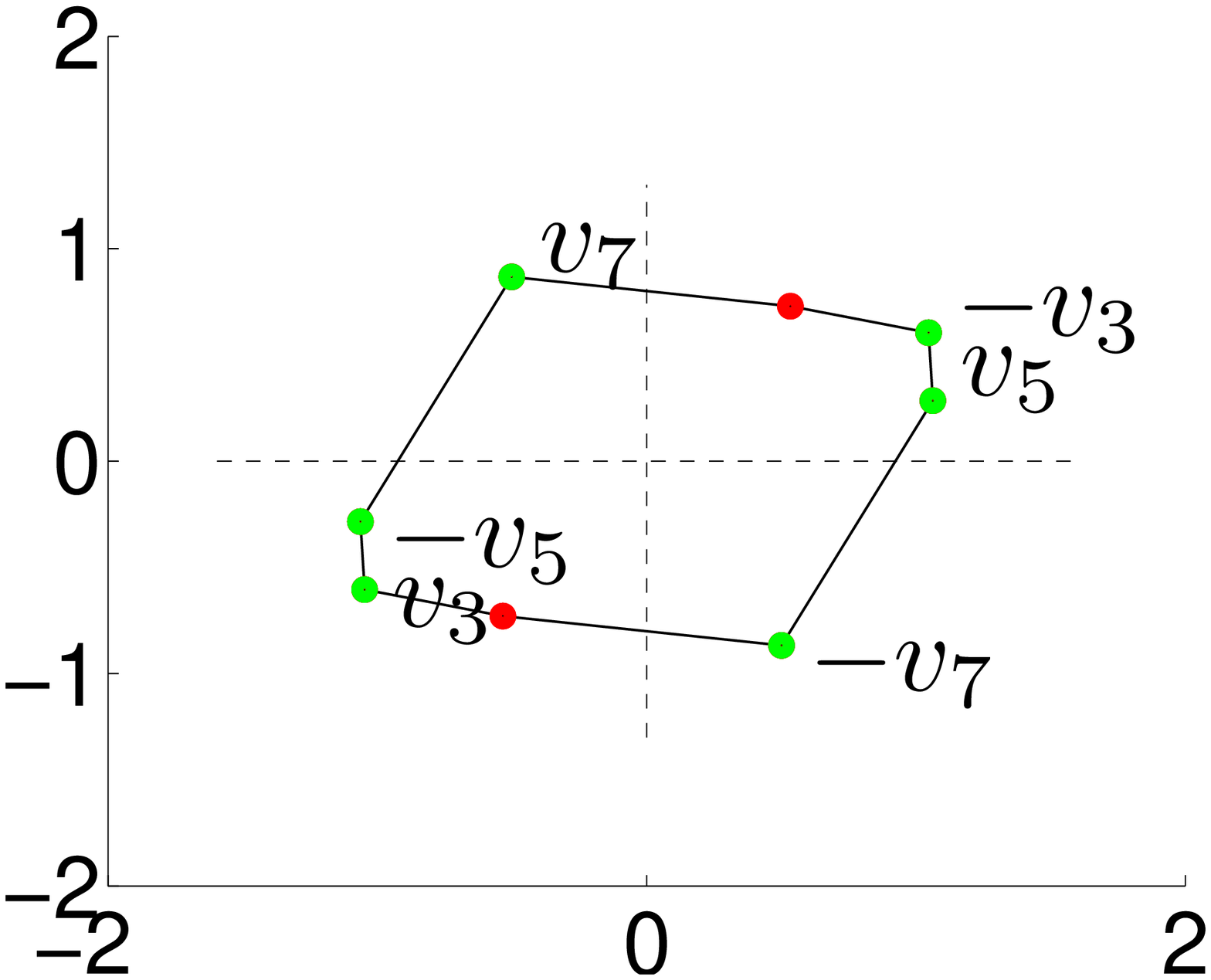}}
\small\centerline{$L_{1}$}\medskip
\end{minipage}
\hfill
\begin{minipage}[b]{0.3\linewidth}
  \centering
  \centerline{\includegraphics[width=\linewidth]{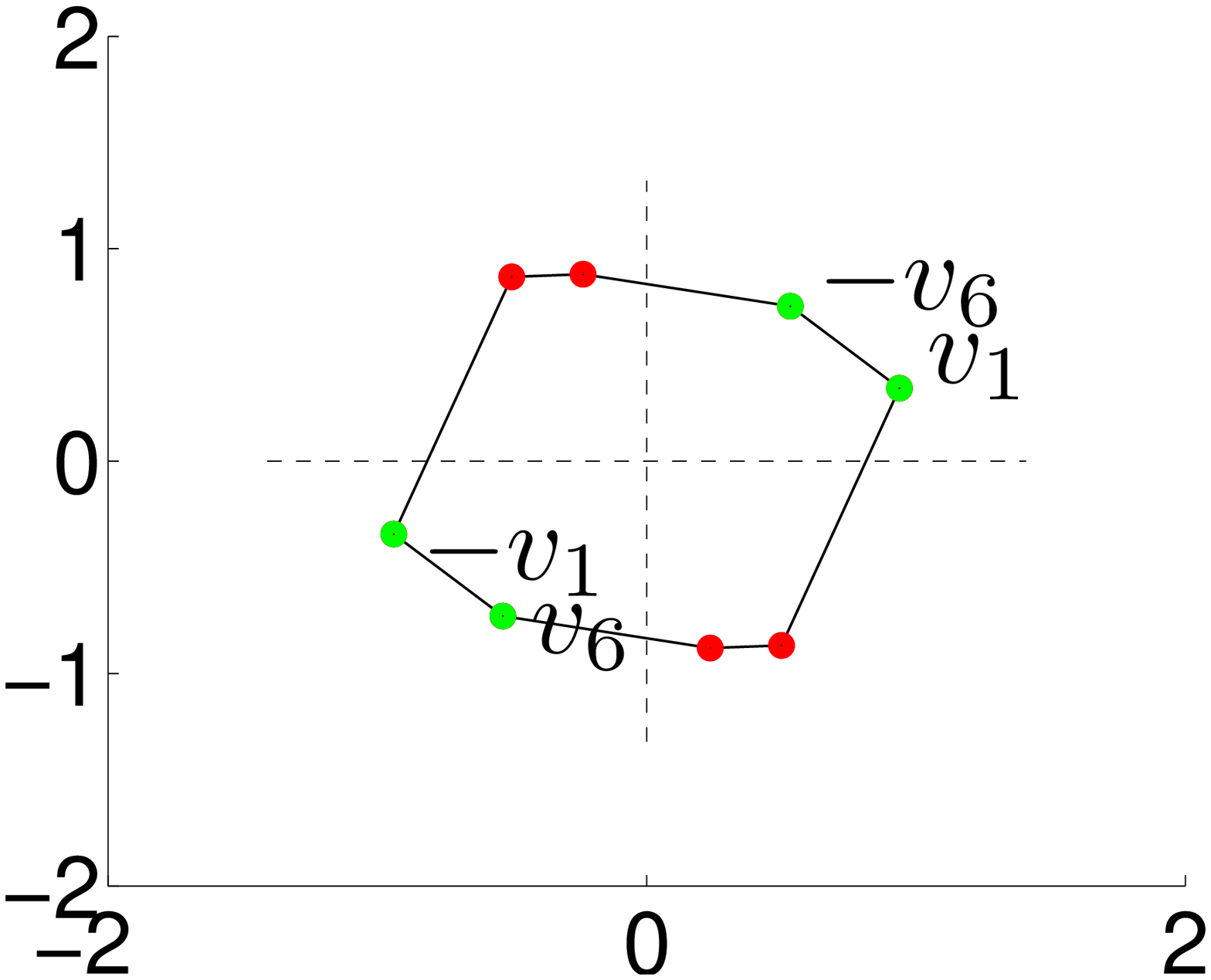}}
\small\centerline{$L_{2}$}\medskip
\end{minipage}
\hfill
\begin{minipage}[b]{0.3\linewidth}
  \centering
  \centerline{\includegraphics[width=\linewidth]{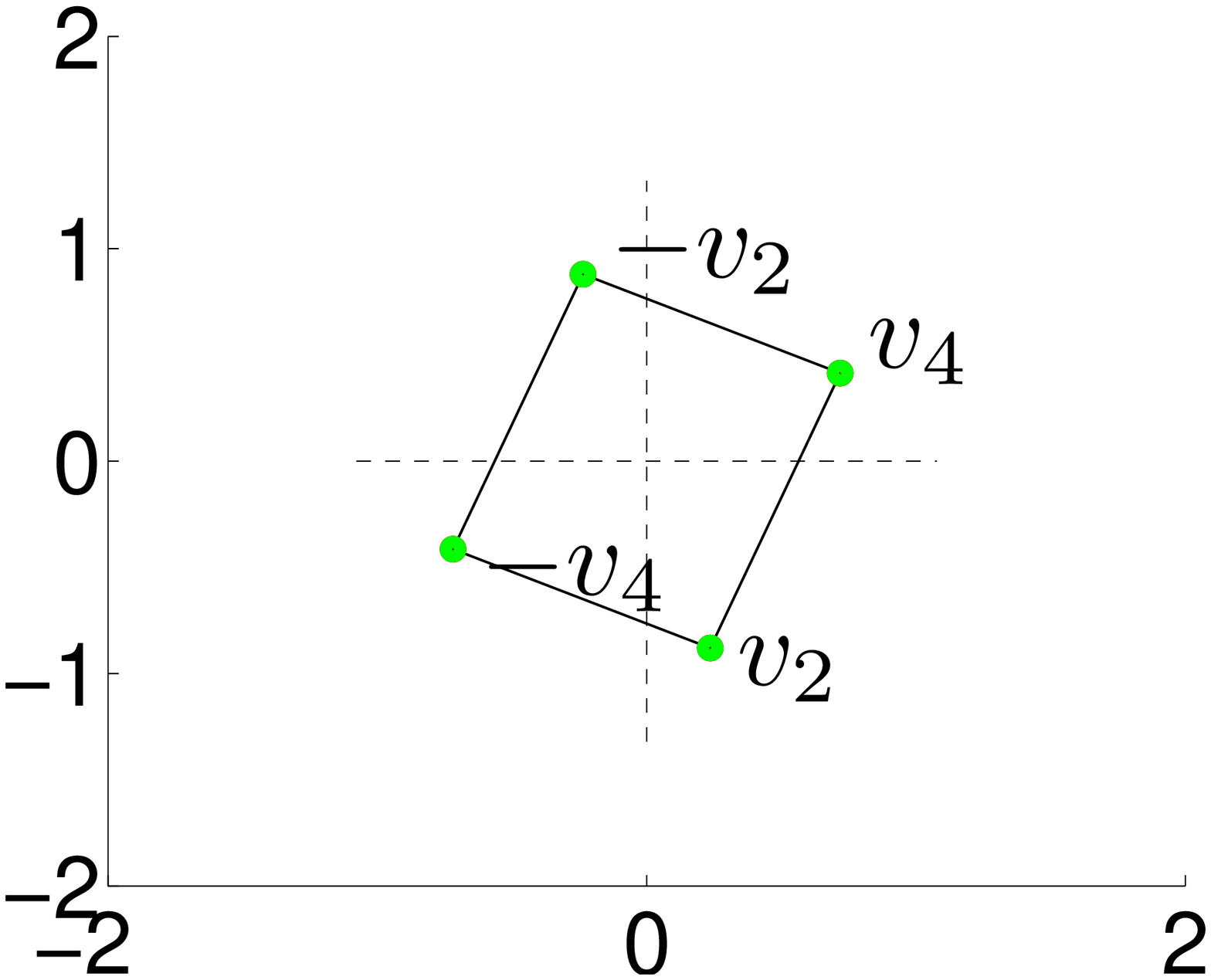}}
\small\centerline{$L_{3}$}\medskip
\end{minipage}
%
\caption{Vertices and their symmetrized convex hulls after one step of the algorithm}\label{fig:Ex1_Step1}
\end{figure}

\begin{figure}
\centering
\begin{minipage}[b]{.3\linewidth}
  \centering
  \centerline{\includegraphics[width=\linewidth]{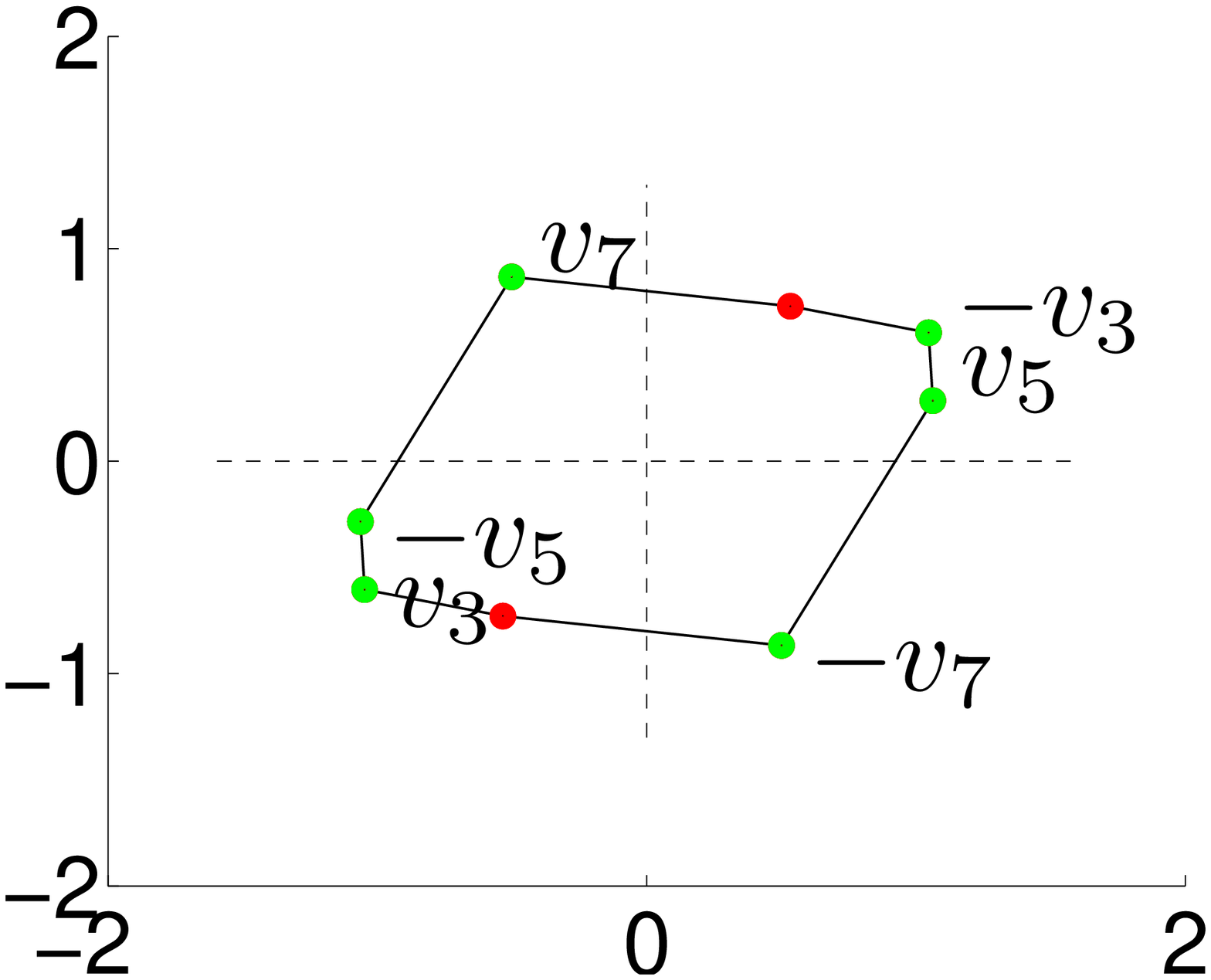}}
\small\centerline{$L_{1}$}\medskip
\end{minipage}
\hfill
\begin{minipage}[b]{0.3\linewidth}
  \centering
  \centerline{\includegraphics[width=\linewidth]{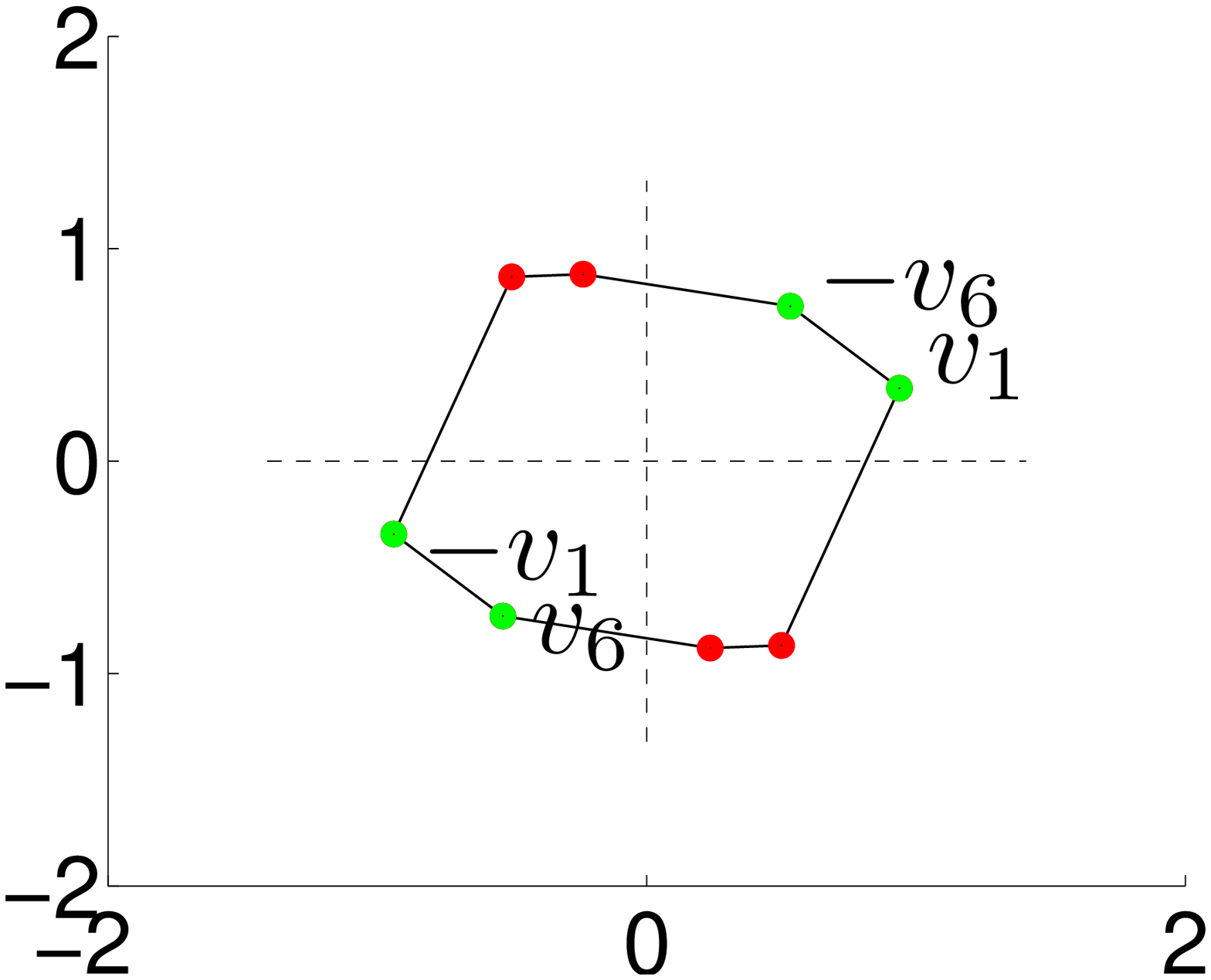}}
\small\centerline{$L_{2}$}\medskip
\end{minipage}
\hfill
\begin{minipage}[b]{0.3\linewidth}
  \centering
  \centerline{\includegraphics[width=\linewidth]{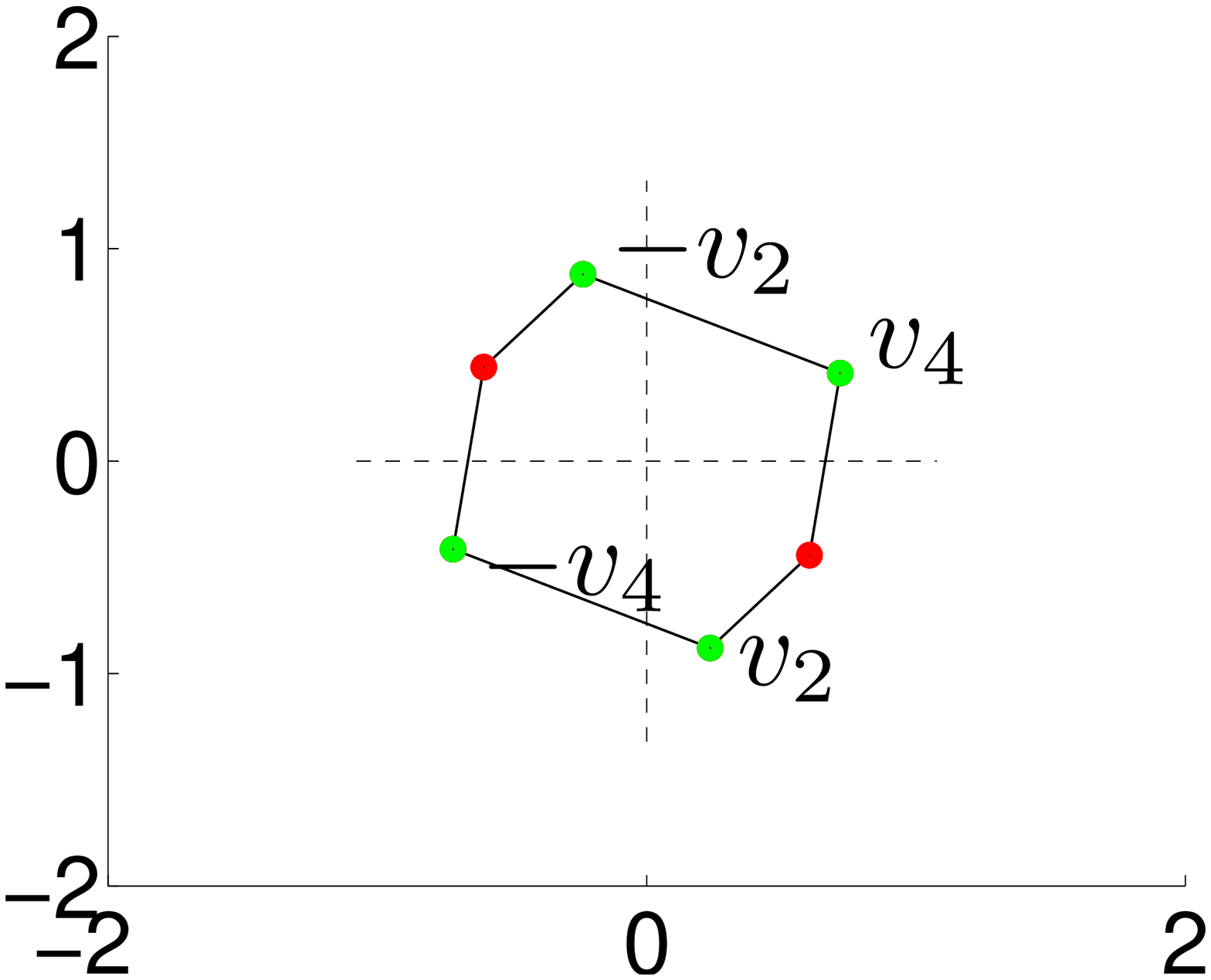}}
\small\centerline{$L_{3}$}\medskip
\end{minipage}
%
\caption{Invariant polytopes for $\xi$}\label{fig:Ex1_Step2}
\end{figure}

After two steps of Algorithm \ref{algo} we have invariant polytopes which are plotted in Figure \ref{fig:Ex1_Step2}.

So we can conclude that $\Pi$ is an  s.m.p. for $\xi$ and $\hr(\xi)=\rho(\Pi)^{\frac{1}{7}}= 1.456846\ldots$.

If we study the unconstrained problem using, for instance, the Matlab code presented in \cite{CP15} and posted on MatlabCentral\footnote{\url{http://www.mathworks.com/matlabcentral/fileexchange/36460-joint-spectral-radius-computation}}
we find the candidate s.m.p. $Q=A_4 A_3 A_4 A_4 A_2$. After scaling the set $\cA$ by $\rho(Q)^{1/5}$ to get $\bar{\cA}$, we can use the technique described in \cite{GP13} to construct in 4 steps a extremal polytopic norm whose unit ball contains the following 7 vectors: $ {v}_0$ which is the leading eigenvector of $\bar Q$, $\bar A_2  {v}_0$,  $\bar A_4  {v}_0$,  $\bar A_4 \bar A_2  {v}_0$,  $\bar A_4 \bar A_4 \bar A_2  {v}_0$,   $\bar A_2 \bar A_4 \bar A_2  {v}_0$,  $\bar A_3 \bar A_4 \bar A_4 \bar A_2  {v}_0$. This allows to conclude that $Q$ is an s.m.p. for $\cA$ and that $\rho(\cA)~=~\rho(Q)^{1/5}~=~ 1.693476\ldots$.

We observe that the matrix $A_1$, which in this example is the negative identity matrix, clearly does not count towards the computation of the unconstrained JSR. However the same matrix is fundamental for the computation of the constrained JSR. As a matter of fact it does appear in the constrained s.m.p. $\Pi$.

\end{ex}

\begin{ex}\label{ex:two}

We consider now the case where the dimensions of the spaces can be different each other and some matrices appear in edges corresponding to different spaces. For the sake of simplicity, we consider only one and two dimensional spaces $L_j$, however we recall that the proposed algorithm works with any dimension. The triplet is $\xi=(G, \cL, \cA)$ with $\cL=\left\{L_i\right\}_{i\in\{1,\ \ldots,\ 3\}}$, $L_1 = L_3 = \re^2, \, L_2 = \re^1$, the operators
$$
\cA\ =\ \left\{A_i\right\}_{i\in\{1,\ldots 4\}}
            =\left\{\left[
              \begin{array}{cc}
                0 & 1 \\
                -1 & 0 \\
              \end{array}
            \right],\
            \left[
              \begin{array}{c}
                1 \\
                -1\\
              \end{array}
            \right],\
            \left[
              \begin{array}{cc}
                1 & 2 \\
              \end{array}
            \right],\
            \left[
              \begin{array}{cc}
                1 & -1 \\
                1 & 1 \\
              \end{array}
            \right]
\right\}
$$
and the graph $G$ depicted in Figure \ref{fig:GraphGprime}

\begin{figure}[ht]
\begin{center}
\begin{tikzpicture}[->,>=stealth',shorten >=1pt,auto,node distance=2cm,
                    thick,main node/.style={circle,draw,font=\sffamily\bfseries}]

  \node[main node] (1) {$L_{1}$};
  \node[main node] (2) [below left of=1] {$L_2$};
  \node[main node] (3) [below right of=1] {$L_3$};

  \path[every node/.style={font=\sffamily\footnotesize}]
    (1) edge [loop above] node [right] {$A_1$} (1)
        edge  [bend left] node  [right]  {$A_1,A_4$} (3)
    (2) edge [bend left] node[below right] {$A_2$} (1)
    (3) edge [bend left] node[below] {$A_4$} (1)
        edge [bend left] node [below left] {$A_3$} (2)
        edge [loop right] node [above] {$A_1$} (3);
\end{tikzpicture}
\caption{Graph $G$}\label{fig:GraphGprime}
\end{center}
\end{figure}

By an exhaustive search among all closed path of length $l_0\leq 10$ we identify the candidate s.m.p. $ \Pi= A_3     A^3_4    A_2$ corresponding to the closed path $\alpha$ shown in Figure \ref{fig:pathAlphaPrime}.

\tikzset{%
   peer/.style={draw,circle,black,bottom color=white, top color= white, text=black, minimum width=25pt},
   superpeer/.style={draw, circle,  left color=white, text=black, minimum width=25pt},
   point/.style = {fill=black,inner sep=1pt, circle, minimum width=5pt,align=right,rotate=60},
   }

\def \n {5}

\begin{figure}[ht]
\begin{center}
\begin{tikzpicture}[->,>=stealth',shorten >=1pt,auto,node distance=2.8cm, semithick]
\node[superpeer,label={$ {v}_3$}]    (A)   at ({360/\n * (1 - 1)}:2.5){$L_3$}; 
\node[peer,label={$ {v}_2$}]         (B)   at ({360/\n * (2 - 1)}:2.5){$L_1$}; 
\node[peer,label={$ {v}_1$}]         (C)   at ({360/\n * (3 - 1)}:2.5){$L_{2}$}; 
\node[peer,label={$ {v}_5$}]         (B1)  at ({360/\n * (4 - 1)}:2.5){$L_3$}; 
\node[peer,label={$ {v}_4$}]         (C1)  at ({360/\n * (5 - 1)}:2.5){$L_1$}; 

\path (B1) edge [color=black, bend left=20, above, sloped] node[] {$ A_3$}(C);

\path (C) edge [color=black, bend left=20, above, sloped] node[] {$ A_2$}(B);

\path (B) edge [color=black, bend left=20, above, sloped] node[] {$ A_4$}(A);

\path (A) edge [color=black, bend left=20, above, sloped] node[] {$ A_4$} (C1);

\path (C1) edge [color=black, bend left=20, above, sloped] node[] {$ A_4$}(B1);
\end{tikzpicture}
\caption{Closed path $\alpha$ of the candidate s.m.p. $ \Pi$ for the triplet $ \xi$}\label{fig:pathAlphaPrime}
\end{center}
\end{figure}
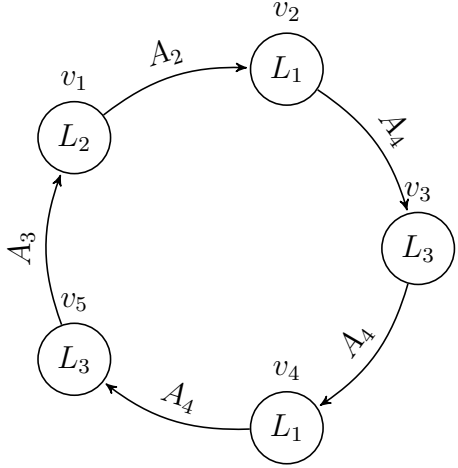

We scale the set ${\cA}$ and get ${\widetilde{\cA}} = \{ \widetilde{A}_i={A_i}\rho(\Pi)^{-\frac{1}{5}}\}_{i}$, so that $\hr(\widetilde\xi)\geq 1$, with $\widetilde \xi=({\widetilde{\cA}},\ G,\ \cL)$, and $\rho(\widetilde {\Pi})= 1$.

We denote the leading eigenvectors of the candidate s.m.p. $\widetilde \Pi$ and its cyclic permutations by $\left\{ {v}_i\right\}_{i\in\{1,\ldots 5\}}$. Assuming that $ {v}_1 =\widetilde \Pi  {v}_1$, then $ {v}_2=\widetilde{A}_2   {v}_1,\  {v}_3=\widetilde{A}_4   {v}_2,\ {v}_4=\widetilde{A}_4  {v}_3,\  {v}_5=\widetilde{A}_4   {v}_4,\  {v}_1=\widetilde{A}_3   {v}_5$.

We have ${v}_2,\ {v}_4\in L_1$, ${v}_1\in L_2$, and ${v}_3,\ {v}_5\in L_3$, see fig. \ref{fig:Ex2_Step0}.

\begin{figure}
\centering
\begin{minipage}[b]{.3\linewidth}
  \centering
  \centerline{\includegraphics[width=\linewidth]{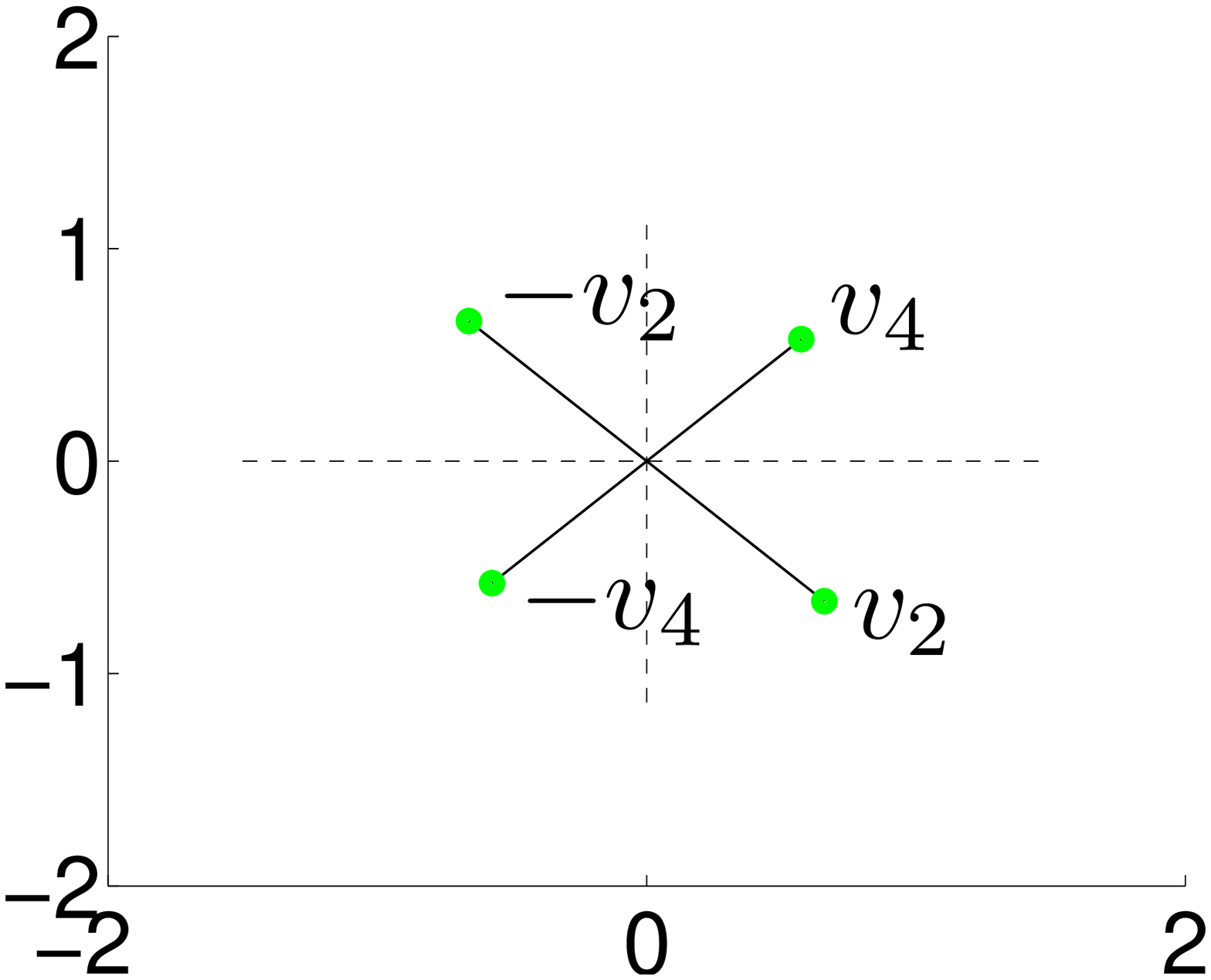}}
\small\centerline{$L_{1}$}\medskip
\end{minipage}
\hfill
\begin{minipage}[b]{0.3\linewidth}
  \centering
  \centerline{\includegraphics[width=\linewidth]{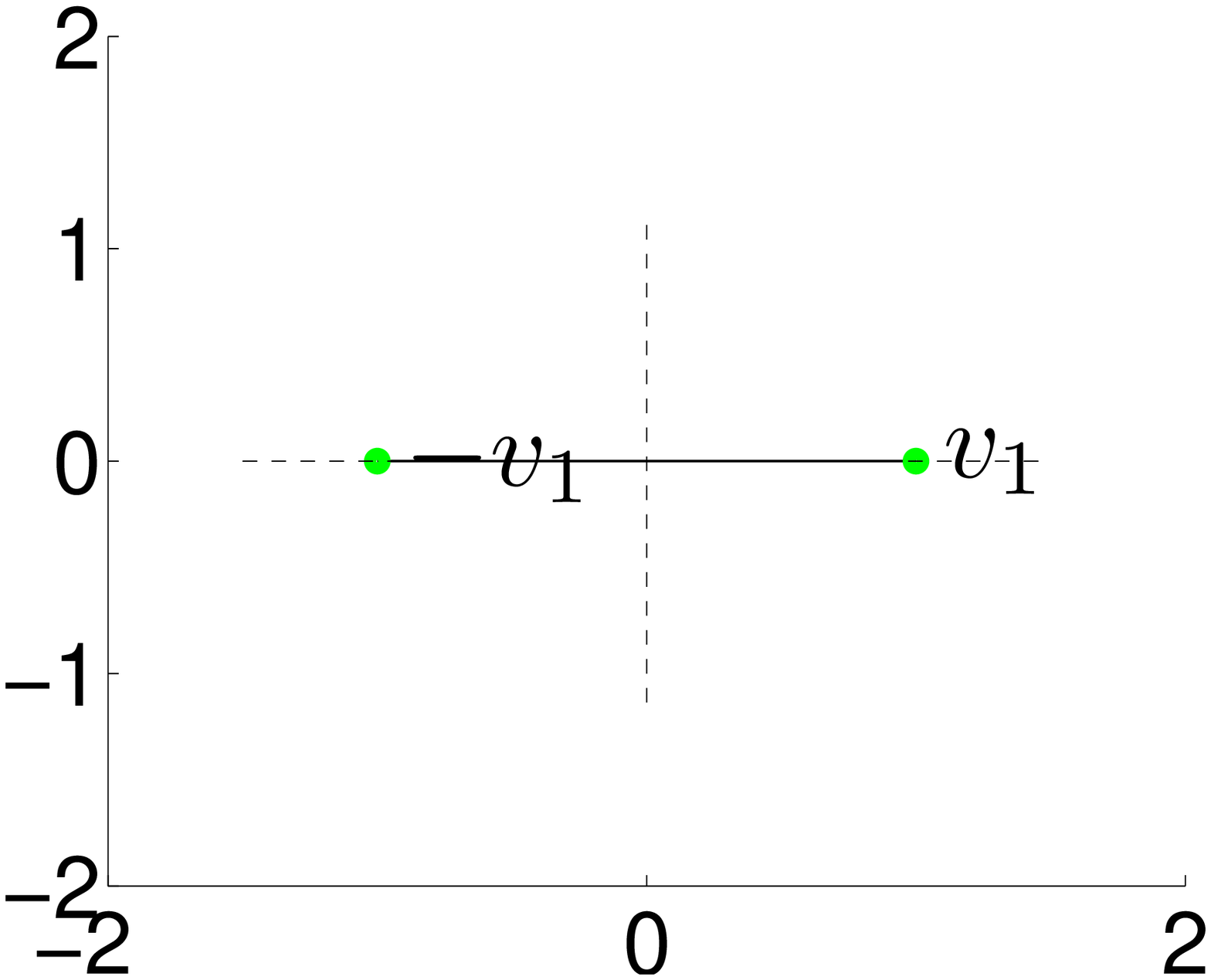}}
\small\centerline{$L_{2}$}\medskip
\end{minipage}
\hfill
\begin{minipage}[b]{0.3\linewidth}
  \centering
  \centerline{\includegraphics[width=\linewidth]{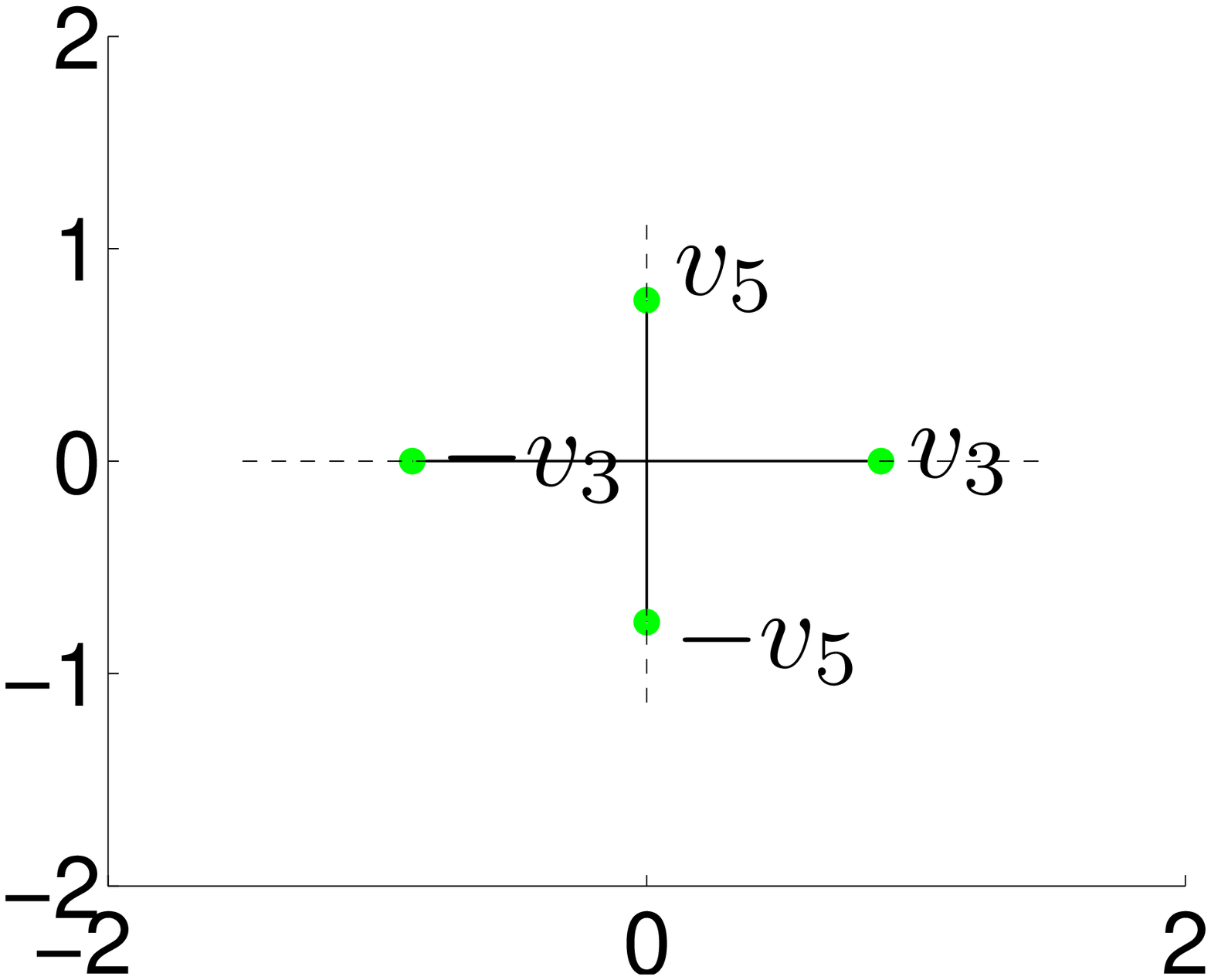}}
\small\centerline{$L_{3}$}\medskip
\end{minipage}
%
\caption{Leading eigenvectors of $\widetilde\Pi$ and its cyclic permutations each of them assigned to the corresponding space.}\label{fig:Ex2_Step0}
\end{figure}

\begin{figure}
\centering
\begin{minipage}[b]{.3\linewidth}
  \centering
  \centerline{\includegraphics[width=\linewidth]{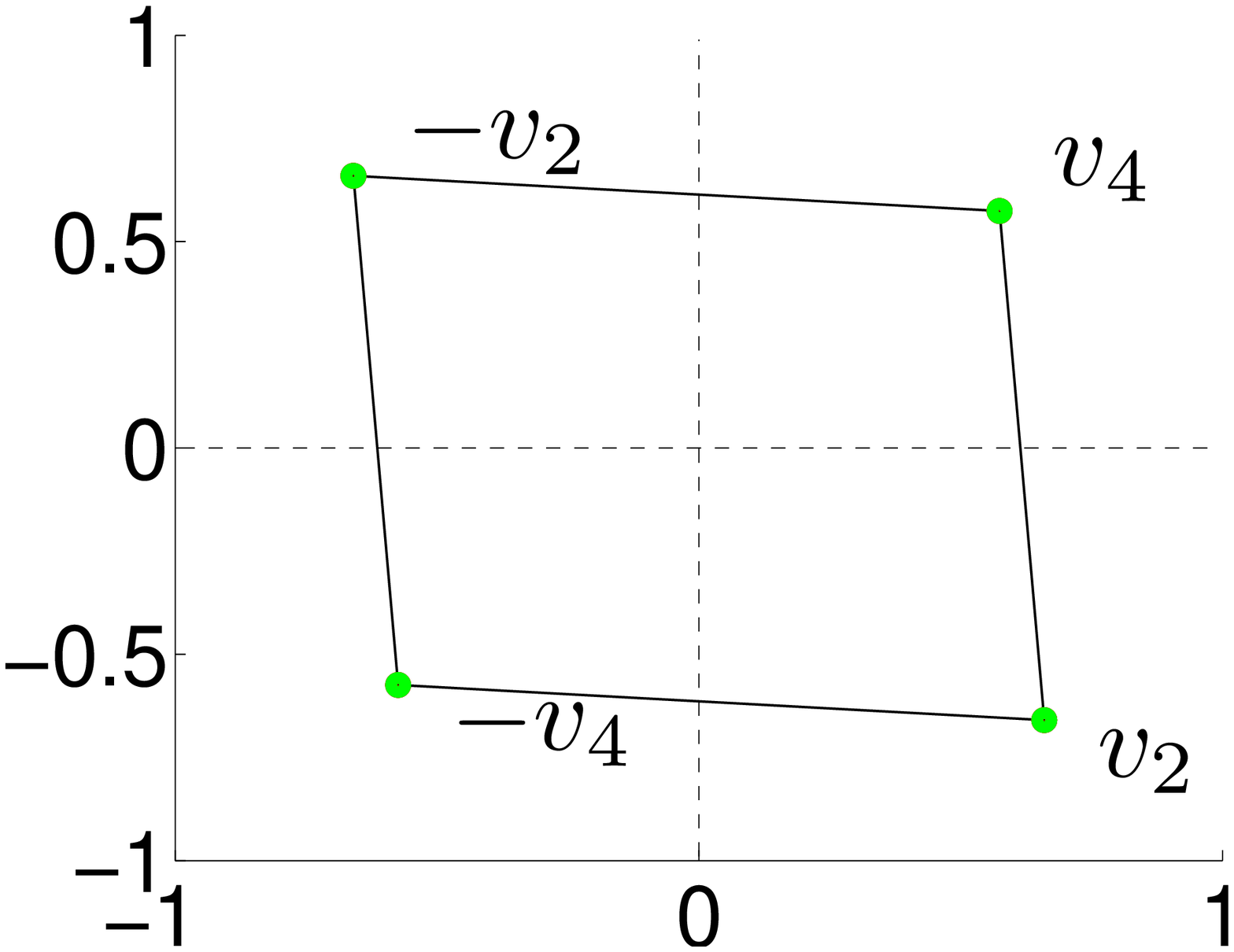}}
\small\centerline{$L_{1}$}\medskip
\end{minipage}
\hfill
\begin{minipage}[b]{0.3\linewidth}
  \centering
  \centerline{\includegraphics[width=\linewidth]{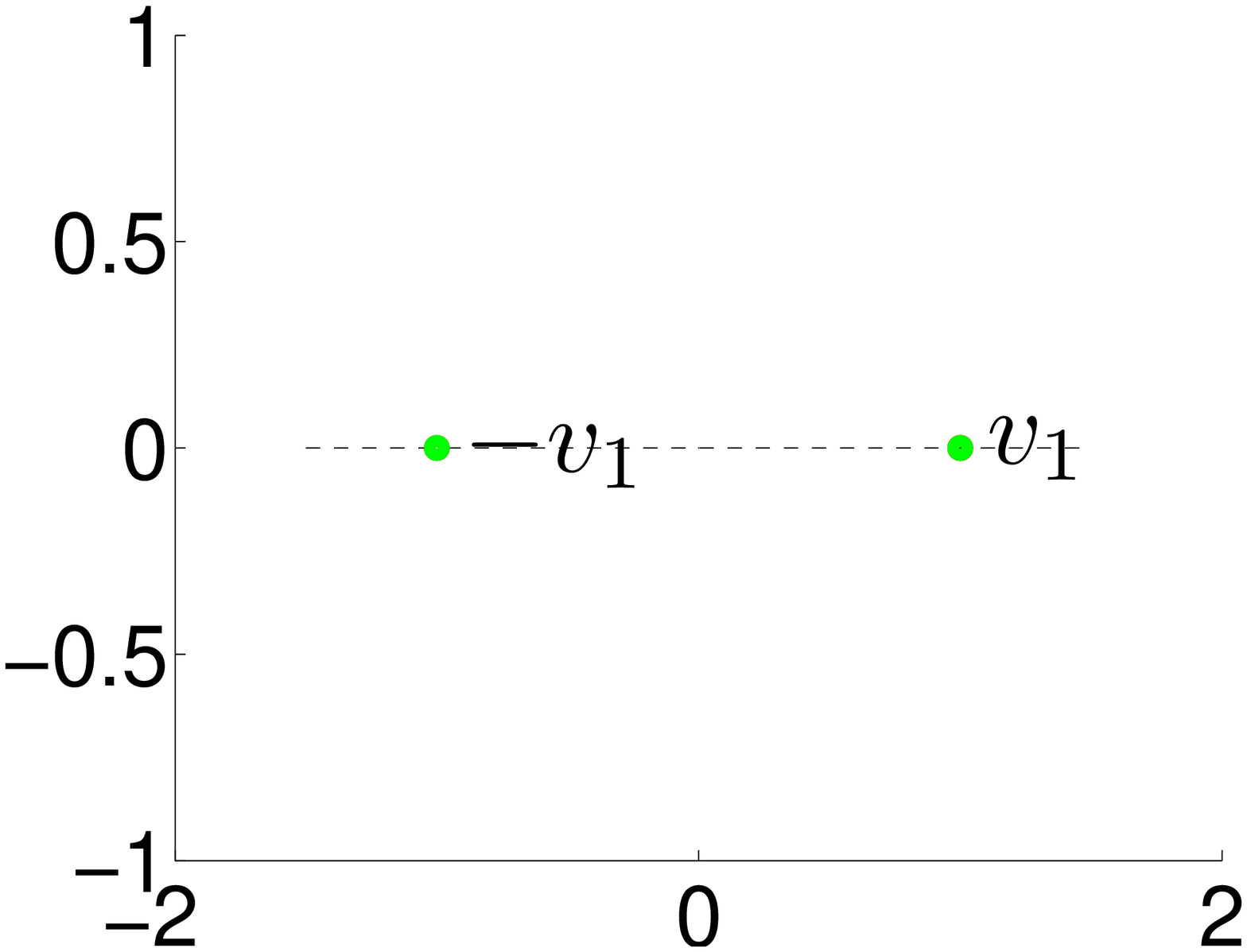}}
\small\centerline{$L_{2}$}\medskip
\end{minipage}
\hfill
\begin{minipage}[b]{0.3\linewidth}
  \centering
  \centerline{\includegraphics[width=\linewidth]{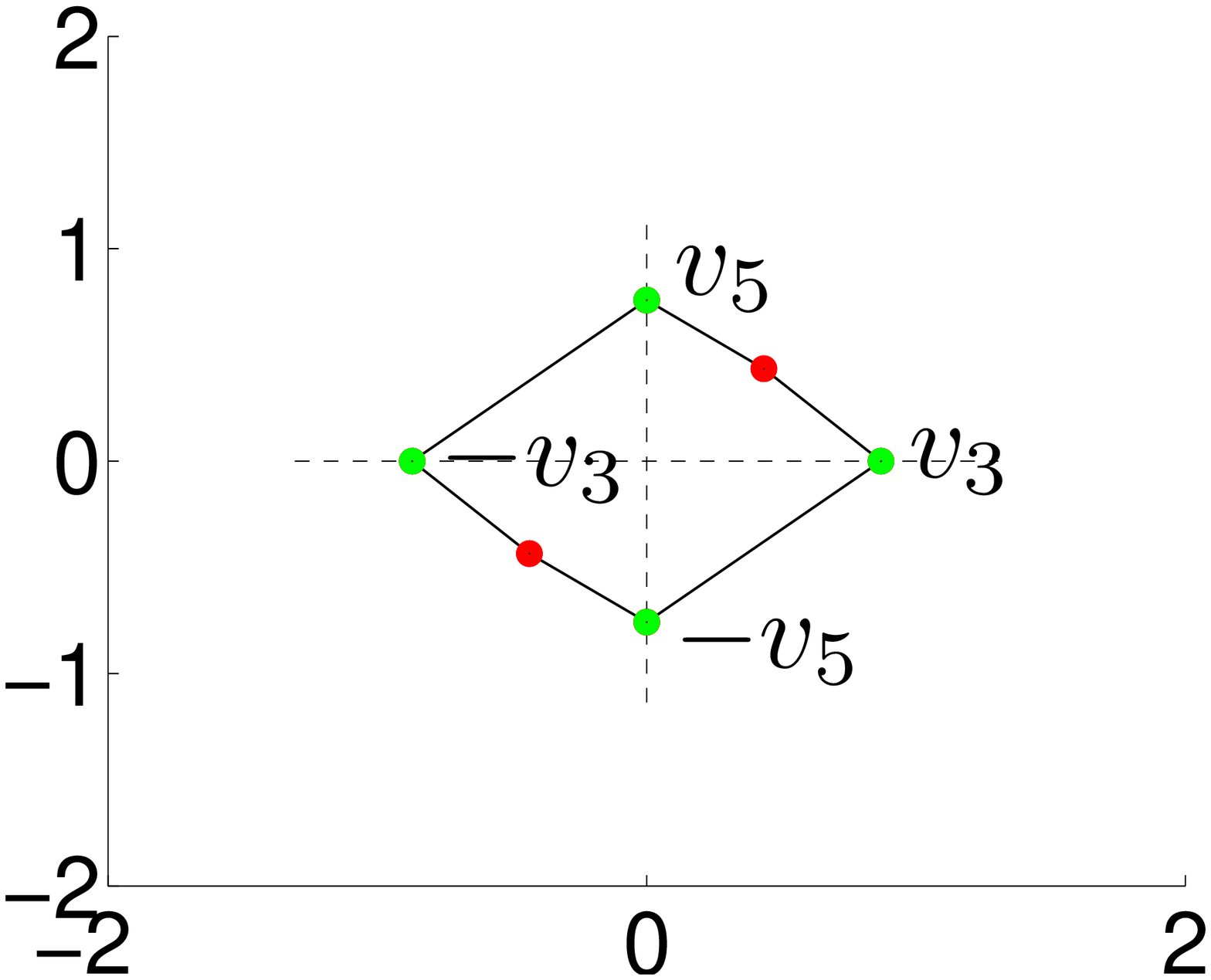}}
\small\centerline{$L_{3}$}\medskip
\end{minipage}
%
\caption{Invariant polytopes for $\xi$}\label{fig:Ex2_Step1}
\end{figure}

After the first step of Algorithm \ref{algo} we have an extremal polytope  multinorm whose unit balls (invariant polytopes) are plotted in Figure \ref{fig:Ex2_Step1}.

So we can conclude that $\Pi$ is an s.m.p. for $\xi$ and $\hr(\xi)=\rho(\Pi)^{\frac{1}{5}}= 1.515717\ldots$.

In this case we cannot compute the unconstrained JSR since not all products of matrices are allowed.

%
%
%

\end{ex}
\begin{ex}\label{ex:three}

In this example we provide some statistics on the performance of the proposed method.
To this aim we consider two fixed graph structures $G^{(1)}$ and $G^{(2)}$, depicted in Figure \ref{fig:GraphG1} and \ref{fig:GraphG2} respectively, we fix the dimension $d$ and  generate $d\times d$ matrices $A_1$ and $A_2$,  as either uniformly distributed random matrices, which we study in Case A, or normally distributed random ones, which is analyzed in Case B.

First of all, we observe that focusing on cases in which the set of matrices contains only two elements may appear restrictive, however, as it has been shown for the computation of the joint spectral radius \cite{CG10}, already in this setting arises the complexity and variety of possible problems which we encounter in the computation of these kind of quantities.

Secondly we point out that, considering the two matrices $A_1$ and $A_2$ as letters of an alphabet, the graphs $G^{(1)}$ and $G^{(2)}$ correspond to two dictionaries with forbidden subword $A_1A_2A_1$ and subwords $A_1A_2A_1$ and $A_1^2$ respectively, as explained in Section \ref{ss-words}, Example \ref{ex10}.

In producing the statistics we make the following assumptions:
\begin{itemize}
  \item We deal with irreducible triplets
  \item The s.m.p. is a dominant product
  \item The s.m.p. has a unique and simple real leading eigenvalue
\end{itemize}



\paragraph{Case A}

We fix the dimension $d\in\{3,\ 5,\ 8,\ 10,\ 15,\ 20\}$, and we produce random matrices whose entries are real numbers drawn from the standard normal distribution by means of the Matlab command \verb"randn". Then we scale  each matrix so that they have spectral radii equal 1. In this way we produce 20 sets $\cA^{(1)}=\left\{ A_1,\ A_2\right\}$ and 20 sets $\cA^{(2)}$. For both triplets $\xi^{(1)}=\left(G^{(1)}, \cL^{(1)}, \cA^{(1)}\right)$ and $\xi^{(2)}=\left(G^{(2)}, \cL^{(2)}, \cA^{(2)}\right)$, where sets $\cL^{(1)}$ and $\cL^{(2)}$ contain four and three spaces of dimension $d$ respectively, we compute first a candidate s.m.p., by means of the method described in \cite{CP15}, then, using the algorithm presented in \cite{GP13} we build an extremal polytopic norm computing the JSR of the unconstrained problem. Afterward, we identify a candidate constrained s.m.p. and we use the proposed method to find an extremal polytopic multinorm. We observe that for the triplet $\xi^{(1)}$ we can identify two vertices in the graph $G^{(1)}$. Hence we consider a set $\cL^{(1)}$ which contains three spaces of dimension $d$.

Mean values statistics on the performance are given in Table \ref{tab:Ex3CaseA_1} and \ref{tab:Ex3CaseA_2} for triplets $\xi^{(1)}$ and $\xi^{(2)}$ respectively.

We point out that, for the triplet $\xi^{(1)}$, we skip all cases in which both $A_1$ and $A_2$ are candidates s.m.p. In fact in this case the hypothesis of a dominant s.m.p. becomes clearly false. Such cases require an ad hoc balancing of the multinorms, as explained in Section \ref{ss-balance}, which is out of the scope of this example.

Furthermore we observe that, while in the unconstrained joint spectral radius problem it is non generic to deal with a reducible set of matrices, in the setting under study it is common to have reducibile triplets. The percentage of cases of reducible triplets out of all the tested cases goes from approximately $20\%$ when $d=3$ to more than $60\%$ when $d=20$, both for triplets $\xi^{(1)}$ and $\xi^{(2)}$.

\begin{table}[ht]
\begin{tabular}{|c|c|c|c|c|c|c|}\hline
  Size   & Steps &  Length SMP & CPU time (s) & Ver. Pol. 1 & Ver. Pol. 2  &  Ver. Pol. 3  \\ \hline
\hline
    3   &    8   &    6   &    5   &    7   &    7   &    8  \\ \hline
    5   &   14   &    5   &   35   &   26   &   25   &   28  \\ \hline
    8   &   18   &    7   &  148   &   83   &   83   &   89  \\ \hline
   10   &   21   &    6   &  349   &  139   &  139   &  151  \\ \hline
   15   &   26   &    5   &  18738   &  470   &  470   &  493  \\ \hline
   20   &   29   &    6   &  16331   &  933   &  938   &  978  \\ \hline
\end{tabular}
\caption{Example \ref{ex:three}, Case A, mean values statistics on triplets $\xi^{(1)}$}\label{table:Test1} 
\label{tab:Ex3CaseA_1}
\end{table}

\begin{table}[ht]
\begin{tabular}{|c|c|c|c|c|c|c|}\hline
  Size   & Steps &  Length SMP & CPU time (s) & Ver. Pol. 1 & Ver. Pol. 2  &  Ver. Pol. 3  \\ \hline
\hline
    3   &    9   &    2   &    3   &    6   &    6   &    6  \\ \hline
    5   &   13   &    2   &   13   &   14   &   14   &   14  \\ \hline
    8   &   15   &    2   &   72   &   34   &   34   &   34  \\ \hline
   10   &   19   &    2   &  315   &   62   &   62   &   62  \\ \hline
   15   &   22   &    2   &  21295   &  195   &  195   &  195  \\ \hline
   20   &   25   &    2   &  51076   &  308   &  308   &  308  \\ \hline
\end{tabular}
\caption{Example \ref{ex:three}, Case A, mean values statistics on triplets $\xi^{(2)}$}\label{table:Test4}  
\label{tab:Ex3CaseA_2}
\end{table}

We observe also that in the case of triplet $\xi^{(2)}$ the number of vertices contained in each extremal polytope tends to be equal for a fixed dimension in many cases, but not all the time, thanks to the cyclic structure of the graph $G^{(2)}$. This is the reason why the mean values of the number of such vertices rounded to the closest integer, which we reported in the last three columns of Table \ref{table:Test4}, tend to be the same for a fixed dimension. In the case of triplet $\xi^{(1)}$, instead, the mean values of the number of vertices contained in each extremal polytope tends to be similar each other due to the averaging.

Finally we point out that during the numerical tests (which were performed with MATLAB R2011a installed on a 64--bit Windows 7 Professional computer equipped with a core i3-3227U 
processor and 8GB RAM) we run the proposed algorithm in some cases did not compute an extremal polytopic multinorm. Either because the initial guess for the candidate s.m.p. was wrong or because the number of iterations become bigger than a maximal number we set a priori, in our computation is set to 40. The percentage of such cases is ranging from $2\%$ to $10\%$, as we increase the dimension of the matrices, for the triplets $\xi^{(1)}$, whereas is around $1\%$ for any dimension of the matrices in the triplets $\xi^{(2)}$.
In these cases we end up having anyway an interval of approximation for the constrained j.s.r. given by the spectral radius of the candidate s.m.p. and the maximal value of polytopic multinorm of the matrices in the set under study. The length of such intervals range from a maximum value of order $10^{-3}$ to a minimum of order $10^{-10}$.

\paragraph{Case B}

This time we consider $d\in\{5,\ 10,\ 20,\ 50,\ 100\}$, and by means of the Matlab command \verb"rand" we produce matrices whose entries are uniformly distributed random real numbers in the interval $(0,\ 1)$. Then we scale each of them so that they have spectral radii equal to one.  We repeat this process to produce 20 sets $\cA^{(1)}$ and 20 sets $\cA^{(2)}$. As for the previous case, for both triplets $\xi^{(1)}$ and $\xi^{(2)}$, we compute an s.m.p. and an extremal polytopic norm for the unconstrained problem, using both the method described in \cite{CP15} and the algorithm presented in \cite{GP13}. Afterward, we identify a candidate s.m.p. and we use the proposed method to find an extremal polytopic multinorm.

Mean values statistics for this second case are given in Table \ref{tab:Ex3CaseB_1} and \ref{tab:Ex3CaseB_2} for triplets $\xi^{(1)}$ and $\xi^{(2)}$ respectively.

As for case A, for the triplet $\xi^{(1)}$ we skip all cases in which both $A_1$ and $A_2$ are candidates s.m.p.

We point out also that in this case, since we produce matrices whose entries are uniformly distributed random number in the interval $(0,\ 1)$ as $d$ increases, before the scaling, the average of each row and column of such matrices tends to $0.5$ therefore the vector of all ones becomes an eigenvector corresponding to the eigenvalue 0.5. After the  scaling the spectrum of such a matrix does contain an eigenvalue 1 and, the bigger the dimension $d$, the smaller all the other eigenvalues are going to be.

\begin{table}[ht]
\begin{tabular}{|c|c|c|c|c|c|c|}\hline
  Size   & Steps &  Length SMP & CPU time (s) & Ver. Pol. 1 & Ver. Pol. 2  &  Ver. Pol. 3  \\ \hline
\hline
    5   &    2   &    3   &  0.0764   &    2   &    1   &    1  \\ \hline
   10   &    3   &    3   &  0.1922   &    3   &    2   &    2  \\ \hline
   20   &    2   &    3   &  0.2092   &    2   &    2   &    2  \\ \hline
   50   &    1   &    3   &  0.2219   &    1   &    1   &    1  \\ \hline
  100   &    1   &    3   &  0.1079   &    1   &    1   &    1  \\ \hline
\end{tabular}
\caption{Example \ref{ex:three}, Case B, mean values statistics on triplets $\xi^{(1)}$}\label{table:Test3}  
\label{tab:Ex3CaseB_1}
\end{table}

\begin{table}[ht]
\begin{tabular}{|c|c|c|c|c|c|c|}\hline
  Size   & Steps &  Length SMP & CPU time (s) & Ver. Pol. 1 & Ver. Pol. 2  &  Ver. Pol. 3  \\ \hline
\hline
    5   &    3   &    2   &  0.0708   &    1   &    1   &    1  \\ \hline
   10   &    3   &    2   &  0.0726   &    1   &    1   &    2  \\ \hline
   20   &    3   &    2   &  0.0642   &    1   &    1   &    2  \\ \hline
   50   &    2   &    2   &  0.2083   &    1   &    1   &    1  \\ \hline
  100   &    2   &    2   &  0.0676   &    1   &    1   &    1  \\ \hline
\end{tabular}
\caption{Example \ref{ex:three}, Case B, mean values statistics on triplets $\xi^{(2)}$}\label{table:Test2}  
\label{tab:Ex3CaseB_2}
\end{table}
\end{ex}

Finally we observe that also in these tests the proposed algorithm in some cases did not compute an extremal polytopic multinorm. The percentage of such cases is ranging from $0\%$ to almost $30\%$, roughly as we decrease the dimension of the matrices, for the triplets $\xi^{(1)}$, whereas is always $0\%$ for any dimension of the matrices in the triplets $\xi^{(2)}$.
Also in these cases we end up having a length of such intervals which range from a maximum value of order $10^{-3}$ to a minimum of order $10^{-10}$.

\section{Construction of a polytope Barabanov multinorms}\label{s-bar}

By Theorem~\ref{th60} Algorithm~1 terminates within finite time if and only if
the triplet~$\xi$ possesses a dominant product with a unique and simple leading eigenvalue. Moreover, the algorithm
produces an invariant family of polytopes $\cP = \{P_i\}_{i=1}^n$. The corresponding Minkowski multinorm
$\|\cdot \|_{\cP} = \{\|\cdot \|_{P_i}\}_{i=1}^n$ is extremal, i.e., $\|A_{ji}\|_{\cP} \le \hr(\xi)$ for all
$A_{ji} \in \cA$. However, it is not necessarily invariant. Nevertheless, the triplet $\xi$ does have a polytope  invariant norm, which can be constructively found. It turns out that Algorithm~1  applied to the dual triplet~$\xi^*$
also converges within finite time producing an invariant family of polytopes. The dual of those polytopes generate
Barabanov norm for~$\xi$.  This is guaranteed by Theorem~\ref{th65} below. To formulate it we need to define first
the dual triplet.
\begin{defi}\label{d80}
A triplet $\xi^* = (G^*, \cL^*, \cA^*)$ is dual to a triplet $\xi = (G, \cL, \cA)$ if

1) the multigraph $G^*$ has the same vertices as~$G^*$ and the reverses edges;

2) $L_i^*$ is a dual space for $L_i\, , \, i = 1, \ldots , n$;

3) $\cA_{\, ij}^*$ consists of operators adjoint to the operators of the family~$\cA_{ji}$.
\end{defi}
Since all the spaces $L_i$ are finite-dimensional, we identify $L_i^*$ and $L_i$.
Clearly, every path $\alpha: i_1\to \cdots \to i_n$ along $G$ corresponds to the
reverse path~$\alpha^*: i_n \to \cdots \to i_1$ on $G^*$ and the corresponding matrix products are
adjoint to each other. Hence, $\hr(\xi^*) = \hr(\xi)$. For a given polytope $P$
we denote  by $\cV(P)$ the set of its vertices.
\begin{theorem}\label{th65}
If Algorithm~1 applied for  a triplet $\xi$ converges within finite time,
then it also does for the dual triplet $\xi^*$. Moreover, the multinorm
\begin{equation}\label{dual-bar}
\|x\|_i \ = \ \max_{u \in \cV(P_i')}  \, \bigl(u\, , \, x   \bigr) \, , \qquad i=1 ,\ldots , n\, ,
\end{equation}
where $\{P_i'\}_{i=1}^n$ is the invariant family of polytopes  produced by Algorithm~1 for $\xi^*$, is
a polytope Barabanov multinorm for
$\xi$.
\end{theorem}
{\tt Proof.} Assume $\hr(\xi) = 1$, and hence $\hr(\xi^*) = 1$. By Theorem~\ref{th10},
if Algorithm~1 converges within finite time, then the chosen candidate constrained s.m.p. product~$\Pi \in \cC(G)$ is dominant
and has a unique simple eigenvalue. Therefore, the dual product $\Pi^*$ (the product of adjoint  operators in the inverse order)
possess the same properties for the family~$\xi^*$. Hence, Algorithm~1 applied for $\xi^*$ with the candidate product
$\Pi^*$ converges as well producing some invariant family of polytopes $P_i' \subset L_i^*, \, i = 1, \ldots , n$.
By the construction of the algorithm, for each~$i$, the polytope $P_i'$ coincides with the convex hull
of images $A_{ij}'P_j', \, A_{ij}' \in \cA_{ij}^*$,  taken over all incoming edges $l_{ij}^* \subset \ell_{ij}^*, \, i = 1, \ldots , n$.
On the other hand, each operator $A_{ij}'$ from the family~$\cA_{ij}^*$ is adjoint
to the corresponding operator $A_{ji} \in \cA_{ji}$, i.e., $A_{ij}' = A_{ji}^*$.
Thus,
$$
P_i' \ = \ {\rm co} \, \Bigl\{ A_{\, ij}'\, P_j' \quad \Bigl| \ A_{ij}' \in \cA_{ij}^*, \,
j = 1, \ldots , n\, \Bigr\} \ = \ {\rm co} \, \Bigl\{ A_{ji}^*\, P_j' \quad \Bigl| \ A_{ji} \in \cA_{ji}, \,
j = 1, \ldots , n\, \Bigr\}\, .
$$
Consequently, the sets
$$
\cV(P_i') \qquad \mbox{and} \quad   \bigcup\limits_{A_{ji} \in \cA_{ji}, \,
j = 1, \ldots , n} \ A_{ji}^*\, \bigl( \cV(P_j')\, \bigr)\, ,
$$
have the same convex hulls. Therefore, the multinorm $\{\|\cdot \|_i\}_{i=1}^n$defined by~(\ref{dual-bar})  satisfies
$$
\|x\|_i \quad = \quad \max \, \Bigl\{   \bigl(v'\, , \, x   \bigr) \ \Bigl| \
v' \in \cV(P_i') \Bigr\} \quad  = \
$$
$$
\max \, \Bigl\{   \bigl(v'\, , \, x   \bigr) \ \Bigl| \
v' \, \in \, A_{ji}^*\, \bigl(\,  \cV(P_j')\, \bigr), \
A_{ji} \in \cA_{ji}, \,
j = 1, \ldots , n\Bigr\}\ = \
$$
$$
\max_{A_{ji} \in \cA_{ji}, \, j = 1, \ldots , n}\
\max \, \Bigl\{  \, \bigl(A_{ji}^*u'\, , \, x   \bigr) \, , \
u' \in \cV(P_j')\, \Bigr\} \ = \
$$
$$
\max_{A_{ji} \in \cA_{ji}, \, j = 1, \ldots , n}\
\max \, \Bigl\{  \, \bigl(w'\, , \, A_{ji}x   \bigr) \, , \
w' \in \cV(P_j')\, \Bigr\}
 \ = \ \max_{A_{ji} \in \cA_{ji}, \, j = 1, \ldots , n}\
\bigl\| A_{ji}x\bigr\|_j \,.
$$
Thus, $\|x\|_i \ = \max\limits_{A_{ji} \in \cA_{ji},  j = 1, \ldots , n}\
\bigl\| A_{ji}x\bigr\|_j\, $, hence $\{\|\cdot\|_i\}_{i=1}^n$ is an invariant multinorm.

  {\hfill $\Box$}
\medskip

\begin{cor}\label{c20}
If a triplet possesses a dominant product with a unique and simple leading eigenvalue,
then it possesses a polytope Barabanov multinorm.
\end{cor}
The single-space version а Theorem~\ref{th65} was established in~\cite{GZ15}.

\section{Generalizations and special cases}\label{s-special}

\subsection{Positive systems}\label{ss-pos}

A triplet $\xi$ is called {\em positive} if one can introduce a basis in each space $L_i\, , \, i = 1, \ldots , n$,
such that all operators from $\cA$ are written by nonnegative matrices. Dealing with positive systems we
will assume that such a collection of bases is fixed and identify the operators with their matrices.
It was observed in the literature that some of methods of computing of the joint spectral radius
in the classical
(unconstrained) case work more efficiently for nonnegative matrices. For instance, the Invariant polytope algorithm works
effectively for nonnegative matrices of dimensions $d = 100$ and higher (see examples and statistics in~\cite{GP13} along with the discussion of this phenomenon). That is why the positive systems deserve a special analysis. First of all,
the irreducibility assumption can be relaxed to positive irreducibility. This notion is directly extended
from the single-space case, where it is well-known. To define it we need some extra notation.

A coordinate subspace of $\re^d$ is a subspace spanned by several basis vectors.
A triplet
$\xi' = (G, \cL', \cA')$ is {\em positively embedded} in $\xi$, if both these triplets are positive,
$L_i'$ is a coordinate subspace of~$L_i$
for each $i$, and every operator $A_{ji}' = A_{ji}|_{L_i'}$ maps $L_i'$ to $L_j'$, whenever
$l_{ji} \in G$. The embedding is strict if $L_i'$ is a proper subspace of $L_i$ at least for one~$i$.
\begin{defi}\label{d90}
A triplet $\xi = (G, \cL, \cA)$ is positively reducible if it is positive and has a strictly embedded triplet. Otherwise, it is called positively irreducible.
\end{defi}
The factorization of a positively reducible triplet to positive triplets of smaller dimensions
is realized in precisely the same way as in Section~6. Thus, analysing positive triplet we can concentrate on the
  positively irreducible case.

The following analogue of Theorem~\ref{th10} holds for positive triples.
\begin{theorem}\label{th70}
An positively irreducible triplet is non-defective.
\end{theorem}
{\tt Proof} actually repeats the proof of Theorem~\ref{th10} with several different points.
First, we denote by $S_i$ not the unit sphere in $L_i$ but its intersection with the positive orthant.
The sets $U_{i,k}$ are defined in the same way, the case
$\cup_{k \in \n} U_{i,k} \, = \, S_i$ is considered in the same way as for Theorem~\ref{th10}.
In the converse case,  there exists $z \in S_q$ which does not belong to any of the sets $U_{q, k}$, i.e.,
for every path $\alpha$ staring at the vertex $q$ we have  $\|\Pi_{\, \alpha} \, x\| \le 2$.
Hence, the point~$z$ has a bounded orbit.
For every $i$, denote by $M_i$ the set of nonnegative points from
$L_i$ that have bounded orbits.
The linear span of each $M_i$ is a linear subspace of $L_i$. Moreover, it is a coordinate subspace.
 Otherwise, there is $a \in M_i, a > 0$. For every $x \in L_i, x\ge 0$, there is a number $\lambda > 0$
 such that $\lambda x \le a$. Hence the orbit of $x$ is bounded as well. Thus, $M_i$ coincides with the positive orthant of~$L_i$.
 Now the positive irreducibility implies that $M_j$ is the positive orthant of $L_j$ for all $j$ which implies non-defectivity.
 The remainder of the proof is the same as for Theorem~\ref{th10}.

   {\hfill $\Box$}
\medskip

A norm is called {\em monotone} if $\|x\| \ge \|y\|$, whenever $x \ge y \ge 0$. A monotone multinorm is a connection
of monotone norms. The following theorem sharpens
Theorem~\ref{th20} for the case of positive systems. Its single-space version was established in~\cite{GP13}.
\begin{theorem}\label{th80}
A positively  irreducible triplet possesses a monotone invariant multinorm.
\end{theorem}
{\tt Proof} is actually the same as for Theorem~\ref{th20}, with two modifications.
First of all, we take an arbitrary initial monotone multinorm. Then the function
$f(x) = \limsup\limits_{|\alpha| \to \infty} \|\Pi_{\, \alpha}\, x\|$ defined for nonnegative $x$
is obviously monotone.
To show  that $f(x) > 0$ for all $x\ne 0$, we note that
 $f(x) = 0$ implies $\lim\limits_{|\alpha| \to \infty} \|\Pi_{\, \alpha}\, x\| = 0$.
Let $M_i$ be the set of nonegative points $x \in L_i$ satisfying this equality.
The linear span of $M_i$ is a linear subspace of $L_i$. The collection $\{M_i\}_{i=1}^n$
possesses the invariance property. If some $M_i$ contains a strictly positive point~$a$, then
it coincides with the positive orthant of~$L_i$.
Indeed, an arbitrary point $x \in L_i, x\ge 0$, satisfies $\lambda x \le a$ for some $\lambda > 0$, hence
 $\lim\limits_{|\alpha| \to \infty} \|\Pi_{\, \alpha}\, x\| = 0$, and so $x \in M_i$.
 Thus, the linear span of $M_i$ is a coordinate subspace of $L_i$. Then we repeat the proof of Theorem~\ref{th20}
 and come to the contradiction with the positive irreducibility.

   {\hfill $\Box$}
\medskip

Algorithm~1 is modified for positive case as follows. First of all, we omit all
extra assumptions on the candidate constrained s.m.p. $\tilde \Pi_{\alpha}$. Indeed, since this matrix is now
positive, then by the Perron-Frobenius theorem its leading eigenvalue is real and positive.
Second, in the $k$th iteration of Algorithm~1 for the general case, we check whether $\tilde A_{ji}v$
is an interior point of the set ${\rm absco} \, \bigl(\cV_{j}^{(k-1)}\bigr)$. Now we replace this set by
${\rm co}_{-} \, \bigl(\cV_{j}^{(k-1)}\bigr)$. Thus, the version of Algorithm~1 for a positive system
constructs a collection of monotone polytopes $\{P_i\}_{i=1}^n$ instead of symmetric polytopes.
This version works much faster than the algorithm for general matrices, because
for any set $M \subset \re^d_+$, the set ${\rm co}_{-}(M)$ is bigger than ${\rm absco}(M)$
in the positive orthant (usually it is much bigger). Hence, in the positive case, each iteration of Algorithm~1
sorts our more vertices than in general case.

A complete analogue of Theorem~\ref{th65} holds for positive case and gives a monotone  invariant  polytope norm.

\subsection{Stabilizability, lower spectral radius and antinorms}\label{ss-lower}

The notion of stailizability well-known for dynamical systems is also extended for
systems on graphs in a direct manner.
 \begin{defi}\label{d100}
The system $\xi$ is called stabilizable  if it has at least one infinite path such that all corresponding  trajectories
converge to zero.
 \end{defi}
 If the multigraph $G$ is strongly connected  and has an edge associated to a zero operator, then
 every path going through this edge produce vanishing trajectories. Hence, the system is stabilizable in this case.
We see that analysing stabilizability one cannot identify an empty edge with a zero operator,
otherwise, all systems will be stabilizable.
In the sequel we assume that $G$ has at least one infinite path, i.e., has a cycle.

 See~\cite{BS08, D95, FV12, LA09, SDP08} for properties of stabilizable systems in the classical (unconstrained) case
 and for criteria of stabilizability.
 Most of those properties and criteria are extended to systems on graphs. In particular,
 the stabilizability is expressed  in terms of the lower spectral radius.
  \begin{defi}\label{d110}
 The lower spectral radius (LSR) of a triplet $\xi = (G, \cL, \cA)$ is
 \begin{equation}\label{lsr}
 \check \hr(\xi) \ = \ \lim_{k \to \infty}\, \min_{|\alpha| = k} \|\Pi_{\alpha}\|^{\, 1/k}\, .
 \end{equation}
 \end{defi}
 \begin{prop}\label{p60}
The system is stabilizable if and only if $\check \hr < 1$.
\end{prop}
This fact in the single-space case is well-known and originated with~\cite{D95}.
The proof for general systems on graphs is different.
\smallskip

{\tt Proof}. The sufficiency is obvious. To prove the necessity we assume that
there is a path~$\alpha$ such that $\|P_{\alpha_j}x_0\| \to 0$ as $j\to \infty$
for every $x\in L_{i_0}$. Here $\alpha_j$ denotes the prefix of $\alpha$ of length~$j$.
There is at least one vertex through which $\alpha$ passes infinitely many times.
Without loss of generality we assume that this is the starting vertex~$x_0$.
For each basis vector  $e_s$ of the space $L_{i_0}$ we have
$\|P_{\alpha_j}e_s\| \to 0$,  hence $\|P_{\alpha_j}\| \to 0$ as $j\to \infty$.
Therefore, there exist arbitrarily long closed paths $\alpha_j$ such that $\|P_{\alpha_j}\| < 1$.
 Take one of them and denote $q = \bigl[P_{\alpha_j}\bigr]^{1/|\alpha_j|} < 1$. Since $\alpha_j \in \cC(G)$ it follows that all powers $(\alpha_j)^k$ are well-defined. Consequently,  $\, \check \hr \le \bigl[\, \rho(P_{\alpha_j}^k)\, \bigr]^{1/k|\alpha_j|} = q < 1$.

   {\hfill $\Box$}
\medskip

From the computational point of view, the lower spectral radius is still worse than JSR.
For instance, it is, in general, a discontinuous function of matrices. Nevertheless,
some algorithms of approximate computing of LSR for the classical case (unconstrained systems)
exist~\cite{PJB}. In~\cite{GP13} an algorithm of exact computation
of LSR for positive systems (i.e., with all matrices from $\cA$ nonnegative) was presented. Under some mild assumptions, it
gives the exact value of LSR for a vast majority of  families  of nonnegative matrices. The idea is similar to Algorithm~1, with some modificaions. First of all, it uses the notion of {\em antinorm} instead of norm.
An antinorm is a nonnegative, nontrivial (not identical zero), positively homogeneous, and  concave function defined on the positive orthant.
In a sense, this is a ``concave norm'' on the positive orthant. A collection of antinorms
$\{f_i\}_{i=1}^n$ can be called ``multi-antinorm'', but we drop the prefix and call it just antinorm.
The notions of extremal and invariant antinorms are the same as for norms (Definition~\ref{d60}), with the replacement of $\max$ by $\min$.
See~\cite{GP13} for the existence results for extremal and invariant antinorms in the single-space case.
In particular, every positive system has a monotone extremal antynorm~\cite[Theorem~5]{GP13}. This result is extended to general systems on graphs without any change.

The LSR version of Algorithm~1 constructs an extremal polytope antinorm. We begin with
exhaustion of all closed paths of lengths bounded by a given number~$l_0$ and find the
candidate $\Pi$ for {\em spectrum minimizing product}, i.e., a product for which the value
${\rho_{\alpha} = \|\Pi_{\alpha}\|^{1/|\alpha|}}$
is minimal. Then we repeat the routine of Algorithm~1, replacing the symmetrized convex hull ${\rm absco}(M)$
by ${\rm co}_+(M)$. So, the LSR version of Algorithm~1 deals with infinite polytopes $P_j^{(k)} =  {\rm co}_+\bigl(\cV_j^{(k)}\bigr)$ (see Introduction for the definition).
The algorithm halts when no new vertices appear. In this case the constructed infinite polytopes $\bigl\{ P_j^{(k)}\bigr\}_{j=1}^n$ in invariant and generate
an extremal antinorm, and the lower spectral radius is found: $\check \rho(\xi) = \rho_{\alpha}$.

\subsection{The case of several spectral maximizing  products. Balancing method}\label{ss-balance}

In applications of the Invariant polytope algorithm in the single-space case,
we sometimes meet the following trouble: the constrained s.m.p. candidate product is not unique.
There are several products $\Pi_1 , \ldots , \Pi_r$ (not powers or cyclic  permutations of each other)
 that have the same maximal value $\rho(\Pi_j)^{1/|\alpha_j|}, \, j = 1, \ldots , r$. In this case,
 the system~$\xi$ does not have a dominant product because the dominant product must be  (by definition!)  unique!
 Hence,  by Theorem~\ref{th60},  Algorithm~1 cannot converge. Of course, the situation when spectral radii of some
 products coincide is not generic. Nevertheless, it is sometimes emerges in applications, when the operators from~$\cA$
 have some symmetries or relations to each other. In the classical (unconstrained) case the method of {\em balancing}
 is presented in~\cite{GP16} to extend Algorithm~1 to this case. The idea is to multiply the leading eigenvectors
 of the products $\Pi_1 , \ldots , \Pi_r$  by certain positive coefficients $a_1, \ldots , a_n$
 respectively, and then apply Algorithm~1 starting with all those multiplied initial vectors simultaneously. The coefficients
 $\{a_j\}_{j=1}^r$ realizing the balancing of eigenvectors can be  found as solutions of an optimization problem.
 See~\cite{GP16} for details, examples, and for the criterion of convergence.  This method is extended to the systems on graphs without any change. The optimization problem and the way of computing of the balancing coefficients remain the same.

\section{Applications}\label{s-appl}

We elaborated  a factorization procedure of an arbitrary system on graph to several irreducible systems of smaller dimensions, proved the theorem of existence of invariant multinorm for irreducible system, and presented  the Invariant polytope algorithm (Algorithm~1)  for computing that antinorm and the value of JSR. Now let us discuss possible areas of applications of our results.

\subsection{Constrained linear switching systems}\label{ss-lss}

One of the main applications  is in the study of constrained linear switching systems. The constraints are usually  generated by graphs, regular languages, or finite automata.
Philippe and  Jungers in~\cite{PhJ1} considered  systems that have finitely many stages; a transfer from one stage to another is either impossible or is realized by a given linear operator. Thus, we have a graph and a system of linear operators corresponding to its edges. This is a special case
of a triplet $\xi(G, \cL\, \cA)$, with a graph (not mutigraph)~$G$ and with the same linear spaces
 $L_1 = \cdots = L_n = \re^d$ in all its vertices. Sufficient conditions
for the stability of such systems presented in~\cite{PhJ1} establish the stability by constructing an ellipsoidal
multinorm, in which all operators from~$\cA$ become contractions. As a rule, those conditions are not necessary. The existence of such an antinorm can be verified
by solving the corresponding s.d.p. problem. Algorithm~1 (Section~\ref{s-algor}) makes it possible to find a precise value of JSR, which gives the criterion of stability and the corresponding Barabanov multinorm.

Special cases of the constrained switching systems from~\cite{PhJ1} are Markovian systems, when
each operator~$A_i \in \cA$ has ist own list of operators from~$\cA$ that are allowed to follow~$A_i$~\cite{D, K, WRDV}.
In~\cite{SFS} this concept was extended to arbitrary graphs and interpreted as ``systems with memory''.
Another extension was studied in~\cite{PhJ2} as systems with  switching sequences defined by a regular
language generated by a finite automata, see also~\cite{OPJ}. We considered Markovian systems and its generalization, with an arbitrary dictionary of prohibited words, in subsections~\ref{ss-markov} and~\ref{ss-words}.

The results of subsections~\ref{ss-pos} and \ref{ss-lower} are also directly applicable to the
constrained switching systems. The modification of Algorithm~1 to nonnegative matrices (\S~\ref{ss-pos}) makes it more efficient for positive systems. The results on the lower spectral radius and invariant antinorms~(\S~\ref{ss-lower}) give a criterion of stabilizability of constrained switching systems.

 \subsection{Fractals}\label{ss-fract}

We consider a nonlinear generalization of our construction,  to a set of arbitrary metric spaces
 and arbitrary maps between them. Thus, we have a triplet $\xi = (G, \cM, \cF)$ with a multigraph~$G$
 with vertices~$\{g_i\}_{i=1}^n$,
 complete metric spaces $\{M_i\}_{i=1}^n$ associated to them,  and finite sets $\cF_{ji}$ of maps $F_{ji}: \, M_i \to M_j$.
 A collection of compact sets $\cK = \{K_i\}_{i=1}^n, \, K_i \subset M_i, \, i = 1, \ldots, n$,
  is called a fractal if it possesses the property
  $$
  K_j \ = \ \bigcup\limits_{F_{ji} \in \cF_{ji}, \, i=1, \ldots , n}\ F_{ji}\, K_i\, , \qquad j = 1, \ldots, n\, .
   $$
 Thus, each set $K_j$ is the union of images of sets $\{K_i\}_{i=1}^n$ by the operators
 associated to all incoming edges of the vertex~$g_j$. This is a straightforward  generalization of the classical concept  of fractal by J.Hutchinson~\cite{H}
 in case of a single metric space. Similarly to the classical situation, one can show that if all
 maps $F_{ji}$ are contractions, then the system has a unique fractal. If all the spaces $M_i$ and all the maps $F_{ji}$ are affine, then we have an {\em affine system (triplet)~$\xi$}. In this case the fractal (if it exists) is called affine.
  In most cases the contraction property of an affine system is too restrictive.
  Indeed, an affine map may be a non-contraction, but become contraction in a different norm introduced in the space.
 That is why, a more general existence and uniqueness result is the following:
 \smallskip

 {\em A system possesses a unique affine  fractal  provided there exists
 a multinorm in the affine spaces $M_i$ in which all the maps $F_{ji}$ become contractions, i.e.,
 the operator norm of linear parts of all those maps are smaller than one.}
  \smallskip

  In view of Proposition~\ref{p6}, we obtain the following tight and affinely-invariant sufficient condition:
 \begin{prop}\label{p70}
 An affine system possesses a unique fractal whenever the joint spectral radius of the associated  linear system
 is smaller than one.
  \end{prop}
 Thus, the stability of a linear system implies the existence of a fractal for the affine system.

 \subsection{Attractors of hyperbolic dynamical systems}\label{ss-dyn}

Systems of maps of metric spaces along edges of a graph arise naturally  in the study of attractors of dynamical systems.
A hyperbolic dynamical system satisfying {\em A Axiom} admits the so-called {\em Markov partition} to sets $\{M_i\}_{i=1}^n$ with a collection of diffeomorphisms $F_{ji}: M_i \to M_j$~(see~\cite{B, KH}). An attractor of that dynamical system is defined similarly to the definition of a fractal above. See~\cite{EB, HLZ} for the analysis of general hyperbolic attractors and~\cite{L} for extensions  to Lipschitz maps. If all those diffeomorphisms are affine, then one can use Proposition~\ref{p70} to prove the existence of an attractor. Applying Algorithm~1 one can construct a polytope multinorm
in which all $F_{ji}$ are contractions (in case $\rho < 1$). Moreover, the corresponding operator norms $\|F_{ji}\|$
can be used to estimate the dimension of the attractor~\cite{L}.

\subsection{Application to numerical ODEs}

An important application of our approach is the possibility to determine sharp bounds for the stepsize ratio
in the zero stability analysis of $k$--step BDF--formulas for the numerical approximation of initial value problems for ODEs, on grids with variable stepsize.

\subsection*{$k$--step BDF--method with variable stepsize} \ \\

For the initial value problem
$$
y'(t)=f(t,y(t)), \ \ \ y(t_0)=y_0,
$$
consider the grid $\Delta=\{t_0, t_1, \ldots , t_n, \ldots \}$ characterized by
the stepsizes $h_j=t_{j+1}-t_j$, $j=0,1, \ldots $, and the $k$--step BDF--method
defined by
$$
y_{j+k}=\alpha_{j,k-1}y_{j+k-1}+\ldots+\alpha_{j,1}y_{j+1}+\alpha_{j,0}y_j+
h_j\beta_jf(t_{j+k},y_{j+k}).
$$
The coefficients $\alpha_{j,s}$, $s=0,1,k-1$, depend on the step ratios
\begin{equation}
\omega_{j,1}=\frac{h_{j+1}}{h_j},\ldots,
\omega_{j,k-1}=\frac{h_{j+k-1}}{h_{j+k-2}}.
\label{2}
\end{equation}

It is well--known (see e.g. \cite{HNW}) that zero--stability is equivalent to the uniform boundedness
of all the solutions of the homogeneous linear difference equation
\begin{equation}
u_{j+k}=\alpha_{j,k-1}u_{j+k-1}+\ldots+\alpha_{j,1}u_{j+1}+\alpha_{j,0}u_j,
\label{3}
\end{equation}
which, in turn, is equivalent to the uniform boundedness of the sequence of products
$$
A_{\nu}\ldots A_0, \ \ \ \nu=0,1, \ldots ,
$$
where the $k\times k$--matrices
$$
A_j=\left( \begin{array}{ccccc}
\alpha_{j,k-1} & \alpha_{j,k-2} & \hdots & \alpha_{j,1} & \alpha_{j,0} \\
1 & 0 & \hdots & 0 & 0 \\
0 & 1 & 0 & \hdots & 0 \\
0 & 0 & 1 & \hdots & 0 \\
0 & \hdots & 0 & 1 & 0
\end{array} \right), \ \ \ j=0,1, \ldots ,
$$
are the companion matrices associated to the difference equation \eqref{3}.

The presence of the common eigenvector $[1,1,\ldots,1]^T$, corresponding
to the common eigenvalue $\lambda=1$, makes it possible to reduce the
dimension of one unit by using a suitable similarity transformation, 
$$
T^{-1} A_j T =\left( \begin{array}{ccccc}
\gamma_{j,k-2} & \gamma_{j,k-3} & \hdots & \gamma_{j,0} & 0 \\
1 & 0 & \hdots & 0 & 0 \\
0 & 1 & 0 & \hdots & 0 \\
0 & 0 & 1 & \hdots & 0 \\
0 & \hdots & 0 & 1 & 1
\end{array} \right), \ \ \ j=0,1, \ldots .
$$

Letting
$$
C_j=\left( \begin{array}{cccc}
\gamma_{j,k-2} & \gamma_{j,k-3} & \hdots & \gamma_{j,0} \\
1 & 0 & \hdots & 0  \\
0 & 1 & 0 & \hdots \\
0 & \hdots & 1 & 0
\end{array} \right),
$$
where $\gamma_{j,0}, \ldots, \gamma_{j,k-2}$ are certain coefficients
depending on $\alpha_{j,0},\ldots,\alpha_{j,k-1}$,
the zero--stability of the BDF--method is guaranteed if
$$
\lim_{\ell\rightarrow \infty}C_{\ell}\ldots C_0=O.
$$

It can be shown that the matrix $C_j$ can be written in terms of the ratios
$$
\omega_{j,1},\ldots,\omega_{j,k-1}
$$
so that a natural question is that of finding sharp bounds for such ratios,
which guarantee zero--stability of the formula.

This would determine an infinite dimensional family, whose analysis is difficult
in general. For this reason, in order to simplify the analysis,
we allow here only a finite number of values for the ratios $\omega_{j,t}$,
say -- for example --
\begin{equation}
\omega_{j,t} \in \{ \theta,1,1/\theta\}
\label{eq:ex33}
\end{equation}
(with $\theta>1$).
This is a sort of discretization of the set of possible stepsize variations.

Our aim is to prove that for a certain given $\theta$,
every sequence of matrices $\{C_{\theta}^{(j)}\}_{j\geq 0}$ of the family
\begin{equation}
{\cal F}_{\theta}=\{C(\omega_1,\ldots,\omega_{k-1})\}_{\omega_1,\ldots,\omega_{k-1}},
\label{15}
\end{equation}
with $\omega_1,\ldots,\omega_{k-1} \in \{ \theta,1,1/\theta \}$, and under the constraint:
\begin{equation}
\omega_1^{(\ell+1)} = \omega_2^{(\ell)}, \qquad \omega_2^{(\ell+1)} = \omega_3^{(\ell)}, \ldots,
\omega_{k-2}^{(\ell+1)} = \omega_{k-1}^{(\ell+1)}
\label{eq:constraint}
\end{equation}
is such that
\begin{equation}
\lim_{\ell\rightarrow \infty}C_{\theta}^{(\ell)}\ldots C_{\theta}^{(0)}=O.
\label{12}
\end{equation}
Moreover finding the supremum $\Theta$ such that \eqref{12} holds for any $\theta \le \Theta$
would provide a sharp limit for uniform boundedness of the solutions of the variable stepsize
BDF formula with the assigned values for stepsize ratios.

In view of the constraint \eqref{eq:constraint} this leads to a Markovian joint spectral radius problem.

\subsection*{The $3$--step BDF formula}

This case has been extensively studied in the literature (see e.g. \cite{GZ}) but it gives a natural benchmark
for our approach.
In this case we make use of the similarity transform determined by the matrix
$$
T=\left( \begin{array}{ccc}
1 & 1 & 1 \\
0 & 1 & 1 \\
0 & 0 & 1
\end{array} \right)
$$
which gives
$$
C_j=\left( \begin{array}{ccc}
\gamma_{j,1} & \gamma_{j,0} \\
1 & 0
\end{array} \right),
$$
$$
\gamma_{j,0}=-\alpha_{j,0} \ \ \ {\rm and} \ \ \ \gamma_{j,1}=-1+\alpha_{j,2}.
$$

Straightforward calculations give:
\begin{equation}
C(\omega_1,\omega_2)=\left( \begin{array}{ccc}
\gamma_1(\omega_1,\omega_2) & \gamma_0(\omega_1,\omega_2) \\
1 & 0
\end{array} \right),
\label{9}
\end{equation}
where
\begin{equation}
\gamma_0(\omega_1,\omega_2)=-\frac
{\omega_1^3\omega_2^2(\omega_2+1)^2}
{(\omega_1+1)(3\omega_1\omega_2^2+4\omega_1\omega_2+\omega_1+2\omega_2+1)}
\label{11}
\end{equation}
and
\begin{equation}
\gamma_1(\omega_1,\omega_2)=\frac
{\omega_2^2(\omega_1^2\omega_2^2+4\omega_1^2\omega_2+2\omega_1\omega_2+
3\omega_1^2+3\omega_1+1)}
{(\omega_1+1)(3\omega_1\omega_2^2+4\omega_1\omega_2+\omega_1+2\omega_2+1)}.
\label{10}
\end{equation}
The constraint is now
\[
\omega_1^{(j+1)} = \omega_2^{(j)}.
\]
In the considered illustrative case \eqref{eq:ex33}, we have to determine the constrained JSR of the family of nine matrices

\begin{equation}\label{eq:C9}
\{ C_i \}_{i=1,\ldots,9}
\end{equation}

where
\[
\begin{array}{lll}
C_1 = C(1/\theta,1/\theta) & C_2 = C(1/\theta,1) & C_3 = C(1/\theta,\theta) \\
C_4 = C(1,1/\theta) & C_5 = C(1,1) & C_6 = C(1,\theta) \\
C_7 = C(\theta,1/\theta) & C_8 = C(\theta,1) & C_9 = C(\theta,\theta)
\end{array}
\]

Based on the constraint \eqref{eq:ex33} we build the graph of admissible products which is depicted in Figure \ref{fig:GraphBDF3}.

\tikzset{%
   peer/.style={draw,circle,black,bottom color=white, top color= white, text=black, minimum width=25pt},
   superpeer/.style={draw, circle,  left color=white, text=black, minimum width=25pt},
   point/.style = {fill=black,inner sep=1pt, circle, minimum width=5pt,align=right,rotate=60},
   }

\begin{figure}[ht]
\begin{center}
\begin{tikzpicture}[->,>=stealth',shorten >=1pt,auto,node distance=2cm,
                    thick,main node/.style={circle,draw,font=\sffamily\bfseries}]

  \node[main node] (E) {$C_1$};
  \node[main node] (C) [right of=E] {$C_{2}$};
  \node[main node] (B) [right of=C] {$C_3$};
  \node[main node] (A)  [below of=E] {$C_4$};
  \node[main node] (E1) [below of=C] {$C_5$};
  \node[main node] (D1) [below of=B] {$C_6$};
  \node[main node] (C1) [below of=A] {$C_7$};
  \node[main node] (B1) [below of=E1] {$C_8$};
  \node[main node] (A1) [below of=D1] {$C_9$};

\path (E) edge [loop left] node[] {}(E);
\path (A) edge [color=black, bend right, above, sloped] node[] {}(E);
\path (C1) edge [color=black, bend left=20, above, sloped] node[] {}(E);

\path (E) edge [color=black, bend left=20, above, sloped] node[] { }(C);
\path (A) edge [color=black, bend left=20, above, sloped] node[] { }(C);
\path (C1) edge [color=black, bend left=20, above, sloped] node[] { }(C);

\path (E) edge [color=black, bend left=20, above, sloped] node[] { }(B);
\path (A) edge [color=black, bend left=20, above, sloped] node[] { }(B);
\path (C1) edge [color=black, bend left=20, above, sloped] node[] { }(B);

\path (C) edge [color=black, bend left=20, above, sloped] node[] { }(A);
\path (E1) edge [color=black, bend left=20, above, sloped] node[] { }(A);
\path (B1) edge [color=black, bend left=20, above, sloped] node[] { }(A);

\path (C) edge [color=black, bend left=20, above, sloped] node[] { }(E1);
\path (E1) edge [loop left] node[] { }(E1);
\path (B1) edge [color=black, bend left=20, above, sloped] node[] { }(E1);

\path (C) edge [color=black, bend left=20, above, sloped] node[] { }(D1);
\path (E1) edge [color=black, bend left=20, above, sloped] node[] { }(D1);
\path (B1) edge [color=black, bend left=20, above, sloped] node[] { }(D1);

\path (B) edge [color=black, bend left=20, above, sloped] node[] { }(C1);
\path (D1) edge [color=black, bend left=20, above, sloped] node[] { }(C1);
\path (A1) edge [color=black, bend left=20, above, sloped] node[] { }(C1);

\path (B) edge [color=black, bend left=20, above, sloped] node[] { }(B1);
\path (D1) edge [color=black, bend left=20, above, sloped] node[] { }(B1);
\path (A1) edge [color=black, bend right, above, sloped] node[] { }(B1);

\path (B) edge [color=black, bend left=20, above, sloped] node[] { }(A1);
\path (D1) edge [color=black, bend right, above, sloped] node[] { }(A1);
\path (A1) edge [loop left] node[] { }(A1);

\end{tikzpicture}
\caption{Graph $G_{3}$ of admissible products corresponding to the $3$--step BDF formula when only three values for the ratios $\omega_{j,t}$ are allowed}\label{fig:GraphBDF3}
\end{center}
\end{figure}
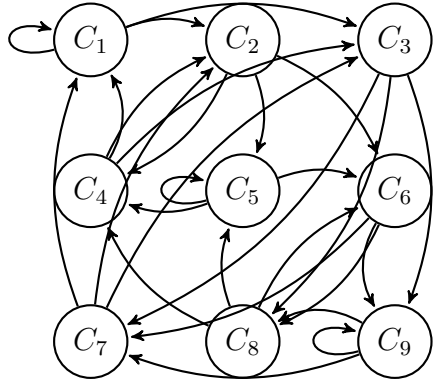

\subsection*{Numerical results}

 We apply Algorithm \ref{algo} to the $3$--step BDF formula when $\omega_{j,t}$ may assume only three values: $\theta$, $1$, and $1/\theta$.
 In particular we deal with the triplet $\xi = (G_{3}, \cL, \mathcal{C})$ where $G_{3}$ is the graph of admissible products plotted in Figure \ref{fig:GraphBDF3}, $\cL=\{L_i\}_{i=1}^9$ is a set of linear spaces in $\co^{2 \times 2}$, and $\mathcal{C}$ is the set of the nine matrices $C_i$ described in \eqref{eq:C9}.

 To run the calculation we use as candidate s.m.p. the matrix $C_9 = C(\theta,\theta)$ whose spectral radius is given by the modulus of a couple of complex conjugate eigenvalues for values of $\theta$ in the interval $[1,2]$.

 Following what described in Section \ref{ss-ident}, it is possible to reduce the number of vertices of the graph  $G_{3}$ to three.

 Furthermore, if we set $\theta=\Theta=\frac{1+\sqrt{5}}{2}$, then $C_9 = C(\theta,\theta)$ has spectral radius equal to 1 and, using the proposed algorithm, after three steps we construct an invariant complex polytope multinorm which contains two vertices in each linear space 
$L_i$.
For values $1 < \theta < \Theta$ we prove that $C_9$ is still an s.m.p. and the joint spectral radius of the family $\{ C_i \}_{i=1,\ldots,9}$ 
is smaller than $1$.

\subsection*{The $4$--step BDF formula}

This case is unexplored in the literature so that our results give an indication
about zero stability of such formula, which might be useful in a code implementing
it.

For this formula a suitable transformation is determined by the matrix
$$
T=\left(
\begin{array}{llll}
 1 & 1 & 1 & 1 \\
 0 & 1 & 1 & 1 \\
 0 & 0 & 1 & 1 \\
 0 & 0 & 0 & 1
\end{array} \right),
$$
which gives
\[
C(\omega_1,\omega_2,\omega_3) = \left( \begin{array}{ccc}
\gamma_2(\omega_1,\omega_2,\omega_3) & \gamma_1(\omega_1,\omega_2,\omega_3) & \gamma_0(\omega_1,\omega_2,\omega_3) \\
1 & 0 & 0 \\
0 & 1 & 0
\end{array}
\right)
\]
\vskip 5mm
\noindent The coefficients $\gamma_2, \gamma_1, \gamma_0$ are the following: \\ \scriptsize
\begin{eqnarray*}
\hskip -1.3cm
&& \gamma_2(\omega_1,\omega_2,\omega_3) =
\\
\hskip -1.3cm
&&\frac{(\omega_3+1)^2 (\omega_3 \omega_2+\omega_2+1)^2 (\omega_1 (\omega_3 \omega_2+\omega_2+1)+1)^2}{(\omega_2+1) (\omega_2
   \omega_1+\omega_1+1) \left(3 \omega_2 \omega_3^2+4 \omega_2 \omega_3+2 \omega_3+\omega_2+\omega_1 \left(\omega_2^2 (4 \omega_3+1)
   (\omega_3+1)^2+2 \omega_3+2 \omega_2 \left(3 \omega_3^2+4 \omega_3+1\right)+1\right)+1\right)}-1
\\[20mm]
\hskip -1.3cm
&& \gamma_1(\omega_1,\omega_2,\omega_3) =
\\
\hskip -1.3cm
&&-\frac{\omega_3^2 (\omega_3 \omega_2+\omega_2+1)^2 (\omega_1 (\omega_3 \omega_2+\omega_2+1)+1)^2}{(\omega_1+1) \left(3 \omega_2
   \omega_3^2+4 \omega_2 \omega_3+2 \omega_3+\omega_2+\omega_1 \left(\omega_2^2 (4 \omega_3+1) (\omega_3+1)^2+2 \omega_3+2 \omega_2 \left(3
   \omega_3^2+4 \omega_3+1\right)+1\right)+1\right)}+
\\[2mm]	
\hskip -1.3cm
&& \frac{(\omega_3+1)^2 (\omega_3 \omega_2+\omega_2+1)^2 (\omega_1 (\omega_3
   \omega_2+\omega_2+1)+1)^2}{(\omega_2+1) (\omega_2 \omega_1+\omega_1+1) \left(3 \omega_2 \omega_3^2+4 \omega_2 \omega_3+2
   \omega_3+\omega_2+\omega_1 \left(\omega_2^2 (4 \omega_3+1) (\omega_3+1)^2+2 \omega_3+2 \omega_2 \left(3 \omega_3^2+4
   \omega_3+1\right)+1\right)+1\right)}-1
\\[20mm]
\hskip -1.3cm
&& \gamma_0(\omega_1,\omega_2,\omega_3) =
\\
\hskip -1.3cm
&&\frac{\omega_1^4 \omega_2^3 \omega_3^2 (\omega_3+1)^2 (\omega_3 \omega_2+\omega_2+1)^2}{(\omega_1+1) (\omega_2 \omega_1+\omega_1+1) \left(3
   \omega_2 \omega_3^2+4 \omega_2 \omega_3+2 \omega_3+\omega_2+\omega_1 \left(\omega_2^2 (4 \omega_3+1) (\omega_3+1)^2+2 \omega_3+2
   \omega_2 \left(3 \omega_3^2+4 \omega_3+1\right)+1\right)+1\right)}
\end{eqnarray*} \normalsize

The constraint are now
\begin{equation}\label{eq:Const_BDF4}
\omega_1^{(j+1)} = \omega_2^{(j)}, \qquad \omega_2^{(j+1)} = \omega_3^{(j)}
\end{equation}

In the illustrative case \eqref{eq:ex33} we have to determine the constrained JSR of the family
$$
\{ C_i \}_{i=1,\ldots,27}
$$
where

\begin{equation}\label{eq:Matrices_BDF4}
\begin{array}{lll}
C_1 = C(1/\theta,1/\theta,1/\theta) & C_2 = C(1/\theta,1/\theta,1) & C_3 = C(1/\theta,1/\theta,\theta) \\
C_4 = C(1/\theta,1,1/\theta) & C_5 = C(1/\theta,1,1) & C_6 = C(1/\theta,1,\theta) \\
\ldots & \ldots & \ldots
\end{array}
\end{equation}

\subsection*{Numerical results}

As for the case of the $3$--step BDF, we consider first the illustrative case \eqref{eq:ex33} where $\omega_{j,t}$ is allowed to attain only three values: $\theta$, $1$, $1/\theta$, for some $\theta\geq 1$.

To run the calculation we use as candidate s.m.p. the matrix $C_{27} = C(\theta,\theta,\theta)$ whose spectral radius is given by the modulus of a pair of complex conjugate eigenvalues for values of $\theta$ in the interval $[1,2]$.

Similarly to the case studied for the $3$--step BDF formula, we can construct the graph of admissible products of matrices $\{C_i\}_{i=1}^{27}$ following the constraint \eqref{eq:Const_BDF4}. Furthermore, also in this case, we can reduce the number of vertices of such graph to four, applying the procedure described in Section \ref{ss-ident}.

If we set $\theta=\Theta\approx 1.2807$, such that $C_{27}$ has spectral radius equal to 1, then after seven steps of Algorithm \ref{algo} we construct an invariant complex polytope multinorm which contains six vertices in each of the four linear spaces.

We tested also the case of five possible values of $\omega_{j,t}$: $\theta$, $\sqrt{\theta}$, $1$, $1/\theta$, $1/\sqrt{\theta}$. The set $\{C_i\}$ contains now 125 matrices and, based on Section \ref{ss-ident}, the graph of admissible products can be reduced to a graph containing only 25 vertices.
Using as candidate s.m.p. the matrix corresponding to the case of $\omega$'s all equal to $\theta$, which is $C_{125}=C(\theta,\theta,\theta)$. If we set $\theta\approx 1.2807$ we construct an invariant complex polytope multinorm which contains eight vertices, except two of them that contain nine vertices.
The same procedure works for $1 < \theta < \Theta$.

Therefore on the basis of the numerical experiments we conjecture that the $4$--step BDF formula is zero stable
for $\omega_{j,t} \le 1.2807\ldots$, i.e. the value associated to a maximal constant increase of the stepsize
ratio.

\bigskip

\bigskip

\textbf{Acknowledgements.} A large part of this work was carried out when the third author visited
University of L'Aquila and Gran Sasso Science Institute (L'Aquila, Italy).  He is grateful for their hospitality.

\bigskip


\medskip


\begin{thebibliography}{}



\bibitem{AS}
T.~Ando and M.-H.~Shih,
\newblock {\em Simultaneous contractibility},
\newblock SIAM J. Matrix Anal. Appl. 19, (1998), No 2, 487--498.
\smallskip

\bibitem{AP}
D.\,Arapura and C.\,Peterson,
\newblock  {\em The common invariant
subspace problem: an approach via Gr\"obner bases},
\newblock Linear Alg. Appl., 384 (2004), 1–-7.
\smallskip

\bibitem{B1}
N.\,E.~Barabanov,
\newblock {\em Lyapunov indicator for discrete inclusions, I--III},
\newblock Autom. Remote Control, 49 (1988), No 2, 152--157.
\smallskip

\bibitem{BW}
M.~A.\,Berger and Y.\,Wang,
\newblock{\em Bounded semigroups of matrices},
\newblock Linear Alg. Appl.,
166 (1992) 21-27.
\smallskip





\bibitem{BS08}
F.\,Blanchini,  C.\,Savorgnanb,
\newblock {\em Stabilizability of switched linear systems does not imply the existence of convex Lyapunov functions},
\newblock Automatica, 44 (2008), no 4, 1166--1170.
\smallskip

\bibitem{BTV}
V.D.\, Blondel, J.\,Theys,\, and A.A.\,Vladimirov,
\newblock {\em An elementary counterexample to the finiteness conjecture},
\newblock SIAM Journal on Matrix Analysis, 24 (2003), no 4, 963-–970.
\smallskip

\bibitem{B}
R.\,Bowen,
\newblock {\em Markov partitions for Axiom A diffeomorphisms},
\newblock  American J. Math., 92 (1970), no 3, 725--747.
\smallskip

\bibitem{CG10}
A. Cicone, N. Guglielmi, S. Serra-Capizzano, and M. Zennaro,
\newblock {\em  Finiteness property of pairs of 2 x 2 sign--matrices via real extremal polytope norms},
\newblock Linear Alg. Appl., 432 (2010),   796--816. doi: 10.1016/j.laa.2009.09.022
\smallskip

\bibitem{CP15}
A.\,Cicone and V.Yu.\,Protasov,
\newblock {\em Fast computation of tight bounds for the joint spectral radius.},
\newblock preprint.
\smallskip

\bibitem{CMS}
Y.\,Chitour, P.\,Mason, and M.\,Sigalotti,
\newblock {\em On the marginal instability of linear switched systems},
\newblock Syst. Cont. Letters, 61 (2012), 747--757
\smallskip

\bibitem{D}
X.~Dai,
\newblock {\em Robust periodic stability implies uniform exponential
stability of Markovian jump linear systems and random
linear ordinary differential equations},
\newblock J. Franklin Inst., 351 (2014), 2910–-2937.
\smallskip

\bibitem{DL92}
I.~Daubechies and J.~Lagarias,
\newblock {\em Two-scale difference equations.
II. Local regularity, infinite products of matrices and fractals},
\newblock SIAM J. Math. Anal. 23 (1992), 1031--1079.
\smallskip

\bibitem{D95}
M.~Dogruel, U.~Ozguner,
\newblock {\em Stability of a set of matrices: a control theoretic approach},
\newblock Proceedings of the 34th IEEE Conference on Decision and Control,
13-15 Dec. 1995, vol. 2, 1324 - 1329.
\smallskip

\bibitem{EB}
D.B.\,Ellis and  M.G.\,Branton,
\newblock {\em Non-self-similar attractors of hyperbolic iterated function
systems},
\newblock Lecture Notes in Mathematics 1342 (Springer, Berlin, 1988),  158--171.
\smallskip

\bibitem{E95}
L.~Elsner,
\newblock {\em The generalized spectral-radius theorem: an analytic-geometric proof.}
\newblock Linear Alg. Appl. 220 (1995), 151--159.
\smallskip

\bibitem{Fek}
M.~Fekete,
 \newblock {\em \"Uber die Verteilung der Wurzeln bei gewissen algebraischen Gleichungen mit ganzzahligen Koeffizienten},
 \newblock  Mathematische Zeitschrift 17 (1923), no 1,  228–-249.
 \smallskip

\bibitem{FV12}
E.\,Fornasini, M.E.\,Valcher,
\newblock {\em Stability and stabilizability criteria for discrete-time positive switched systems}
\newblock IEEE Trans. Automat. Control  57  (2012), no 5, 1208--1221.
\smallskip

\bibitem{GVL13}
G.\,Golub and C.\,Van Loan.
\newblock {\em Matrix computations},
{\em Johns Hopkins Studies in the Mathematical Sciences}.
\newblock Johns Hopkins University Press, Baltimore, MD, 2013.
\smallskip

\bibitem{G}
G.~Gripenberg,
\newblock {\em  Computing the joint spectral radius},
\newblock Linear Alg. Appl., 234 (1996),   43--60.
\smallskip

\bibitem{GP13}
N.~Guglielmi and V.Yu.~Protasov,
\newblock {\em Exact computation of joint spectral characteristics of matrices},
\newblock Found. Comput. Math.,
13 (2013), no 1,  37--97.
\smallskip

\bibitem{GP16}
N.~Guglielmi and V.Yu.~Protasov,
\newblock {\em Invariant polytopes of linear operators with applications to regularity of wavelets and of subdivisions},
\newblock  SIAM J. Matrix Anal., 37 (2016), no 1,   18--52.
\smallskip

\bibitem{GWZ05}
N.~Guglielmi, F.~Wirth, and M.~Zennaro,
\newblock {\em Complex polytope extremality results for families of matrices},
\newblock SIAM J. Matrix Anal. Appl. 27 (2005),  721--743.
\smallskip

\bibitem{GZ15}
N.\,Guglielmi and M.\,Zennaro,
\newblock {\em Canonical construction for Barabanov polytope norms and antinorms
for sets of matrices},
\newblock SIAM J. Matrix Anal. Appl. 36 (2015),  634--655.
\smallskip



\bibitem{GZ01}
N.~Guglielmi and M.~Zennaro.
\newblock {\em On the asymptotic properties of a family of matrices}.
\newblock { Linear Alg. Appl.},  322 (2001), 169--192.

\bibitem{GZ07}
N.~Guglielmi and M.~Zennaro,
\newblock {\em Balanced complex polytopes and related vector and matrix norms,}
\newblock J. Convex Anal. 14 (2007), 729--766.
\smallskip

\bibitem{GZ}
N. Guglielmi, and M. Zennaro,
\newblock {\em  An algorithm for finding extremal polytope norms of matrix families.},
\newblock Linear Alg. Appl., 428 (2008),   2265--2282. doi: 10.1016/j.laa.2007.07.009
\smallskip

\bibitem{HNW}
E. Hairer, S. P. Norsett, and G. Wanner
\newblock {\em  Solving Ordinary Differential Equations I},
\newblock Springer--Verlag, Berlin Heidelberg, 1993.
\smallskip

\bibitem{HLZ}
R.\,Hildebrand, L.\,Lokutsievskiy, and M.\,Zelikin,
\newblock {\em Typicalness of chaotic fractal behaviour of integral vortexes in Hamiltonian systems with discontinuous right hand side},
\newblock 	arXiv:1506.02320
\smallskip

\bibitem{H}
J.\,E.~Hutchinson,
\newblock {\em Fractals and self-similarity},
\newblock Indiana Univ. Math.~J., 30 (1981), no 5, 713--747.
\smallskip

\bibitem{J09}
R.\,M.~Jungers,
\newblock {\em The joint spectral radius: theory and applications, }
\newblock Lecture Notes in Control and Information Sciences, vol. 385, Springer-Verlag, Berlin Heidelberg, 2009.
\smallskip


\bibitem{KH}
A.\,Katok and B.\,Hasselblatt,
\newblock {\em Introduction to the Modern Theory of Dynamical
Systems}, \\
\newblock Cambridge University Press, Cambridge, 1996.
\smallskip

\bibitem{K}
V.\,Kozyakin,
\newblock {\em The Berger-Wang formula for the Markovian joint spectral radius},
\newblock Linear Alg. Appl., 448 (2014), 315--328.
\smallskip

\bibitem{LA09}
H.\,Lin and P.J.\,Antsaklis,
\newblock {\em Stability and stabilizability of switched linear
systems: a survey of recent results},
\newblock IEEE Trans. Autom. Contr., 54 (2009), no 2, 308--322.
\smallskip

\bibitem{L}
L.V.\,Lokutsievskii,
\newblock {\em Fractal structure of hyperbolic  Lipschitzian dynamical systems},
\newblock Russian J. Math. Physics, 19 (2012), no 1, 27--44.
\smallskip

\bibitem{MR14}
C.~M\"oller and U.~Reif,
\newblock {\em A tree-based approach to joint spectral radius determination},
\newblock Linear Alg. Appl., 563 (2014), 154--170.
\smallskip

\bibitem{OPJ}
M.\,Ogura, V.M.\,Preciado, and R.M.\,Jungers,
\newblock {\em Efficient method for computing lower bounds on the $p$-radius of switched linear systems},
\newblock arXiv preprint (2015),  arXiv:1503.03034
\smallskip


\bibitem{PJ15}
P.\,A.~Parrilo and    A.\,Jadbabaie,
\newblock {\em Approximation of the joint spectral radius using sum of squares},
\newblock  Linear Alg. Appl.  428 (2008), no 10, 2385--2402.
\smallskip

\bibitem{PhJ1}
M.~Philippe and R.\,M.~Jungers,
\newblock {\em Converse Lyapunov theorems for discrete-time linear
switching systems with regular switching sequences},
\newblock 	arXiv:1410.7197
\smallskip

\bibitem{PhJ2}
M.\,Philippe and  R.M.\,Jungers,
\newblock {\em A sufficient condition for the boundedness of matrix products accepted by an automaton},
\newblock Proceedings of the 18th International Conference on Hybrid Systems: Computation and Control, ACM New York, NY, USA (2015), 51--57.
\smallskip

\bibitem{P96}
V.\,Yu.~Protasov,
\newblock {\em The joint spectral radius and invariant sets of linear operators.}
\newblock Fundam. Prikl. Mat. 2 (1996),  205--231.
\smallskip

\bibitem{P97}
V.\,Yu.~Protasov,
\newblock {\em The generalized spectral radius. A geometric approach},
\newblock Izvestiya Math. 61 (1997), 995--1030.
\smallskip

\bibitem{P06}
V.\,Yu.~Protasov,
\newblock {\em Fractal curves and wavelets},
\newblock Izvestiya Math., 70 (2006), No 5, 123--162.
\smallskip

\bibitem{PG}
V.\,Yu.~Protasov and N.\,Guglielmi,
\newblock {\em Matrix approach to the global and local
regularity of wavelets},
\newblock Poincare J. Anal. Appl., (2015), no 2, 77--92.
\smallskip

\bibitem{PJ}
V.\,Yu.~Protasov, R.\,M.\,Jungers,
\newblock {\em Resonance and marginal instability of switching systems}, \\
\newblock Nonlinear Analysis: Hybrid Systems,  17 (2015), 81--93.
\smallskip


\bibitem{PJB}
V.\,Yu.~Protasov, R.\,M.~Jungers, and V.\,D.~Blondel,
\newblock {\em Joint spectral characteristics of matrices:
a conic programming approach},
\newblock SIAM J. Matrix Anal. Appl., 31 (2010), no 4,  2146--2162.
\smallskip

\bibitem{RS60}
G.\,C.~Rota and G.~Strang,
\newblock {\em A note on the joint spectral radius},
\newblock Kon. Nederl. Acad. Wet. Proc. 63 (1960),  379--381.
\smallskip

\bibitem{SDP08}
E.\,De Santis, M.D.\,Di Benedetto, G.\,Pola,
\newblock {\em Stabilizability of linear switching systems},
\newblock Nonlinear Analysis: Hybrid Systems, 2 (2008), no 3, 750--764.
\smallskip

\bibitem{SFS}
M.\,Souza, A.R.\,Fioravanti, and R.N.\,Shorten,
\newblock {\em Dwell-time control of continuous-time switched linear systems},
\newblock Proceedings of the IEEE Conference on Decision and Control (2015), 4661--4666.
\smallskip

\bibitem{SS90}
G.W.~Stewart, J.G.~Sun,
\newblock {\em Matrix perturbation theory},
\newblock Academic Press, New York, 1990.
\smallskip

\bibitem{WRDV}
Y.Wang, N.\,Roohi, G.E.\,Dullerud, and M.\,Viswanathan,
\newblock {\em Stability of linear autonomous systems under regular switching sequences},
\newblock Proceedings of the IEEE Conference on Decision and Control (2015), 5445--5450
\smallskip



\end{thebibliography}
\end{document}